\documentclass[11pt,a4paper,reqno]{amsart}

\usepackage[colorlinks,bookmarksopen,
bookmarksnumbered,citecolor=red,urlcolor=red]{hyperref}
\usepackage{amsaddr}
\usepackage{epsfig,amsmath,amssymb,amsfonts,mathrsfs}
\usepackage{graphicx}
\usepackage{float} 
\usepackage{textcomp, eufrak}
\usepackage{verbatim}
\usepackage{fancyhdr}
\usepackage{enumitem}
\usepackage{soul}

\numberwithin{equation}{section}

\newtheorem{theorem}{Theorem}[section]
\newtheorem{definition}[theorem]{Definition}
\newtheorem{lem}[theorem]{Lemma}
\newtheorem{prop}[theorem]{Proposition}
\newtheorem{rem}[theorem]{Remark}

\newtheorem{exa}[theorem]{Example}

\numberwithin{equation}{section}
\newtheorem{Assum}[theorem]{Assumption}

\newcommand{\R}{\mathbb{R}}

\newcommand{\N}{\mathbb{N}}

\newcommand{\Rp}{\mathbb{R}^+}
\newcommand{\F}{\mathcal{F}}

\newcommand{\PS}{(\Omega, \mathcal{F}, \mathbb{P})}

\newcommand{\I}{\mathbb{I}}

\newcommand{\Bb}{\mathcal{B}}

\newcommand{\D}{\mathcal{D}}
\newcommand{\Z}{\mathbb{Z}}

\newcommand{\T}{\mathbb{T}}

\newcommand{\p}{\mathbb{P}}

\newcommand{\bn}{\begin{definition}}
\newcommand{\en}{\end{definition}} 
\newcommand{\bt}{\begin{theorem}}                
\newcommand{\et}{\end{theorem}}
 \newcommand{\bnm}{\begin{enumerate}}              
\newcommand{\enm}{\end{enumerate}}
\newcommand{\br}{\begin{rem}} 
\newcommand{\er}{\end{rem}}

\newcommand{\om}{\omega}

\newcommand{\Om}{\Omega}

\newcommand{\btm}{\begin{itemize}}
 \newcommand{\etm }{\end{itemize}}

\newcommand{\E}{\mathbb{E}}

\newcommand{\Rd}{\R^d}

\newcommand{\Ic}{\mathcal{I}}
\newcommand{\Sc}{\mathbb{S}^1}
\newcommand{\Scd}{\mathbb{S}^{d-1}}

\newcommand{\LG}{\mathcal{L}}

\renewcommand\leftmark{{Stochastic averaging for SDEs in the random periodic regime}}

\renewcommand\rightmark{{Stochastic averaging for SDEs in the random periodic regime }}
\begin{document}

\author[]{Kenneth Uda}
\address{\vspace*{-0.2cm} Department of Mathematics, University of Edinburgh, Scotland, UK}
\email{K.Uda@ed.ac.uk}
\title{ Averaging principle for stochastic differential equations in the random periodic regime}

\date{}
\maketitle

\begin{abstract}
We present the validity of stochastic averaging principle for non-autonomous slow-fast stochastic differential equations (SDEs) whose fast motions admit random periodic solutions. Our investigation is motivated by some problems arising from multi-scale stochastic dynamical systems, where configurations are time dependent due to nonlinearity of the underlying vector fields and the onset of time dependent random invariant sets. Averaging principle with respect to uniform ergodicity of the fast motion is no longer available in this scenario. Lyapunov second method together with synchronous coupling and strong Feller property of Markovian flows of SDEs are used to prove the ergodicity of time periodic measures of the fast motion on certain minimal Poincar\'e section and consequently identify the averaging limit. 
\end{abstract}
{\bf Keywords:} Periodic measures; PS-ergodicity; random periodic solutions;  stochastic averaging;  H\"omander's Lie bracket conditions; strong Feller property.

\newcommand{\Xt}{\mathbb{R}^d}
\newcommand{\Bt}{\mathcal{M}}
\newcommand{\Yt}{\R^N}
 \newcommand{\Nb}{\mathbb{K}}
\tableofcontents
 \section{Introduction, notations and setup}
 \subsection{Introduction}
 Many nonlinear deterministic and stochastic dynamical models in science and engineering consist of motions at disperate spatiotemporal scales, giving rise to slow-fast systems. Multi-scale methods may be employed to analyse some of the desirable structures of these slow-fast motions in the long run. Multi-scale methods are well-established theoretical and computational techniques in nonlinear stochastic dynamical systems and related areas, its applications range from signal processing, climate dynamics, material science, mathematical finance, molecular dynamics, etc.  Averaging principle is an important multi-scale technique, it suggests that the slow subsystem of a slow-fast system may be approximated by a simpler system obtained by averaging over the fast motion. The active usage of  various averaging could be traced back to 18th century as a well-known approximation procedure in celestial mechanics, for example, in the study of perturbation of planetary motion from the periodic variation in the orbital elements. Rigorous results justifying averaging principle go back to series of papers by Krylov, Bogolyubov and Mitropol'skii (e.g.,\cite{Kry-Bog}). Generalisation of averaging principle to chaotic and stochastic dynamical systems with important applications in science and engineering are among the most investigated problems in nonlinear dynamical systems and related areas (see e.g.,~\cite{Hasm66, Fredlinbook, Kifer09, Stuart08}).  In probabilistic language, averaging principle is a variant of weak law of large numbers.
 
 The primary objective of the present paper is to investigate averaging principle in the random dynamical systems (RDS) framework (e.g.,~\cite{Arnold01, Kifer09, Kifer04, Kifer03, Kifer95, Stuart08, YU99}) for fully coupled slow-fast SDEs with time dependent vector fields. Our investigation stems from some problems arising from long time behaviour of stochastic dynamical systems, where dynamical configurations are time dependent due to the existence of time dependent random invariant sets. Models whose evolution rule are time dependent slow-fast SDEs are particularly natural in numerous applications and are speedily becoming cradle for the development new mathematical ideas. Two applicable examples are the following:
 \btm[leftmargin = 0.7cm]
\item Climate-weather interactions when diurnal cycle and seasonal cycle are taken into consideration and one has to encode the time dependent statistics of the fast variables (the weather) when making predictions (e.g.,~\cite{Bran12, Checkron11, Andy10, Andy16}).  \item Neuronal learning models, where the strength of neurons connections evolve at a slower rate to account for synaptic events (synaptic plasticity), leading to a time dependent slow-fast SDE (e.g,~\cite{Wain12}), etc.
\etm
The study of these class of systems are important and challenging problems in modern applied mathematics. One of the technical challenges is to identify the right notion of convergence of the slow motion to the averaged motion. There are no many theoretical results in this direction, however, significant progress have been made over the last decade in the  study of the long time behaviour of time dependent RDS (e.g., \cite{Feng11, Feng12, Feng18, Uda14, Uda16, Wan151,Wan152, Wan14, Zhao09}). Based on these results, one can justify averaging principle for time dependent RDS once a suitable notion of ergodicity is known.
 
The purpose of ergodicity is to study invariant measures and related problems, it is one of the well studied problems in dynamical systems, stochastic analysis, statistical physics and related areas. Intuitively, ergodic dynamical system is the one that  the behaviour of observables averaged over long time is the same as the behaviour of observables averaged over phase space. In RDS, ergodicity is the cornerstone in the study of long time behaviour (e.g., part II and III of \cite{Arnold} and references therein) and in particular, in the investigation of averaging. Various caveats of ergodicity such as uniform/unique or non-uniform ergodicity exist in the literature, depending on the mode of convergence and or existence of unique or several extremal invariant measure(s). Averaging principle for SDEs when the fast subsystem is independent of the slow variable holds true whenever the fast motion is uniform or unique ergodic (e.g., \cite{Hasm66, Stuart08, Fredlinbook}). In general, uniform or unique ergodicity is not available in the fully coupled case and in particular, time dependent setting, so one has to resort to other form of ergodicity to identify the averaging limit. In the time dependent framework, some form of ergodicity like ergodicity on Poincar\'e sections (PS-ergodicity \cite{Feng18, Uda18}) or semi-uniform ergodicity (e.g,~\cite{Sturm00}) of the fast subsystem could be employed to identify the averaging limit. In this paper, we restrict ourselves to when the fast subsystem is SDE with time periodic forcing, then under some weak dissipative conditions, PS-ergodicity is proved and consequently, used to identify the averaging limit. Extra essential difficulty in our case is that we avoid some analytic assumption such as uniform ellipticity of the fast subsystem, at the expense of proving the existence of unique stable random periodic process and non-degeneracy of Malliavin covariance of solution flows under certain weak dissipative assumptions, restricted H\"omander's Lie bracket condition and Lyapunov second method.
 
 Lyapunov second method is a powerful technique for the investigation of stability of solutions of nonlinear dynamical systems and control systems in finite and infinite dimensions. Its extension to random dynamical systems generated by SDEs is essentially due to Hasminskii (e.g.,~\cite{Hasm80}). This method was explored further by X. Mao (e.g.,~\cite{Mao94}) and extended further by Schmalfuss in (\cite{Schmal01}) to the case of nontrivial stationary solutions and random attractors.  The method involves the study of random invariant sets without explicit knowledge of solutions of the system. The framework of Lyapunov second method is simple in the sense that its verification is based explicitly on the underlying vector fields without explicit representation of flows of solutions. Here, we employ Lyapunov second method together with  synchronous coupling to prove the existence of unique stable random periodic solution and ergodicity of the corresponding periodic measures. These techniques are valid even when the drift of the SDE is only one-sided Lipschitz continuous. Importantly, it makes it easier to prove results in the narrow topology generated by the dual of the Markov evolutionary semigroup of the flows of SDEs.

The rest of the paper is organised as follows. In \S\ref{Nota}, we fix some revelant function spaces used in the paper and in \S\ref{Ave set-Up}, we formulate and describe the problem considered in the paper. The existence of random periodic solutions and ergodic periodic measures for fast subsystem are considered in \S\ref{Ran_p}. The validity of the averaging principle in the random periodic regime is proved in \S\ref{Aver_pp}, and we conclude with a toy example arising from kinetic theories of turbulent flows.
  \subsection{Notations}\label{Nota}
 Let $(\mathbb{M}, \text{d})$ be a complete separable metric space, we denote
 the set of bounded continuous real-valued functions by
\begin{align*} 
 \mathcal{C}_b(\mathbb{M}) =\Big\{f: \mathbb{M}\rightarrow\R \; \text{continuous }\;: \Vert f\Vert_{\infty}<\infty\Big\},\qquad
\Vert f\Vert_{\infty} := \sup_{x\in \mathbb{M}}\vert f(x)\vert.
\end{align*}
The set of bounded Lipschitz continuous real-valued functions is denoted by
 \begin{align*}
  \text{Lip}_b(\mathbb{M})&=\{f\in\mathcal{C}_b(\mathbb{M}): \Vert f\Vert_{BL}<\infty\},\\
  \Vert f\Vert_{BL} &= \max\{\Vert f\Vert_{\infty},\Vert f\Vert_L\},\quad \text{and} \quad \Vert f\Vert_{L} = \sup\bigg\{\frac{\vert f(y)-f(z)\vert}{\text{d}(y,z)}: y\neq z\bigg\}.
 \end{align*}
 
  Let $(\mathbb{M}, \Bb(\mathbb{M}), \mu)$ be a Borel measure space, we denote $L^p(\mathbb{M}), \; 1\leqslant p<\infty$ as the set of real-valued Lebesgue integrable functions $f: \mathbb{M}\rightarrow \R$ defined by 
\begin{align*}
L^p(\mathbb{M}):=\Big\{f: \mathbb{M}\rightarrow\R \; \text{measurable}\; : \Vert f\Vert_{p}<\infty\Big\}, \qquad \Vert f\Vert_{p}=\left(\int_{\mathbb{M}}\vert f\vert^pd\mu\right)^{1/p}.
\end{align*}
 Let $\Bt\subseteq \mathbb{M}$ be a non-empty Borel measurable set, we shall most at times restrict the set of functions above on $\Bt$, for example, $f\in \mathcal{C}_b(\Bt)$ means $f\!\!\upharpoonright_{\Bt}\ \in \mathcal{C}_b(\mathbb{M})$, where  $ f\!\!\upharpoonright_{\Bt}: \Bt\cap \text{Dom}f\rightarrow\R.$
\bn[\cite{Arnoldp, Kunita}]\;\rm
Let $l\in \N$ and $0\leq \delta \leq 1.$
\btm
\item[(a)] Let $n, m\in \N\setminus\{0\}$, denote $\mathcal{C}^{l,\delta}(\R^n;\R^m)$ as the Fr\'echet space of functions $f: \R^n\rightarrow\R^m$ which are continuously differentiable and $l$-th derivative is $\delta$-H\"older continuous, with the semi norms 
\begin{align*}
&\Vert f\Vert_{l,0,\Nb}:= \sum_{0\leq \vert \alpha\vert \leq l}\sup_{x\in \Nb} \vert D^{\alpha}f(x)\vert,\\
& \Vert f\Vert_{l,\delta;\Nb}: = \Vert f\Vert _{l,0;\Nb}+\sum_{\vert \alpha\vert =l}\sup_{x,y\in \Nb, x\neq y}\frac{\vert D^{\alpha}f(x)-D^{\alpha}f(y)\vert}{\vert x-y\vert^{\delta}}, \qquad 0<\delta \leq 1.
\end{align*} 
Here $\Nb$ is a compact convex subset of $\R^n,$ $\alpha = (\alpha_1,\cdots, \alpha_n)\in \N^n$ and $\vert \alpha\vert = \alpha_1+\cdots+\alpha_n$ and 
\begin{align*}
D^{\alpha}f(x) = \frac{\partial^{\vert \alpha\vert }f}{(\partial x_1)^{\alpha_1}\cdots (\partial x_n)^{\alpha_n}}.
\end{align*} 
\item[(b)] Let $\mathcal{C}_b^{l,\delta}(\R^n;\R^m)$ be the Banach space of the functions $f: \R^n\rightarrow \R^m$ which are in the Fr\'echet space $\mathcal{C}^{l,\delta}(\R^n;\R^e)$ with the norm 
\begin{align*}
&\Vert f\Vert_{l,0}:= \sup_{x\in\R^n}\frac{\vert f(x)\vert}{1+\vert x\vert}+\sum_{1\leq \vert \alpha\vert \leq l}\sup_{x\in \R^n}\vert D^{\alpha}f(x)\vert,\\
& \Vert f\Vert_{l,\delta}:= \Vert f\Vert_{l,0}+\sum_{\vert \alpha\vert =l}\sup_{x,y,\in \R^n, x\neq y}\frac{\vert D^{\alpha}f(x) -D^{\alpha}f(y)\vert}{\vert x-y\vert^{\delta}}<\infty, \quad 0<\delta\leq 1.
\end{align*}
We endow the space $\mathcal{C}_b^{l}(\R^n;\R^{m})$ with the norm 
\begin{align*}
\Vert f\Vert_{l} = \sup_{x\in \R^n}\frac{|f(x)|}{1+|x|}+\sum_{1\leqslant |\alpha|\leqslant l}\sup_{x\in \R^n}|D^\alpha f(x)|.
\end{align*}
\etm
\en
 
Let $\PS$ be a complete probability space and $\mathcal{G}\subseteq\F,$ we denote $L^p(\Om, \mathcal{G}, \p),\;p\geqslant 1$ as the space of $\mathcal{G}$-measurable random variables $X:\Om\rightarrow\R^m$ such that $\E\vert X\vert^p<\infty $ equiped with the $L^p$ norm 
$\Vert X\Vert_p = \left(\E\vert X\vert^p\right)^{1/p}.$\\
We shall fix the probability space $\PS$ as the classical Wiener space, i.e., $\Om=\mathcal{C}_0(\R;\R^m), \; \break m\in \N$, is a linear subspace of continuous functions that take zero at $t=0$, endowed with compact open topology defined via
\begin{align*}
\text{d}(\om,\hat{\om})= \sum_{n=0}^\infty\frac{1}{2^n}\frac{\Vert \om -\hat{\om}\Vert_n}{1+\Vert \om -\hat{\om}\Vert_n}, \qquad \Vert \om -\hat{\om}\Vert_{n} = \sup_{t\in [-n, n]}\vert \om(t) -\hat{\om}(t)\vert.
\end{align*}
The sigma algebra $\F$ is the Borel sigma algebra generated by open subsets of $\Om$ and $\p$ is the Wiener measure, i.e., that the law of the process $\om\in\Om$ with $\om(0) =0.$
\subsection{The averaging setup}\label{Ave set-Up}
 We are concerned with the following slow-fast SDE
\begin{align}\label{Slow fast}
\begin{cases}
dX^{\varepsilon}_t = \varepsilon F(X^{\varepsilon}_t, Y_t^{\varepsilon})dt, \quad & X_0^\varepsilon = x\\
dY^{\varepsilon}_t = b(t,X^{\varepsilon}_t, Y^{\varepsilon}_t)dt+ \sum_{k=1}^{N}\sigma_k(t,X^{\varepsilon}_t, Y_t^{\varepsilon})dW^k_t, \quad &Y_0^\varepsilon = y,
\end{cases}
\end{align}
where $0<\varepsilon\ll 1$ is the time scales separation between $X^\varepsilon$ and $Y^\varepsilon$. The vector fields $F: \Xt\times\Yt\rightarrow \Xt,$  $b: \R\times \Xt\times\Yt\rightarrow \Yt$ and $\sigma_k: \R\times\Xt\times \Yt\rightarrow  \Yt, \; 1\leqslant k\leqslant N,$ are such that  $F\in \mathcal{C}^{l-1}_b(\Rd\times\Yt; \Rd)$ and $b, \sigma_k\in \mathcal{C}^{l, \delta}_b(\R\times\Rd\times\Yt; \Yt), \; l\geqslant 2$.
Given this regularity of the vector fields, the solutions of the SDE (\ref{Slow fast}) exist and generate stochastic flow of $\mathcal{C}^{l-1,\alpha}$ diffeomorphisms with $0<\alpha\leqslant \delta$ (e.g.,~\cite{Elw78, Ikeda81, Hkunita, Kunita}). We define the flow map  $\Xi_t^{\varepsilon}(\om)\equiv \Xi_{0,t}^{\varepsilon}(\om): \Rd\times\Yt\rightarrow\Rd\times\Yt$ by $\Xi_{0,t}^{\varepsilon,x,y}(\om) =(X_{0,t}^{\varepsilon,x,y}(\om), Y_{0,t}^{\varepsilon, x,y}(\om))\equiv (X_{t}^{\varepsilon,x,y}(\om), Y_{t}^{\varepsilon, x,y}(\om)),$ for almost all $\om\in \Om,$ where $\Om= = \mathcal{C}_0(\R;\Yt)$ is $N$-dimensional Wiener space equiped with compact open topology.

If $\varepsilon =0,$ the flow $\Xi^{0}_t(\om)$ reduces to $\Xi^{0,x,y}_t(\om) = (x, Y_t^{x,y}(\om))$ where $Y_t^{x,y}$ is the solution of the following fast subsystem with the slow variable frozen 
\begin{align}\label{Fast}
dY^{x,y}_t = b(t, x, Y^{x,y}_t)dt + \sum_{k=1}^{N}\sigma_k(t, x, Y^{x,y}_t)dW_t, \quad Y_0^{x,y} = y.
\end{align}
It is natural to think of the slow-fast SDE (\ref{Slow fast}) as a perturbation of the parametrised SDE (\ref{Fast}). In this way, the flow map $\Xi^{0,x,y}_t =(x, Y^{x,y}_t)$ describes idealised physical system where $x= (x_1,\cdots, x_d)$ are fixed parameters, whereas the perturbed flow $\Xi^{\varepsilon,x,y}_t$ is describing real system with the evolution of the parameters taken into consideration. 

In some previous works (e.g., \cite{Kifer01, Kifer04b, Kifer09, Kifer04, Kifer95, YU99}), the validity of averaging principle for fully coupled ODEs and SDEs were investigated. The averaging principle for the fully coupled system is significantly different from the one when the fast motion is independent of the slow variables. The required assumptions in the fully coupled case depend on the specific class of systems and convergence of the slow motion to averaged motion is in a different sense (e.g.,~\cite{Fredlinbook, Kifer09, Loch88}). Consider the limit $\bar{F}$ defined by
\begin{align}\label{Unierg}
\bar{F}(x)=\bar{F}_y(x) := \lim_{T\rightarrow\infty}\frac{1}{T}\int_0^T F(x, Y_t^{x,y})dt.
\end{align}
For example, if $\mu^x$ is an ergodic invariant measure of the fast motion $Y^x=\{Y_t^{x,y}: \;(t, y)\in\Rp\times \Yt\}$, then the limit (\ref{Unierg}) exists for $\mu^x$-almost all $y,$ due to ergodicity and it coincides with the following integral
\begin{align}
\bar{F}(x)=\bar{F}_\mu(x) := \int_{\Yt}F(x, y)\mu^x(dy).
\end{align}
In the case, where the vector fields $ b, \sigma_k, \; 1\leqslant k\leqslant N$ are independent of $(t,x),$ the fast flow defines a closed system, i.e., $Y_t^x = Y_t $ and $\mu^x = \mu,$ in this case, one recovers the averaging principle proved long ago by Hasminskii (\cite{Hasm66}). So the averaged vector field $\bar{F} $ inherit the regularities of slow vector field $F,$  thus, there exists a unique solution $\bar{X}^{\varepsilon, x}_t$ to the following averaged equation 
\begin{align}\label{Av1}
\frac{d\bar{X}^{\varepsilon,x}_t}{dt} = \varepsilon\bar{F}(\bar{X}^{\varepsilon,x}_t), \quad \bar{X}^{\varepsilon, x}_0 = x.
\end{align}
In this case, the averaging principle states (e.g.,~\cite{Hasm66, Kifer03}), for $\delta>0,$
\begin{align}\label{Hasm_Av}
\lim_{\varepsilon\rightarrow 0}\p\left(\sup_{0\leqslant t\leqslant T/\varepsilon}\vert X_t^{\varepsilon, x,y}-\bar{X}_t^{\varepsilon,x}\vert>\delta\right) =0, \quad \text{for all}\; (x,y)\in \Lambda,
\end{align}
where $\Lambda :=\left\{(x,y)\in\Rd\times \Yt: \bar{F}_y(x) = \bar{F}_\mu(x)\right\}$. \\ If $b, \sigma_k, \; 1\leqslant k\leqslant N,$ is autonomous but depend on $(x,y),$ the averaged vector field $\bar{F}(x)$ need not be continuous in $x,$ let alone Lipschitz continuous; so the averaged equation (\ref{Av1}) may admit infinitely many solutions or no solution at all. This is partly because in this case, the fast motions $Y^x=\{Y_t^{x,y}: (t, y)\in\Rp\times \Yt\}$ do not define a closed system, as $Y_t^x$ may have different properties for different $x;$ hence, no well-defined family of invariant measures $\{\mu^x: x\in \Rd\}.$ Even in the rare case when $\mu^x$ is the same for all $x\in \Rd,$ averaging in fully coupled system may not satisfy the limit (\ref{Hasm_Av}), for example in the presence of resonance (e.g.,~\cite{Kifer09, Loch88}). It is proved in (e.g.,\cite{Arnold01, Kifer01}) that if the convergence in (\ref{Unierg}) is uniform in $x$ and $y,$ that any limit points $\bar{X}_t^x$ in probability of $X^{\varepsilon, x,y}_{t/\varepsilon}$ as $\varepsilon\rightarrow 0,$ solves the averaged equation 
\begin{align}
\frac{d\bar{X}_t^x}{dt} = \bar{F}(\bar{X}_t^x), \quad \bar{X}_0^x = x.
\end{align}
Note that the limit (\ref{Unierg}) is uniform with respect to $y$ if and only if the fast flow $Y^x$ is uniquely ergodic, this occurs rather rarely and exclude some important models. In general, averaging in the fully coupled case requires different setup with stronger and specific assumptions (c.f.~\cite{Kifer09}). Precisely, it is noted in (\cite{ Kifer09, Kifer03, Kifer95}) that in the fully coupled case, it is useful to have necessary and sufficient conditions which ensure the following:\btm[leftmargin=0.7cm] \item[{}]  Given that  $\mathbb{Q}\in \mathcal{P}(\Rd\times\Yt)$ admit the disintegration $\mathbb{Q}(dx,dy) = \mu^x(dy)\nu(dx),$ with $\nu\in \mathcal{P}(\Rd)$ and $\mu^x\in \mathcal{P}(\Yt),$  then for $\Nb\Subset$\footnote{$\Nb\Subset\mathbb{M}$ means that $\Nb$ is a non-empty relative compact subset of $\mathbb{M}$.}$\Rd$ and $C\Subset\Yt$ and any  $\delta>0,$ 
\begin{align}\label{Ad_AVE}
\lim_{\varepsilon\rightarrow 0}\int_{\Nb}\int_{C}\p\left(\sup_{0\leqslant t\leqslant T/\varepsilon}\vert X_t^{\varepsilon,x,y}-\bar{X}_t^{\varepsilon,x}\vert>\delta\right)\mu^x(dy)\nu(dx)=0.
\end{align}
\etm
The averaging principle in the form (\ref{Ad_AVE}) was proved long ago by Anosov in the case when the fast subsystem satisfy Liouville property (see chapter 2 of~\cite{Loch88}). It was later generalised by Kifer to the case when the fast subsystem is {\it axiom A} dynamical systems and to the case of Markovian SDEs (\cite{ Kifer09, Kifer03, Kifer95} and references therein). 
 
In the present paper, the fast motion is time dependent, so ergodic invariant measures may not exist. However, if  for example, the coefficients of the slow-fast SDE (\ref{Slow fast}) are time periodic, quasi-periodic, uniform almost periodic or uniform almost automorphic, one may consider evolution system of probability measures $\{\mu^x_t: t\in \R\}$ (e.g., \cite{Schu11}) generated by the Markov evolution  $\{\mathcal{P}_{s,t}^x: s\leqslant t\}$ of the fast motion, i.e., 
\begin{align}\label{Evolu}
\int_{\Yt}\mathcal{P}_{s,t}^xf(y)\mu_s^x(dy) = \int_{\Yt}f(y)\mu_t^x(dy), \quad s\leqslant t, \quad f\in \mathcal{C}_b(\Yt).
\end{align}
One can exploit the nice property of the above evolution system of probability measures with some dissipative conditions on the coefficients of the fast subsystem to identify the averaging limit. There are existing results in this time dependent case; for example, when the fast subsystem is time periodic and uniform elliptic was considered in (e.g.,~\cite{YU92, Wain13}) and applied to investigate the long time behaviour of neuronal learning models in \cite{Wain12}. The case of uniform almost periodic fast subsystem was recently considered for stochastic reaction-diffusion equation in \cite{Cera17}. 
In this paper, we shall assume that the coefficients of the fast motion is periodic in time, i.e., exists $\tau>0$ such that $b(t+\tau, x,y) = b(t,x,y)$, $\sigma(t+\tau,x,y) = \sigma(t,x,y)$. Under this time periodicity and regularity of the vector fields $b ,\sigma_k, \; 1\leqslant k\leqslant N,$ it have been rigorously established in (e.g., \cite{Feng11, Feng12, Feng18, Uda14, Uda16, Uda18, Zhao09}) that the SDE (\ref{Fast}) admit random periodic solution $S^x: \R\times\Om\rightarrow\Yt$ (see~\S\ref{S_Ranp}). The law of the random periodic solution $\mu_t^x(A) := \p\{\om: S^x(t,\om)\in A\}, \; A\in \Bb(\Yt)$ generate evolution system of probability measures and, in particular, periodic measures (see \S\ref{Ranp_M}). Under some dissipative conditions on the fast SDE (\ref{Fast}), we establish ergodicity of the periodic measures $\mu_t^x, t\in [0, \tau)$ on certain minimal Poincar\'e section, we refer to this notion as {\it PS-ergodicity} of periodic measures as recently introduced in \cite{Feng18} (see~\S\ref{Ranp_M}).

\smallskip
 A useful way to visualise PS-ergodic periodic measures is to lift the stochastic flow to the cylinder $\Sc\times\Yt$ \footnote{i.e., $[0, \tau]\cong\Sc$, with the homeomorphism $t\mapsto e^{2\pi it/\tau}$.} and obtain an autonomous stochastic flow of the form $\tilde{Y}_t^{\varepsilon,x,\tilde{y}} = (t+ s\mod \tau, Y_{s,t+s}^{\varepsilon, x, y}(\theta_{-s}\om))$ with the initial value $\tilde{y} = (s, y)\in \Sc\times \Yt.$ The random periodic solution on $\Sc\times\Yt$ then takes the form $\tilde{S}^x(t,\om) = (s \mod \tau, S^x(t, \om)).$ In this case, the slow-fast SDE on the extended phase space $\Rd\times\Sc\times\Yt,$ reads
 \begin{align}\label{lifed}
 \begin{cases}
 dX^{\varepsilon,x,y}_t = \varepsilon F(X^{\varepsilon,x,y}_t, Y_t^{\varepsilon,x,y})dt,\\
 d\begin{pmatrix} s_t \\ Y_t^{\varepsilon, x,y}\end{pmatrix} = \begin{pmatrix}
  1\\ b(s_t, x,Y_t^{\varepsilon, x,y})
  \end{pmatrix}dt + \begin{pmatrix}
 0& 0\\ 0 &\sum_{k=1}^N\sigma_k(s_t,x, Y_t^{\varepsilon, x,y})
  \end{pmatrix} \begin{pmatrix} dW_t^0\\ dW^k_t \end{pmatrix} \\
  X_0^{\varepsilon,x,y} = x\in \Rd, \quad (s_0, Y_0^{\varepsilon,x,y}) = (s,y)=:\tilde{y}\in \Sc\times\Yt,
 \end{cases} 
 \end{align}
where  $W^0_t$ is a one dimensional Brownian motion independent of the $N$-dimensional Brownian motion $\{W_t^k: 1\leqslant k\leqslant N, t\geqslant 0\}$. Apparently, this lifting does not reduce the complexity of the problem, but  the SDE is now autonomous. In this way, under the regularity of the coefficients, the lifted slow-fast SDE (\ref{lifed}) generate a Markovian RDS on the extended phase space $\Rd\times\Sc\times\Yt$ (see \cite{Arnoldp}, chapter 2 of \cite{Arnold} and references therein). 

 In addition to the PS-ergodicity of the periodic measures, it is important in this fully coupled case to prove that the branches of the periodic measures $x\mapsto \tilde{\mu}_r^x, \; r\in [0, \tau)$ are Lipschitz continuous w.r.t. $x\in \Nb\Subset\Rd$ in the narrow topology of $\mathcal{P}(\Bt),$ where $\Bt\subset \Sc\times\Yt$ is the union of minimal Poincar\'e sections (see \S\ref{Ran_Lp}). This Lipschitz property of periodic measures together with the regularity of the slow vector field $F,$ would guarantee the existence of  $\hat{C}_F>0$ such that
\begin{align}\label{Lbdness}
\begin{cases}
\vert \bar{F}(x)-\bar{F}(z)\vert\leqslant \hat{C}_F\vert x-z\vert,\\
\vert \bar{F}(x)\vert\leqslant \hat{C}_F,
\end{cases}
\end{align}
where 
\begin{align}
\bar{F}(x):=\frac{1}{\tau}\int_0^{\tau}\int_{\Bt}F(x,\pi_2(\tilde{y}))\tilde{\mu}^x_r(\tilde{y})dr
\end{align}
and $\pi_2:\Sc\times\Yt\rightarrow\Yt$ with $\pi_2(\tilde{y}) = y$. The local Lipschitz and boundedness property (\ref{Lbdness}) would yield the existence of a unique solution $\bar{X}_t^{\varepsilon,x}$ to the averaged equation
\begin{align}\label{AVper}
\frac{d\bar{X}_t^{\varepsilon,x}}{dt} = \varepsilon \bar{F}(\bar{X}_t^{\varepsilon,x}),\quad \bar{X}_0^{\varepsilon,x} =x.
\end{align}
Finally, given the PS-ergodicity and Lipschitz dependence of the family $\{\tilde{\mu}^x_r: r\in [0, \tau),\; x\in \Nb\Subset \Rd\}$, we are able to deduce averaging principle in the form of the limit (\ref{Ad_AVE}) with the set $C$ replaced by minimal Poincar\'e section $\Bt\subset\Sc\times\Yt$ (see \S\ref{AV_limit}). 

\medskip

\section{Random periodicity}\label{Ran_p}
In this section, we investigate the existence of unique stable random periodic solution and PS-ergodicity of periodic measures for Markovian RDS generated by some class of SDEs in finite dimensions. To facilitate our presentation, we first recall definitions of stochastic flows and RDS (see, e.g.,~\cite{Arnold,Arnoldp, Hkunita, Kunita}) and random periodic solutions (see, e.g.,~\cite{Feng11, Feng12, Feng18, Zhao09}).
\subsection{Random periodic solutions}
\begin{definition}[Stochastic flow \cite{Hkunita, Kunita}]\rm\;  Let $\mathbb{M}$ be a smooth manifold and let $\{\Xi^z_{t,s}(\om): s,t\in \T\subseteq\R,\; z \in \mathbb{M}\}$ be a random field defined on a complete probability space $(\Om, \F, \p).$ 
The map $\Xi_{t,s}(\om)$ is called a {\it stochastic flow of homeomorphism} if there exists a null set $\mathcal{N}\subset \Om$ such that for any $\om\notin \mathcal{N},$ the family of continuous maps $\{\Xi_{t,s}(\om): s,t\in\T\}$ is a flow of homeomorphism in the following sense:
\btm
\item[(i)] $\Xi_{t,s}(\om) = \Xi_{t,u}(\om)\circ\Xi_{u,s}(\om)$ holds for any $s,t,u.$
\item[(ii)] $\Xi_{s,s}(\om) = \text{Id}_{\mathbb{M}}$, for all $s.$
\item[(iii)] the map $\Xi_{t,s}(\om): \mathbb{M}\rightarrow \mathbb{M}$ is a homeomorphism for any $(s,t)$.
\etm
The map $\Xi_{t,s}(\om)$ is a {\it stochastic flow of $\mathcal{C}^k$-diffeormorphism}, if it is a homeomorphism and $\Xi^z_{t,s}(\om)$ is $k$-times continuously differentiable with respect to $z$ for all $s,t\in \T\subseteq\R$ and the derivatives are continuous in $(s,t,z).$ 
\end{definition}

Let $\F_s^t$ be the smallest sub $\sigma$-field of $\F$ containing $\cap_{\varepsilon>0}\sigma\big(\Xi_{u,v}: s-\varepsilon\leq u, v\leq t+\varepsilon\big) $ and all null sets of $\F$. Then the two parameter filtration $\{\F_s^t: s\leq t\}$  is the filtration generated by the stochastic flow $\{\Xi_{s,t}: s,t\in \T\subseteq\R\}.$
\bn[RDS \cite{Arnold, Arnoldp}]\rm
Let $\mathbb{M}$ be a measurable space and  $\theta$ be an ergodic flow of measurable transformations of a complete probability space $(\Om,\F,\p),$ a {\it random dynamical system} (RDS) over $\theta$ is a mapping $\Phi: \T\times\Om\times\mathbb{M}\rightarrow \mathbb{M}$ such that 
\btm
\item[(a)] $(t,\om,z)\mapsto \Phi(t,\om,z)$ is measurable,
\item[(b)] $\Phi(0,\om) = \text{Id}_{\mathbb{M}}$ for all $\om\in\Om$ and 
\item[(c)] $\Phi(t+s, \om)= \Phi(t,\theta_s\om)\circ\Phi(t,\om),\quad \text{for all}\quad s,t\in\T,$ $\om\in\Om$ \quad (cocycle property),
\item[(d)] if $\mathbb{M}$ is a topological space, then $\Phi$ is continuous if $(t,z)\mapsto\Phi(t,\om,z)$ is continuous,
\item[(e)] if $\mathbb{M}$ is a $\mathcal{C}^k$ manifold $1\leq k\leq \infty,$ then  $\Phi$ is smooth of class $\mathcal{C}^k,$ if $\Phi(t,\om,z)$ $k$-times differentiable with respect to $z,$ and the derivatives are continuous with respect to $(t,z).$  
\etm
\en

\medskip
Let  $\theta: \R\times\Om\rightarrow\Om,$ be a $(\Bb(\R)\otimes\F, \F)$-measurable flow defined by
\begin{align}\label{Mes_fl}
\theta_t\om(\cdot) = \om(t+\cdot)- \om(t).
\end{align}
It is well known that the measurable flow $(\theta_t)_{t\in \R}$ is ergodic and defines a filtered metric dynamical system $\theta = (\Om, \F, \p, (\F^t_s)_{t\geqslant s},(\theta_t)_{t\in\R})$ (e.g., \cite{Arnold, Arnoldp}). Furthermore, under suitable regularity of the coefficients of autonomous SDEs together with appropriate adoption of two-sided stochastic calculus, it is well known that the solutions of autonomous SDEs on finite dimensions generate RDS over $\theta$ (e.g.,~\cite{Arnold, Arnoldp, Elw78, Ikeda81, Kunita}). 
\begin{definition}[Random periodic solution for stochastic flows \cite{Feng11, Feng12, Zhao09}]\label{Chu}\rm\;
 A random periodic solution of period $\tau>0$ of a stochastic flow $\Xi:\Delta\times\Om\times \mathbb{M}\rightarrow \mathbb{M},$ is an $\F$-measurable function $S:\T\times\Om\rightarrow\mathbb{M}$ such that $$ S(t+\tau,\om) = S(t,\theta_{\tau}\om)\quad \text{and} \quad \Xi_{s,t+s}^{S(s,\om)}(\om) = S(t+s,\om),\quad \p\;\text{-a.s.},$$ for any $(t,s) \in \Delta\subset\T\times\T$ and $\om\in\Om,$ where $\Delta:=\{(t,s)\in\T\times\T: s\leq t\}.$
\end{definition}
\begin{exa}\label{ex32}\rm
 Let $\mathbb{M} = \R,$ consider the following SDE 
\begin{equation}\label{ehnl}
 \begin{cases} dZ_t = -\alpha(t)Z_tdt+\beta(t)dt+ dW_t, \quad t\geqslant t_0\\ Z_{t_0}= z_0\in \R. \end{cases}
\end{equation}
 Assume that $\alpha, \beta: \R\rightarrow \R$ are continuous functions and there exists $\tau> 0 $ such that $\alpha(t+\tau) = \alpha(t)$ and $\beta(t+\tau) = \beta(t)$ with $$\int_{-\infty}^te^{-\int_s^t\alpha(u)du}\beta(s)ds+\int_{-\infty}^{t}e^{-2\int_{s}^{t}\alpha(u)du}ds<\infty, \quad \text{for}\quad 0\leq t\leq \tau.$$ 
 The stochastic flow $\{\Xi_{t_0,t}(\om): t\geqslant t_0\}$ generated by the SDE (\ref{ehnl}) is defined by $$\Xi_{t_0,t}^{z_0}= z_0e^{-\int_{t_0}^{t}\alpha(u)du} +\int_{t_0}^te^{-\int_s^t\alpha(u)du}\beta(s)ds+ \int_{t_0}^{t}e^{-\int_{s}^t\alpha(u)du}dW_s,$$   
 and the process $S(t,\om)$ defined by $$ S(t,\om) = \int_{-\infty}^te^{-\int_s^t\alpha(u)du}\beta(s)ds+\int_{-\infty}^{t}e^{-\int_{s}^{t}\alpha(u)du}dW_s(\om)$$ is a random periodic solution of period $\tau$.
\end{exa} 
\bn[Random periodic solution for RDS \cite{Feng18, Zhao09}]\rm
A random periodic path of period $\tau>0$ of an RDS $\Phi: \Rp\times\Om\times\mathbb{M}\rightarrow \mathbb{M}$ is an $\F$-measurable function $S: \R\times\Om\rightarrow \mathbb{M}$ such that for almost all $\om \in \Om,$
\begin{align}
S(t+\tau, \om) = S(t, \theta_\tau\om) \quad \text{and} \quad \Phi(t, \theta_s\om, S(s,\om)) = S(t+s, \om)
\end{align}
for $t\in\Rp, \; s\in \R$.
\en
\begin{exa}[\cite{Feng11, Feng12}]\rm\;
 Let $ b: \Yt\rightarrow\Yt, \; N\geqslant 2 $ be a regular vector field and let $\{\Phi(t,.): \; t\in\T\subseteq\R\}$ be a deterministic flow generated by the ODE
  \begin{align}\label{ODE1}
  \frac{dY_t}{dt} =b(Y_t).
   \end{align}
If $u: \T\rightarrow\Yt$ is a periodic solution of the ODE (\ref{ODE1}) of period $\tau>0,$ i.e., 
\begin{align*}
u(t+\tau) = u(t) \quad \text{and}\quad \Phi(t+s, u(s)) = u(t+s), \quad t,s\in\T\subseteq\R.
\end{align*}
 Consider the stochastic process $X_t(\om) = u(t)+Z_t(\om),$  where $Z_t$ solves the following time periodic forcing SDE
\begin{align}\label{sde111}
dZ_t   = \hat{b}(t, Z_t)dt+ \hat{\sigma}(t, Z_t)dW_t,
\end{align}
where $\hat{b}(t,z) := b(u(t)+z)-b(u(t))$ and $\hat{\sigma}(t,z) := \sigma(u(t)+z).$ If $v(t,\om)$ is a random periodic solution of the time periodic forcing SDE (\ref{sde111}), then $S(t,\om) = u(t)+v(t,\om)$ is a random periodic solution of the autonomous SDE
\begin{align}
dX_t = b(X_t)dt+\sigma(X_t)dW_t.
\end{align} 
\end{exa}

\subsection{Existence of stable random periodic solutions}\label{S_Ranp}
Here, we present a version of a random periodic results in (\cite{Uda14, Uda18}). This result will be used to prove PS-ergodicity and Lipschitz dependence of periodic measures on parameters.
As pointed in the introduction, the technique here is a variant of Hasminksii's method for stability of solutions of SDEs (see, e.g.,~\cite{Hasm80, Mao94, Schmal01}). The central idea is to investigate the infinitesimal separation of trajectories via the average growth of a Lyapunov function $V: \R\times\Yt\rightarrow\Rp$ along the two-point motions $\{Y_{s,t}^y(\om)-Y^z_{s,t}(\om): t\geqslant s, \; z,y\in \Yt,\; z\neq y\}$ which leads to the existence of (pre)-fixed point via Arzel\'a--Ascoli like argument.

\smallskip
First, we fix $s\in \R$ and recall that the transition probability function $P(s,y; t, .)$ induced by solutions of the SDE (\ref{Fast}) is defined by
 \begin{equation}
 P^x(s,y; t+s, A) = \p\big(\{\om\in\Om: Y^{x,y}_{s,t+s}(\om )\in A\}\big), \quad t\geqslant 0, \quad A\in\Bb(\Yt).
 \end{equation}
 The Markov evolution $\mathcal{P}^x_{s,t+s}: \mathcal{C}_b(\Yt)\rightarrow\mathcal{C}_b(\Yt)$ is defined by 
 \begin{equation}
 \mathcal{P}^x_{s,t+s}f(y) = \int_{\Yt}f(z)P^x(s,y; t+s, dz) = \E[f(Y_{s,t+s}(\om))| Y^{x,y}_{s,s} = y]
 \end{equation}
 and for any probability measure $\mu$ on $(\Yt,\Bb(\Yt))$, the Markov dual $\mathcal{P}_{t,t+s}^{x*}$ is defined 
 \begin{equation}
 (\mathcal{P}_{s,t+s}^{x*}\mu)(A) = \int_{\Yt}P^x(s,y; t+s, A)\mu(dy), \quad t\geqslant 0, \quad  A\in \Bb(\Yt).
 \end{equation}
Next, the average growth of the function $V$ along the trajectory $Y^{x,y}_{t,s}(\om)$ gives rise to the generator $\mathcal{L}_x$,
 \begin{equation}\label{one point}
 \mathcal{L}_xV(t,y) = \lim_{h\rightarrow 0}\frac{\E[V(t+h,Y_{t,t+h}(\om))|Y_{t,t}^{x,y}=y]-V(t, y)}{\vert h\vert}.
 \end{equation}
For $V\in C_b^{1,2}(\R\times\Yt),$ we use It\^o's formula to write $\mathcal{L}_x$ as 
$$\mathcal{L}_xV(t,y) =\frac{\partial V(t,y)}{\partial t}+ \sum_{i=1}^{N}b^i(t, x,y)\frac{\partial V(t,y)}{\partial y^i}+ \frac{1}{2}\sum_{i,j=1}^{N}\sum_{k=1}^N\sigma_{ik}(t,x,y)\sigma_{jk}(t,x,y)\frac{\partial^2V(t,y)}{\partial y^i\partial y^j}.$$  
Now, consider the difference between two solutions of the fast subsystem starting from two different initial values $\{Y^{x,y}_{s,t+s}(\om)-Y^{x,z}_{s,t+s}(\om): y,z\in \Yt,\; y\neq z,\; t\geqslant 0\}.$ Here, we think of $Y_{s,t}^{x,y}(\om) -Y_{s,t}^{x,z}(\om)$ as a special case of coupling two copies of It\^o process namely; $(y,z)\mapsto H(y,z) = y-z,$ so that $H(Y^{x,y}_{s,t}(\om), Y_{s,t}^{x,z}(\om)) = Y_{s,t}^{x,y}(\om)-Y^{x,z}_{s,t}(\om)$.
Then, by generalised It\^o's formula (Theorem 8.1 in \cite{Hkunita}), the average growth of the function $V\in \mathcal{C}_b^{1,2}(\R\times\Yt; \Rp)$ along the trajectory $ \{Y^{x,y}_{s,t+s}(\om)-Y^{x,z}_{s,t+s}(\om): y\neq z, t\geqslant 0\}$, yields the second order differential operator $\mathcal{L}_x^{(2)}$, defined by
\begin{align}\label{2pp}
\mathcal{L}_x^{(2)}V(t,y-z) &=\partial_tV(t,y-z)+ D_yV(t,y-z)\left(b(t,x,y)-b(t,x,z)\right)\nonumber \\[.1cm] &\hspace{.5cm}+\frac{1}{2}\text{trace}\left( [\sigma(t,x,y)-\sigma(t,x,z)]^TD^2_yV(t,y-z)[\sigma(t,x,y)-\sigma(t,x,z)]\right),
\end{align}
where 
$\sigma(t,y)= (\sigma_{ik}(t,y))_{1\leqslant i,k\leqslant N},\quad
D_yV= \Big(\frac{\partial V}{\partial y_i}\Big)_{1\leqslant i\leqslant N}, \quad D_y^2V= \Big(\frac{\partial ^2 V}{\partial y_i\partial y_j}\Big)_{1\leqslant i,j\leqslant N}.$

\smallskip
Finally, we return to the lifted fast subsystem on the cylinder $\Sc\times\Yt,$
  \begin{align}\label{Fast2}
  d\tilde Y_t^{x,\tilde{y}} = \tilde b(x,\tilde{Y}_t^{x,\tilde{y}})dt+ \tilde{\sigma}(x,\tilde{Y}_t^{x,\tilde{y}}) d\tilde{W}_t, \quad \tilde{Y}_0^{x,\tilde{y}} = \tilde{y}.  
  \end{align}
  where $\tilde{y} = (s,y)\in \Sc\times\Yt,\; \tilde{W}_t = (W_t^0, W_t^k)^T, \; 1\leqslant k\leqslant N,$ the vector fields $\tilde{b}, \; \tilde{\sigma}$ are given by $$\tilde{b}(x,(s,y)) = (1, b(s,x,y))^T \quad \text{and}\quad \tilde{\sigma}(x,(s,y))=\begin{pmatrix}
 0& 0\\ 0 &\sum_{k=1}^N\sigma_k(s,x, y)
  \end{pmatrix}. $$
 Throughout the remaining part of this section, we denote the stochastic flow generated by the lifted SDE (\ref{Fast2}) by  $$\Phi^x(t, \om, (s,y)) :=\tilde{Y}_t^{x,(s,y)}(\om) =\Big( t+s\mod\tau, Y^{x,y}_{s,t+s}(\theta_{-s}\om)\Big).$$ In some places, we may skip the dependence on $x,$ when it does not cause confusion.
 
  Next, recall that the transition probability function in space-time $\Rp\times\Yt$ takes the form 
\begin{align}\label{Space-time}
\notag \tilde{P}(0, (s,y); t, \mathbb{T}\times B) &:= \p\left(\om: (t+s, Y^y_{s,t+s}(\theta_{-s}\om))\in \mathbb{T}\times B\right) \\ &= P(s, x; t+s, B)\I_{\mathbb{T}}(t+s), \;\; \mathbb{T}\in \Bb(\Rp),\;B\in \Bb(\Yt).
\end{align}
 In particular, the transition probability function $\tilde{P}(0,\tilde{y},t,.), \; t\geqslant 0$ and the Markov evolution $(\tilde{\mathcal{P}}_t)_{t\geqslant 0}$ (now Markov semigroup) on $\Sc\times\Yt$ is defined for $\tilde{y}\in \Sc\times\Yt$ as  
\begin{align*}
\quad\tilde{P}(0, \tilde{y}; t, \tilde{A})& = \p(\om: \Phi(t,\om,\tilde{y})\in \tilde{A}),\quad \tilde{A}\in \Bb(\Sc)\otimes\Bb(\Yt),\\
\tilde{\mathcal{P}}_t\tilde h(\tilde{y}) &= \int_{\Sc\times\Yt}\tilde{h}(\tilde{z})\tilde{P}(0,\tilde{y}; t, d\tilde{z}), \quad \tilde{h}\in \mathcal{C}_b(\Sc\times\Yt).
\end{align*}
The Markov dual $(\tilde{\mathcal{P}}^*_{t})_{t\geqslant 0}$ on $\Sc\times\Yt$ is 
\begin{align*}
(\tilde{\mathcal{P}}^*_t\tilde{\mu})(\tilde{A}) = \int_{\Sc\times\Yt}\tilde{P}(0, \tilde{y}; t, \tilde{A})\tilde{\mu}(d\tilde{y}), \quad \tilde{\mu}\in \mathcal{P}(\Sc\times\Yt).
\end{align*}
With all the above notations in place, we present our result on the existence and uniqueness of stable random periodic solution on $\Sc\times\Yt$ under the following assumptions.
\begin{Assum}\label{A2.1}\rm\;
 In addition to the regularity and time periodicity of the coefficients of the fast subsystem, we suppose further that the following hold: \btm \item[{}] There exists a function $V\in \mathcal{C}^{1,2}(\R\times\Yt; \Rp)$, $\tau$-periodic in the time component such that $V(t,0) = 0$ and for all $ t\in\R,$  $y, z\in \Yt,$ $x\in \Nb\Subset\Rd,$
\begin{align}\label{Disp1}
\vert y\vert^p\leqslant V(t,y)\leqslant C\vert y\vert^p \quad \text{and}\quad \mathcal{L}_{x}^{(2)}V(t,y-z)\leqslant \lambda(t,x) V(t, y-z)
\end{align}
for $p\geqslant 1,$  $C\geqslant 1$ with 
\begin{align}\label{Disp11}
\limsup_{t\rightarrow \infty}\frac{1}{2t}\int_0^t\lambda(u)du:=\beta<0.
\end{align}
where $\lambda(t):=\sup_{x\in \bar{\Nb}}\lambda(t,x)$.\\ Moreover, let $\Phi^x(t,\om,(s,y)) = (t+s\mod\tau, Y_{s,t+s}^{x,y}(\theta_{-s}\om))$ be the lifted stochastic flow generated by the fast subsystem with the slow variable frozen such that 
\begin{align}\label{Temp1}
\sup_{t\in [0,1]}\log \E\left[V\left(\Phi^x(t,\om,(s,y))-(s,y)\right)\right]<\infty, 
\end{align}
\etm
\end{Assum}

\begin{theorem} [Existence of stable random periodic path \cite{Uda14, Uda18}]\label{R_per1}\rm\;
Let $\{Y_{s,t}^y( \om): t\geqslant s\}$ be a stochastic flow generated by an SDE with the coefficients $b(t,.) ,\sigma_k(t,.)\in \mathcal{C}^{l,\delta}_b(\Yt;\Yt), \; 1\leqslant k\leqslant N,\; l\geqslant 2,$ $\tau$-periodic and continuously differentiable in time component satisfying assumptions \ref{A2.1}.
Then, there exist a unique $\F_{-\infty}^{s}:= \bigvee_{u\leqslant s}\F_u^s$-measurable and exponential stable random periodic process $\tilde{S}: \R\times\Om\rightarrow\Sc\times\Yt$.
 \end{theorem}
 \noindent {\it Proof.}
 Our presentation is in two steps, in the first step, we show that  $\{\Phi(., \om, \Phi(k\tau,\theta_{-k\tau}\om,\tilde{y})\}_{k\in \N}$ is a Cauchy sequence in $\mathcal{C}([0, \tau]; \Sc\times\Yt)$ and converges exponentially fast, the  uniqueness of the limit follows by cocycle property of $\Phi$ on $\Sc\times\Yt$. The second step is to show via the cocylce property again and time periodicity of the coefficients of our SDE that the limit from the first step satisfies the definition of random periodic solution.\\
 {\bf Step I:} 
 We investigate the separation of the trajectories $\Phi(t,\om,(s,y))$ and $\Phi(t,\om, \Phi(1, \theta_{-1}\om, (s,y))).$ Let $(s,y)\neq (u,z)$ for some $(u,z)\in \Sc\times\Rd$, following same construction as in \cite{Schmal01}, we set
 \begin{align*}
(s_1, X_1(\om)) = \tilde{X}_1(\om) :=\begin{cases} (s,x),\; &\text{if} \quad (s,y)\neq (1+s\mod\tau, X(s+1,s,\theta_{-s-1}\om, y))\\
(u,z), &\text{if} \quad (s,y) = (1+s\mod\tau, X(s+1,s,\theta_{-s-1}\om, y)).
 \end{cases}
 \end{align*}
Notice that $\tilde{X}_1$ is $\F_{-\infty}^{s}$-measurable and $V(\Phi(1,\theta_{-1}\om,(s,x) - \tilde{X}_1(s,\om))$ satisfies the condition (\ref{Temp1}), since the random variable $V(\Phi(1, \theta_{-1}\om,(s,y) ) -(s,y))$ satisfies the same condition (\ref{Temp1}). Denote $\tilde{Y}(t,\om) := \Phi(t,\om, \tilde{X}_1(\om))$ and $\tilde{Z}(t,\om) := \Phi(t, \om, \Phi(1, \theta_{-1}\om, \tilde{y})),$ for $t\geqslant 0.$
 We apply It\^o's formula (e.g.,~\cite{Hkunita, Schmal01}) on the logarithmic process $\log V(\tilde{Y}(t,\om) -\tilde{Z}(t,\om))$, to obtain
 \begin{align*}
 \log V(\tilde{Y}(t,\om) -\tilde{Z}(t,\om)) & = \log V(\tilde{X}_1(\om) - \Phi(1,\theta_{-1}\om, \tilde{y})) \\ & \hspace{1cm}+ \int_0^t\frac{\tilde{\LG}^{(2)}V(\tilde{Y}(r,\om) -\tilde{Z}(r,\om))}{V(\tilde{Y}(r,\om) -\tilde{Z}(r,\om))}dr+ N(t,\om) -\frac{1}{2}q(t,\om),
\end{align*} 
\newcommand{\HC}{\mathcal{H}}
where \begin{align*}
q(t,\om)&=\int_{0}^t\big(\HC V(\tilde{Y}(r,\om) -\tilde{Z}(r,\om)))\big)^2dr,\quad
N(t,\om)=\int_{0}^t\big(\HC V(\tilde{Y}(r,\om) -\tilde{Z}(r,\om))\big)dW_r,\\
\HC V(\tilde{x},\tilde{y})&= \frac{D_yV(\tilde{x}-\tilde{y})\big(\bar{\sigma}(\tilde{x})-\tilde{\sigma}(\tilde{y})\big)}{V(\tilde{x}-\tilde{y})}, \quad \text{and}\quad \tilde{\LG}^{(2)} = \LG^{(2)}+\frac{\partial}{\partial s}.
\end{align*} 
Using  a variant of Doob's martingale inequality (c.f.~Lemma 6.2 of~\cite{Mao94}, see also,~\cite{Schmal01}), we have 
 \begin{align*}
 \p\bigg\{\om: \sup_{t\in [0, k]}\Big(N(t,\om) - \frac{1}{2}q(t,\om)\Big)>2 \log k \bigg\}\leqslant \frac{1}{k^2}.
\end{align*}  
Moreover, by the $\p$-preserving property of $\theta,$ we obtain
\begin{align*}
 \p\bigg\{\om: \sup_{t\in [0, k]}\Big(N(t,\theta_{-k+1}\om) - \frac{1}{2}q(t,\theta_{-k+1}\om)\Big)>2 \log k \bigg\}\leqslant \frac{1}{k^2}.
\end{align*}  
Next, we use conditions (\ref{Disp1}) and Borel Cantelli lemma, we have for $k>1,$
\begin{align}\label{Exp_fast}
\nonumber \frac{1}{k-1}\sup_{t\in [k-1, k]}\log\vert \tilde{Y}(t,\om) -\tilde{Z}(t,\om))\vert &\leqslant 
\frac{1}{p(k-1)}\sup_{t\in [k-1, k]}\log V(\tilde{Y}(t,\om) -\tilde{Z}(t,\om)) \\   \nonumber
 &\leqslant \frac{1}{p(k-1)}\sup_{t\in [k-1, k]} \log V(\tilde{X}_1(\om) - \Phi(1,\theta_{-1}\om, \tilde{y}))\\ &  \qquad+\frac{1}{p(k-1)}\sup_{t\in [k-1, k]} \int_0^t\lambda(r)dr+ \frac{2 \log k}{p(k-1)}.
\end{align}
Inequality (\ref{Exp_fast}) and conditions (\ref{Disp1}) together with the cocycle property of $\Phi$, lead to the following estimate
\begin{align}\label{CAUCHY}
&\frac{1}{(k-1)}\sup_{t\in [k-1, k]}\log \vert  \Phi(t,\theta_{-k+1}\om, \tilde{X}_1(\theta_{-k+1}\om)) -  \Phi(t, \theta_{-k+1}\om, \Phi(1, \theta_{-k}\om, \tilde{y}) \vert \\
&= \frac{1}{(k-1)}\sup_{t\in [0, 1]}\log \vert  \Phi(t,\om,\Phi(k-1,\theta_{-k+1} ,\bar{X}_1(\theta_{-k+1}\om)) -  \Phi(t, \om, \Phi(k,\theta_{-k}\om,\tilde{y}) \vert \nonumber \\
&\leqslant  \frac{1}{p(k-1)}\log V(\tilde{X}_1(\theta_{-k+1}\om) - \Phi(1,\theta_{-k+1}\om, \tilde{y}))+ \frac{1}{p(k-1)}\sup_{t\in [k-1, k]} \int_0^t\lambda(r)dr+ \frac{2 \log k}{p(k-1) \nonumber}\\
&\leqslant  \frac{1}{p(k-1)}\log V(\tilde{X}_1(\theta_{-k+1}\om) - \Phi(1,\theta_{-k+1}\om, \tilde{x}))+ \frac{\beta(k-1)}{pk}+ \frac{2 \log k}{p(k-1)}, 
\end{align}
As $\log V(\Phi(1,\theta_{-1}\om,\tilde{y}) - \tilde{X}_1(\om))$ satisfies condition (\ref{Temp1}), there exist $\hat{\Om}$ with $\p(\hat{\Om}) = 1$, a random variable $\gamma(\om)$ with $\E[\gamma(\om)]<\infty$ such that for any $\om\in \hat{\Om}$ and $0<\varepsilon<-\frac{\beta}{2p}$, there exists  $k(\om)\in \N$ such that for $k\geqslant k(\om),$ 
\begin{align}
\notag & \sup_{t\in [0, 1]}\vert  \Phi(t,\om,\Phi(k-1,\theta_{-k+1} ,\tilde{X}_1(\theta_{-k+1}\om)) -  \Phi(t, \om, \Phi(k,\theta_{-k}\om,\tilde{y}) \vert\\ &\hspace{3cm}\leqslant\gamma(\om)\exp\Big(  \big(\frac{1}{2p}\beta+\varepsilon\big)k\Big), \quad \p\; \text{-a.s.}
\end{align}
Finally, the cocycle property of $\Phi$, leads to 
\begin{align}\label{Expfas}
\notag &\sup_{t\in [0, \tau]}\vert  \Phi(t,\om,\Phi(k\tau-\tau,\theta_{-k\tau+\tau}\om,\tilde{X}_1(\theta_{-k\tau+\tau}\om)) -  \Phi(t, \om, \Phi(k\tau,\theta_{-k\tau}\om,\tilde{y}) \vert&  \\ &\hspace{3cm} \leqslant\hat{\gamma}(\om)\exp\Big(  \big(\frac{1}{2p}\beta+\varepsilon\big)k\tau\Big),\quad \p\;\text{-a.s.}
\end{align}
This implies that $\{\Phi(., \om, \Phi(k\tau,\theta_{-k\tau}\om,\tilde{y})\}_{k\in \N}$ is a Cauchy sequence in $\mathcal{C}([0,\tau]; \Sc\times\Rd).$ Letting $\tilde{S}(.,\om)$ be the limit of this sequence, then by the cocycle property of $\Phi,$ i.e., $\Phi(t+k\tau, \theta_{-k\tau}\om, \tilde{y}) = \Phi(t, \om, \Phi(k\tau,\theta_{-k\tau}\om, \tilde{y}))$, we have 
\begin{align}\label{PS_L}
\lim_{k\rightarrow \infty}\Phi(r+k\tau, \theta_{-k\tau}\om, \tilde{y}) = \tilde{S}(r,\om), \quad \forall r\in [0,\tau].
\end{align}
Next, observe that the random variable $\tilde{S}(0,\om)$ is $\F_{-\infty}^{s}$-measurable and given the $\p$-preserving property of $\theta,$ we can derive the same property by replacing $\om$ by $\theta_t\om, \; t\in \T\subseteq\R$. Indeed, we have 
\begin{align*}
\Phi(t,\theta_r\om, \Phi(r,\om, \Phi(k\tau, \theta_{-k\tau}\om ,\tilde{y}))) &= \Phi(r+t,\om, \Phi(k\tau, \theta_{-k\tau}\om,\tilde{y}))\\
&= \Phi(t+r+k\tau, \theta_{-k\tau}\om, \tilde{y})\\
&= \Phi(r, \theta_t\om, \Phi(t+k\tau, \theta_{-t-k\tau}\theta_t\om, \tilde{y})),
\end{align*}
so that 
\begin{align}\label{PS_LL}
\Phi(m\tau, \theta_{r}\om, \tilde{S}(r, \om)) = \tilde{S}(r, \theta_{m\tau}\om), \quad \forall m\in \Z, \quad r\in [0, \tau].
\end{align}
The continuity of the flow map $(t,\tilde{y})\mapsto \Phi(t,\om,\tilde{y})$ for almost all $\om$ and the equality (\ref{PS_LL}), then for any $t, r\in \T\subset \R $ with $t\geqslant r,$ we have 
 \begin{align}\label{PS_L2}
\notag \Phi(t, \theta_r\om, \tilde{S}(r,\om))&= \Phi(t, \theta_r\om, \lim_{k\rightarrow\infty}\Phi(r+k\tau, \theta_{-k\tau}\om,\tilde{y}))\\
 \notag &=\lim_{k\rightarrow\infty}\Phi(t, \theta_r\om, \Phi(r+k\tau,\theta_{-k\tau}\om,\tilde{y}))\\
 &=\lim_{k\rightarrow\infty}\Phi(t+r+k\tau, \theta_{-k\tau}\om, \tilde{y})= \tilde{S}(t+r,\om).
 \end{align}
 Finally,  we show the exponential stability and the uniqueness of the limit $\tilde{S}(r,\om),\; r\in [0, \tau].$ For this, let $\tilde{Y}(\om),$ $\tilde{Z}(\om)$ be any two $\F_{-\infty}^{s}$-measurable random variables, we deduce from (\ref{Exp_fast}) that 
 \begin{align}\label{Uniq}
 \lim_{k\rightarrow\infty }\sup_{t\in [0,\tau]}\vert \Phi(t+k\tau, \om, \tilde{Y}(\om)) - \Phi(t+k\tau,\om, \tilde{Z}(\om))\vert = 0,
 \end{align}
 exponentially fast, $\p$-a.s. In particular, take $\tilde{Z}(\om) = \tilde{S}(r,\om),$ we obtain the  exponential stability and the uniqueness of $\tilde{S}(r,\om)$.
 \medskip
 
 {\bf Step II:}
Now, it only remains to show the random periodicity of $\tilde{S}(r,\om)$. We start by using that the coefficients $b, \sigma$ of our SDE are $\tau$-periodic in time, so that for $t\geqslant s,$
 \begin{align}\label{Uniq1}
 \notag  Y_{s,t}^y(\theta_{\tau}\om) &= y + \int_s^tb(u, Y^y_{s,u}(\theta_{\tau}\om))du+ \sum_{k=1}^N\int_s^t\sigma_k(u,Y_{s, u}(\theta_{\tau}\om))dW^k_{u+\tau}(\om)\\
 &= y+ \int_{s+\tau}^{t+\tau}b(u,Y^y_{s-\tau, u-\tau}(\theta_{\tau}\om))du+\sum_{k=1}^N\int_{s+\tau}^{t+\tau}\sigma_k(u, Y^y_{s-\tau,u-\tau}(\theta_{\tau}\om))dW_u(\om).
 \end{align}
 On the other hand,
 \begin{align}\label{Uniq2}
 Y^y_{s+\tau, s+\tau}(\om) = y+\int_{s+\tau}^{t+\tau}b(u, Y^y_{s,u}(\om))du+\sum_{k=1}^N\int_{s+\tau}^{t+\tau}\sigma(u, Y^y_{s,u}(\om))dW_u(\om).
 \end{align}
 Then by uniqueness of solution of our SDE and the $\p$-preserving property of $\theta,$ the equations (\ref{Uniq1}) and (\ref{Uniq2}) lead to 
 \begin{align}\label{Per_f}
 Y^y_{s+\tau,t+\tau}(\om) = Y^y_{s,t}(\theta_{\tau}\om), \quad t\geqslant s, \quad \p-a.s.
 \end{align} 
 Now, we return to the lifted flow $\Phi(t,\om,\cdot): \Sc\times\Yt\rightarrow \Sc\times\Yt$,
 \begin{align}\label{Lift2}
 \Phi(t+k\tau,\theta_{-k\tau}\om, \tilde{y})= (s+t+k\tau\mod \tau, Y^y_{s,t+k\tau+s}( \theta_{-s-k\tau}\om))), \quad t\geqslant 0.
 \end{align}
Define a projection map $\Pi: \Sc\times\Yt\rightarrow \Yt$,  $\tilde{S}\mapsto\Pi\tilde{S}\in \Yt$. Next, from the limit (\ref{PS_L}) and the random periodic property in (\ref{Per_f}) together with the continuity of the projection $\Pi$, we have 
\begin{align*}
\Pi \tilde{S}(t,\theta_{\tau}\om)& = \lim_{k\rightarrow\infty}\Pi\Phi(t+k\tau, \theta_{-k\tau+\tau}\om, \tilde{y}) \\
&= \lim_{k\rightarrow \infty}Y^y_{s,t+k+s}(\theta_{\tau}\theta_{-k\tau-s}\om)
= \lim_{k\rightarrow \infty} Y^y_{s+\tau,t+k\tau+s+\tau}(\theta_{-k\tau-s}\om)\\
& = \lim_{k\rightarrow \infty}\Pi \Big(s+t+k\tau\mod\tau, Y^y_{s+\tau,t+k\tau+s+\tau}(\theta_{-k\tau-s}\om)\Big)\\
&= \lim_{k\rightarrow \infty}\Pi\Phi(t+k\tau+\tau, \theta_{-k\tau}\om,\tilde{x}) = \Pi \tilde{S}(t+\tau, \om).
\end{align*}
The continuity of the flow map $(t,r,y)\mapsto Y^y_{r,t}(\om)$ for almost all $\om$ and the equality (\ref{PS_L2}), yield
\begin{align}\label{R11}
\notag Y^{\Pi \tilde{S}(r,\om)}_{r,t+r}(\om) &= \lim_{k\rightarrow\infty}Y_{r,t+r}^{\Pi\Phi(r+k\tau, \theta_{-k\tau}\om, \tilde{y})}(\theta_r\theta_{-r}\om)\\
\notag &=\lim_{k\rightarrow\infty}\Pi\Phi(t, \theta_r\om,\Pi\Phi(r+k\tau, \theta_{-k\tau}\om, \tilde{y})) )\\
\notag &= \lim_{k \rightarrow \infty}\Pi^2\Phi(t,\theta_r\om,\Phi(r+k\tau, \theta_{-k\tau}\om,\tilde{y}))\\
\notag &=\lim_{k\rightarrow \infty}\Pi\Phi(t+r+k\tau, \theta_{-k\tau}\om, \tilde{y})\\
&= \Pi\bar{S}(t+r, \om), \quad \p-\text{a.s.}
\end{align}
This implies that $\Pi\tilde{S}(r,\om)$ is a random periodic solution of the $\tau$-periodic flow $\{Y_{r,t}(\om): t,r\in \T\}$ on $\Yt$ (see Definition \ref{Chu}).
Set $\tilde{U}(r,\om) := (r\mod\tau, \Pi\tilde{S}(r,\om))$ for $r\in \T\subset\R,$ it is relative easy to verify that
 \begin{align}
{U}(r+\tau, \om) = \tilde{U}(r, \theta_\tau\om)
\end{align}
and  (\ref{R11}) leads to  
\begin{align}
\notag \Phi(t,\theta_r\om, \tilde{U}(r,\om))&= (r\mod\tau, Y_{r,t+r}^{\Pi\tilde{S}(r,\om)}(\om))\\ &= (r\mod\tau, \Pi\tilde{S}(t+r,\om))= \tilde{U}(t+r,\om), \quad\p-\text{a.s.}
\end{align}
This implies that $\tilde{U}(r,\om)$ is a random periodic solution of $\Phi$ on $\Sc\times\Yt$ and the uniqueness of the limit (\ref{PS_L}), gives us   $\tilde{S}(r,\om) = \tilde{U}(r,\om).$  
\qed
\begin{rem}\rm
 In practise, the construction of Lyapunov function is not simple and depends on a specific problem. However, one can construct a polynomial Lyapunov function growing at infinity like $\vert x\vert^{2m}, \; m\in\N\setminus\{0\}$, for a class of SDEs with coefficients $b\in \mathcal{C}^1(\R\times\Rd\times\Yt; \Yt),\;\; \sigma_{.k}\in \mathcal{C}^{2,\delta}_b(\R\times\times\Rd\times\Yt; \Yt)$ satisfying the following dissipative conditions:
\begin{align}\label{Dep1}
\begin{cases}
\langle b(t,x,y)-b(t,\hat{x},z), (x,y)-(\hat{x},z)\rangle\leqslant -K_t\left(\vert y-z\vert^2+|x-\hat{x}|^2\right),\\
\vert \sigma_k(t,x,y)-\sigma_k(t,\hat{x},z)\vert \leqslant L_t\left(\vert y-z\vert+|x-\hat{x}|\right), \;\; 1\leqslant k\leqslant N,
\end{cases}
\end{align}
such that 
\begin{align}\label{Dep2}
\limsup_{t\rightarrow\infty}\frac{1}{2t}\int_0^t\lambda(s)ds<0,
\end{align}
where $\lambda(t) = -K_t+\frac{(p-1)}{2}NL^2_t$ for some $p>1.$
The function $K_t$ is defined by 
\begin{align*}
K_t = \liminf_{R\rightarrow \infty}K_t(R),
\end{align*}
where $K_t: \R\rightarrow\R$ is a Borel function defined by 
\begin{align*}
K_t(R) = \inf\Big\{-\frac{\langle b(t,x,y) - b(t,\hat{x},z), (x,y)-(\hat{x},z)\rangle}{\vert y-z\vert^2+|x-\hat{x}|^2 }: \vert y-z\vert+|x-\hat{x}| =R\Big\}.
\end{align*}
There are many important class of SDEs in theory and applications satisfying the above dissipative conditons (\ref{Dep1}) -- (\ref{Dep2}). The existence and uniqueness of random periodic solution for a class of SDEs satisfying related conditions was proved in \cite{Uda14, Uda18}. 
\end{rem}

\subsection{Ergodic periodic measures on Poincar\'e sections}\label{Ranp_M}
In this subsection, we discuss the ergodicity of the fast subsystem (\ref{Fast2}) in the random periodic regime. Many studies of invariant measures  for Markovian RDS generated by SDEs are normally via stochastic analysis  or dynamical systems perspectives (e.g.,~\cite{Arnold, Crau2}). From stochastic analysis perspective, invariant measures are investigated via Markov transition function, in this sense, ergodicity is based on the dynamics of Markov evolution. From dynamical description, one investigates random invariant probability measures whose conditional expectation with respect to a sub-algebra of $\F$ has one-to-one correspondence with the invariant measure of the Markov evolution. Here, we are interested in capturing the ergodicity of transition probability function of random periodic solutions (the law of random periodic solutions). Ergodicity  on Poincar\'e sections (PS-ergodcity) is a form of ergodicity suitable in this case. PS-ergodicity of stochastic dynamical systems was recently developed by Feng and Zhao \cite{Feng18}, it gives a new perspective and a generalised form of Poincar\'e-Bendixson theorem for stochastic dynamical systems. We would like to argue that the random periodic measure of the fast subsystem (\ref{Fast}) is PS-ergodic under  Assumption \ref{A2.1} and restricted H\"omander Lie bracket conditions.
\medskip

We start with preparatory definitions leading to the presentation of our ergodicity result (Theorem \ref{Ps_erg}).
\begin{definition}[Periodic measure \cite{Feng18}]\;\label{Rand_pm}\rm\;
Let $\mathbb{M}$ be a Polish space, a  measure $\mu: \R\rightarrow \mathcal{P}(\mathbb{M})$ is called a {\it periodic measure} of period $\tau>0$, on the phase space $(\mathbb{M},\Bb(\mathbb{M}))$ for the Markovian stochastic flow $\{Y_{r,t}(\om): \;t\geqslant r\}$, if  for $B\in \Bb(\mathbb{M})$, the following hold
 \begin{equation}\label{10.11}
 \mu_{r+\tau} = \mu_r\\ \quad \text{and} \quad \mu_{t+r}(B) = \int_{\mathbb{M}}P(r,x; t+r, B)\mu_r(dx) = (\mathcal{P}^*_{r,t}\mu_r)(B), \quad r\in\R, \quad t\in \Rp.
 \end{equation}
 It is called a periodic measure with minimal period $\tau,$ if $\tau>0$ is the smallest number such that (\ref{10.11}) holds. 
 \end{definition} 
 Let $\tilde{S}: \R\times\Om\rightarrow \Sc\times\Yt$ be a random periodic solution of the RDS $\{\Phi(t,\om): t\geqslant 0\}$. Consider a probability measure $\tilde{\mu}_t\in \mathcal{P}(\Sc\times\Yt)$ defined by 
\begin{align}
  \tilde{\mu}_r(\tilde{A}) := \left(\p\circ \tilde{S}^{-1}(r,.)\right)(\tilde{A}) = \p(\{\om: \tilde{S}(r,\om)\in \tilde{A}\}), \;\; \tilde{A}\in \Bb(\Sc)\otimes\Bb(\Yt).
 \end{align} 
Then the measure $\mu_r$ is $\tau$-periodic as
 \begin{align}\label{PM3}
\nonumber \mu_{r+\tau}(A)= \p\{\om: \tilde{S}(r+\tau,\om)\in \tilde{A}\}
&= \p\{\om: \tilde{S}(r, \theta_{\tau}\om)\in \tilde{A}\}\\
& = \p\{\om: \tilde{S}(r,\om)\in \tilde{A}\}
 =\mu_r(\tilde{A}),
 \end{align}
 Moreover, as it was shown in \cite{Feng18}, $\tilde{\mu}_r$ satisfies (\ref{10.11}). Thus, the law of random periodic solution satisfies Definition (\ref{Rand_pm}). 
 
 Next, we define a Poincar\'e section for the lifted Markov semigroup $(\tilde{\mathcal{P}})_{t\geqslant 0}$ on $\Sc\times\Yt.$ 
\bn[Poincar\'e section for lifted Markov semigroup \cite{Feng18}]\;\rm The Borel subsets of  $\{\Bt_r: r\in \T\subseteq\R\}\subset\Bb(\Sc)\otimes\Bb(\Yt)$ are called Poincar\'e sections for the Markov semigroup $(\tilde{\mathcal{P}}_t)_{t\geqslant 0}$ on the cylinder $\Sc\times\Yt$, if 
\begin{align}\label{Poinc}
\Bt_{r+\tau} = \Bt_r, \quad \text{and}\quad \tilde{P}(0, \tilde{y}; t, \Bt_{t+r}) = 1, \quad \forall \tilde{y}\in \Bt_r, \quad t\geqslant 0.
\end{align}
\en
\begin{rem}\rm\; \label{minP}
The choice of Poincar\'e section is not unique, for example $\Bt_r= \Sc\times\Yt$ and $\Bt_r= \text{supp}(\tilde{\mu}_r)$ satisfy the definition of Poincar\'e section. However, the family $\{\Bt_r= \text{supp}(\tilde{\mu}_r): r\in\R\}$ is a minimal Poincar\'e section or $k_r\tau$-irreducible Poincar\'e section \cite{Feng18}. To see this, fix $r\in [0,\tau)$ and any open set  $A\subset\Yt$ such that  $\tilde{A_r}:=\{r\}\times A\subset \Bt_r= \textrm{supp}(\tilde{\mu}_r)$ with $\mu_r(\Bt_r\backslash\tilde{A}_r)>0$, we have for all $\tilde{y}\in \Bt_r,$
\begin{align*}
\tilde{P}(0, \tilde{y}; k\tau+r, \tilde{A}_r)<1.
\end{align*}
This implies that $\tilde{A}_r$ is not a Poincar\'e section for the Markov semigroup $(\tilde{\mathcal{P}}_{t})_{t\geqslant 0}$ on the cylinder $\Sc\times\Yt$.
\end{rem}
Let the periodic measure $(\tilde{\mu}_r)_{r\in\R}\subset\mathcal{P}(\Sc\times\Yt)$ be given, we consider its time average $\bar{\tilde{\mu}}$ over the interval $[0, \tau)$ defined by 
 \begin{align}
 \bar{\tilde{\mu}}(\tilde{A})  = \frac{1}{\tau}\int_0^\tau \tilde{\mu}_r(\tilde{A})dr, \quad \tilde{A}\in \Bb(\Sc)\otimes\Bb(\Yt).
 \end{align}
 It is realtively easy to check that this time average measure $\bar{\tilde{\mu}}$ is invariant under the Markov semigroup $(\tilde{\mathcal{P}}_{t})_{t\geqslant 0}$ (see also~\cite{Feng18}). Moreover, from the definition of random periodic solution (Definition~\ref{Chu}) and $\p$-preserving property of $\theta,$ we have for any $\tilde{A}\in \Bb(\Sc\times\Yt),$ 
 \begin{align*}
 \bar{\tilde{\mu}}(\tilde{A})  = \frac{1}{\tau}\int_0^\tau \tilde{\mu}_r(\tilde{A})dr &= \frac{1}{\tau}\int_0^\tau\p\{\om: \tilde{S}(r,\om)\in \tilde{A}\}dr\\
 &=\frac{1}{\tau}\E\left[\int_0^\tau\I_{\tilde{A}}(\tilde{S}(r, \om)dr\right]\\
&=\E\left[\frac{1}{\tau}m(\{r\in [0, \tau): \tilde{S}(r,\om)\in \tilde{A}\})\right]\\
 &=\E\left[\frac{1}{\tau}m(\{r\in [u, u+\tau): \tilde{S}(r,\om)\in \tilde{A}\})\right].
 \end{align*}
 This implies that the expected time spent in a Borel set $\tilde{A}\in\Bb(\Sc\times\Yt)$ by the random periodic path $r\mapsto\tilde{S}(r,\om)$ over a time interval of exactly one period at any time is independent of the starting point.  More precisely, this leads to the notion of PS-ergodicity of periodic measures.
\bn[PS-ergodicity \cite{Feng18}]\label{PS-eR}\;\rm
 A family of $\tau$-periodic measures $(\tilde{\mu}_r)_{r\in\R}$ is said to be {\it PS-ergodic} if for each $r\in [0,\tau),$ $\tilde{\mu}_r$ as an invariant measure of the Markov semigroup $(\tilde{\mathcal{P}}^x_{k\tau+r})_{k\in\N},$ at the integral multiples of the period $\tau$ on the Poincar\'e section $\Bt_r$ is ergodic.
 \en
 Equivalently, a family of $\tau$-periodic measures $(\tilde{\mu}_r)_{r\in\R}$ is PS-ergodic, if for any $\tilde{A}\in \Bb(\Sc\times\Yt)$ with $\tilde{A}\subset \Bt_r,$ the following holds
 \begin{align}\label{Krybo}
\lim_{m\rightarrow\infty}\int_{\Sc\times\Yt}\left\vert \int_{0}^\tau\Big[ \frac{1}{m}\sum_{k=0}^{m-1}\tilde{P}(0, \tilde{y}; r+k\tau, \tilde{A})-\tilde{\mu}_r(\tilde{A})\Big] dr\right\vert \bar{\tilde{\mu}}(d\tilde{y}) =0.
 \end{align}
 The limit (\ref{Krybo}) is the Krylov-Bogolyubov scheme for periodic measures \cite{Feng18, Hasm80}.
\newcommand{\PT}{\tilde{\mathcal{P}}}
Now, given the random periodic paths $\{\tilde{S}^x(r,\om): r\in [0, \tau)\},$ we want to show that $\{\tilde{\mu}_r^x: r\in [0, \tau)\}$ is ergodic under the discrete Markov semigroup $(\tilde{\mathcal{P}}^x_{r+k\tau})_{k\in \N}$ for fixed $x\in \Nb\Subset \Rd$ under the Assumption \ref{A2.1}. Precisely, we want to prove the convergence of the Krylov--Bogolyubov scheme for periodic measures. To this end, we start with the following simple and useful Lemma.
\begin{lem}\label{Expon}\rm\;
Let $(\tilde{\mu}^x_r)_{r\in [0, \tau)},\; x\in \Nb$ be $\tau$-periodic measures generated by the random periodic solutions of the SDE (\ref{Fast2}). If there exist $0<K<\infty$ and  $\lambda\in L^1(\R)$ such that 
\begin{align}\label{Negint}
\limsup_{t\rightarrow\infty}\frac{1}{2t}\int_{0}^{t}\lambda(u)du<\beta<0,
\end{align} 
and for all $f\in \mathcal{C}_b(\Bt_r), \;\; p\geqslant 1$,
\begin{align}\label{C2}
\left\vert \tilde{\mathcal{P}}^x_{r+k\tau}f-\int_{\Bt_r}fd\tilde{\mu}^x_r\right\vert\leqslant K||f||_{\infty}\bigg[\exp\left(p^{-1}\int_0^{r+k\tau}\lambda(u)du\right) + \exp\Big( \frac{1}{2p}\beta k\tau\Big)\bigg].
\end{align}
Then,
\begin{align}\label{Krybo2}
\lim_{m\rightarrow\infty}\int_{\Sc\times\Yt}\left\vert \int_{0}^\tau\Big[ \frac{1}{m}\sum_{k=0}^{m-1}\tilde{P}^x(0, \tilde{y}; r+k\tau, \tilde{A})-\tilde{\mu}^x_r(\tilde{A})\Big] dr\right\vert \bar{\tilde{\mu}}^x(d\tilde{y}) =0.
 \end{align}
\end{lem}
\noindent {\it Proof.} 
The proof is relatively straightforward, the idea is to employ density argument. For this, let $\tilde{A}_r\subset \Bt_r\subset\Sc\times\Yt$ be a closed set,  take $f_r = \I_{\tilde{A}_r}$ and let the sequence of functions $(f_n)_{n\in \N}$ be defined by 
\begin{align*}
f_n(\tilde{y}) = \begin{cases} 1, & \quad \text{if}\quad \tilde{y}\in \tilde{A}_r,\\
1-2^nd(\tilde{y}, \tilde{A}_r), & \quad \text{if}\quad d(\tilde{y},\tilde{A}_r)\leqslant2^{-n},\\
0, &\quad \text{if} \quad d(\tilde{x},\tilde{A}_r)\geqslant 2^{-n},
\end{cases}
\end{align*}
where $d(\tilde{y},\tilde{A}_r) = \inf\{\vert \tilde{y}-\tilde{z}\vert: \tilde{z}\in \tilde{A}_r \},\; \tilde{y}\in \Bt_r$. Then,
\begin{align*}
f_n(\tilde{y})\rightarrow f_r(\tilde{y}), \quad \text{as}\; n\rightarrow\infty, \quad \tilde{y}\in \Bt_r.
\end{align*}
Consequently, for all $r\in [0,\tau),$ we have 
\begin{align*}
\tilde{\mathcal{P}}^x_{r+k\tau}f_n(\tilde{y})\rightarrow \tilde{\mathcal{P}}^x_{r+k\tau}f_r(\tilde{y}) = \tilde{\mathcal{P}}^x_{r+k\tau}\I_{\tilde{A}_r}(\tilde{y}).
\end{align*}
This implies that   $\tilde{\mathcal{P}}^x_{r+k\tau}f_n\in \mathcal{C}_b(\Bt_r)$, so that by (\ref{C2}) and as $\tilde{\mu}^x_r$ is invariant w.r.t. $\tilde{\mathcal{P}}^x_{r+k\tau},$ we have
\begin{align*}
\left\vert \tilde{P}^x(0,\tilde{y};r+k\tau,\tilde{A}_r) -\tilde{\mu}^x_r(\tilde{A}_r)\right\vert=\left\vert \tilde{\mathcal{P}}^x_{r+k\tau}\I_{\tilde{A}_r}(\tilde{y})-\tilde{\mu}^x_r(\tilde{A}_r)\right\vert \leqslant & K\bigg[ \exp\left(p^{-1}\int_0^{r+k\tau}\lambda(u)du\right)\\ & \hspace{1cm}+\exp\Big( \frac{1}{2p}\beta k\tau\Big)\bigg].
\end{align*}
It then follows by covering lemma that for any $\tilde{A}\in \Bb(\Sc\times\Yt)$  with $\tilde{A}\subset \Bt_r,$ 
\begin{align*}
\int_0^\tau\left\vert \tilde{P}^x(0,\tilde{y};r+k\tau,\tilde{A}) -\tilde{\mu}^x_r(\tilde{A})\right\vert dr&\leqslant K\bigg[\int_0^{\tau}\exp\left(p^{-1}\int_0^{r+k\tau}\lambda(u)du\right)dr + \exp\Big( \frac{1}{2p}\beta k\tau\Big)\bigg].\\
&= K\bigg[\int_0^{\tau}\exp\left(\frac{1}{2pk\tau}\int_0^{r+k\tau}\lambda(u)du\right)^{2k\tau}dr +\exp\Big( \frac{1}{2p}\beta k\tau\Big)\bigg].
\end{align*}
By the condition (\ref{Negint}), there exist $0<\tilde{K}<\infty, \;\; 0<\tilde{\beta}<1$, such that 
\begin{align*}
\int_0^\tau\left\vert \tilde{P}^x(0,\tilde{y};r+k\tau,\tilde{A}) -\tilde{\mu}^x_r(\tilde{A})\right\vert dr\leqslant \tilde{K}\tilde{\beta}^{k\tau},
\end{align*}
leading to 
\begin{align*}
\int_{\Sc\times\Yt}\left\vert \int_{0}^\tau\Big[ \frac{1}{m}\sum_{k=0}^{m-1}\tilde{P}^x(0, \tilde{y}; r+k\tau, \tilde{A})-\tilde{\mu}^x_r(\tilde{A})\Big] dr\right\vert\bar{\tilde{\mu}}^x(d\tilde{y}) \leqslant \tilde{K}\sum_{k=0}^{m-1}\tilde{\beta}^{k\tau}.
\end{align*}
This gives the desired convergence.
\qed

\vspace{.5cm}
It is not very simple to verify the inequality (\ref{C2}), some non-degenerate conditions would be required. In fact, the sufficient condition for inequality (\ref{C2}) to hold is the strong Feller property of the Markov semigroup $(\tilde{\mathcal{P}}^x_{k\tau+r})_{k\in\N}$. Recall, that a Markov semigroup $(\tilde{\mathcal{P}}^x_{k\tau+r})_{k\in\N}$ has strong Feller property, if for any $f\in \Bb_b(\Sc\times\Yt)$, we have $\tilde{\mathcal{P}}^x_{k\tau+r}f\in \mathcal{C}_b(\Sc\times\Yt)$. Equivalently, a Markov semigroup $(\tilde{\mathcal{P}}^x_{k\tau+r})_{k\in\N}$ has strong Feller property, if and only if 
\btm[leftmargin = .7cm]
\item[(i)] $(\tilde{\mathcal{P}}^x_{k\tau+r})_{k\in\N}$ is a Feller semigroup, i.e.,$\tilde{\mathcal{P}}^x_{k\tau+r}: \mathcal{C}_b(\Sc\times\Yt)\rightarrow\mathcal{C}_b(\Sc\times\Yt),$  and 
\item[(ii)] the family $\{\tilde{P}(0,\tilde{y}_m; t+k\tau,.): m\in \N\}$ is equicontinuous.
\etm
The first item follows from the existence of stochastic flows (e.g., \cite{Kunita, Hasm80}) and the second item require non-degeneracy and certain integrability of Malliavin covariance $\tilde{C}^x_t(\om,\tilde{y})$ of the flows of solutions of the fast motion $\tilde{Y}_t^{x,\tilde{y}}$. Intuitively, the strong Feller property states that for any nearby initial data $\tilde{y}$ and $\tilde{z}$ and any realisation $\om$ of the past of the driving noise, it is possible to construct a coupling between two solutions $\Phi(t,\om,\tilde{y})$ and $\Phi(t,\om,\tilde{z})$ such that with a probability close to $1$ as $\tilde{y}\rightarrow \tilde{z},$ one has $\Phi(t,\om,\tilde{y}) = \Phi(t,\om,\tilde{z})$ for $t\geqslant 1$ (e.g.,~\cite{Hairer11a, Hairer_Note16, Hairer11}). As elucidated by Hairer, Mattingly and others (e.g.,~\cite{Hairer11a, Hairer11, Hairer11c, HairerD}), one way of achieving such a coupling is via change of measure on driving process for one of the two solutions such that the noises $W_t^{\tilde{y}}$ and $W_t^{\tilde{z}}$ driving the solutions $\Phi(t,\om,\tilde{y})$ and $\Phi(t,\om,\tilde{z})$ respectively, are related by 
\begin{align*}
dW_t^{\tilde{z}} = dW_t^{\tilde{y}}+ v^{\tilde{y},\tilde{z}}_tdt,
\end{align*}
where $v$ is a control process that steers the solution $\Phi(t,\om,\tilde{y})$ towards the solution $\Phi(t,\om,\tilde{z})$ (cf.,~\cite{Hairer11c}). If one sets $\tilde{z} = \tilde{y}+\varepsilon\eta$ and look for controls of the form $v^{\tilde{y},\tilde{z}} = \varepsilon v,$ then as $\varepsilon\rightarrow 0,$ the random variable $v$ will induce a deformation onto the solution $\Phi(t,\om,\tilde{y})$ after time $t$ in the form of {\it Malliavin derivative} \footnote{ see equation (\ref{Mallia}) below.} of $\Phi(t,\om,\tilde{y})$ at $\om$ in the direction of $v$, i.e.,
\begin{align*}
\D\Phi(t,\om,\tilde{y})\cdot v=\D_v\Phi(t,\om,\tilde{y}) =:\mathcal{A}_tv.
\end{align*}
On the other hand, the effect of the perturbation of the initial condition by $\eta$ is given as the directional deriative of the solution $\Phi(t,\om,\tilde{y})$ at $\tilde{y}$ in the direction of $\eta$, i.e.,
\begin{align*}
D_x\Phi(t,\om,\tilde{y})\eta = J_{0,t}\eta.
\end{align*}
As extensively discussed in \cite{Hairer11a} for a more general and highly degenerate infinite-dimensional stochastic dynamical systems, in order to prove strong Feller property, the main task is to find a control $v$ such that 
\begin{align}\label{Control}
J_{0,t}\mathcal{A}_tv = J_{0,t}\eta, \quad \text{or} \quad \mathcal{A}_tv = \eta.
\end{align}
To be able to construct the above control process satisfying (\ref{Control}) with appropriate integrable condition, we impose the following restricted H\"omander's Lie bracket condition.
 \begin{Assum}\rm\label{A2.1bc}\
Let $b^x\in \mathcal{C}_b^{l-1}(\R\times\R^N; \R^N)$ and $\sigma_{.k}^x\in \mathcal{C}_b^{l, \delta}(\R\times\R^N; \R^N), \; 1\leqslant k\leqslant N,$ for $l\geqslant 2$ such that 
\begin{align}
\text{dim}(\text{span Lie}\left\{\sigma^x_{.k}: 1\leqslant k\leqslant N\right\})(0,x) = N,
\end{align}
or equivalently, there exists $M(y)=:M\in \N, \; C_M(y) =: C_M >0,$ such that for all $\eta\in \Scd$ we have 
\begin{align}
\sum_{\ell=0}^{M}\sum_{Z\in \Sigma^x_\ell}(\eta\cdot Z)^2(0, y)\geqslant C_M
\end{align}
where $\Sigma^x_0 = \{\sigma^x_{.k}: 1\leqslant k\leqslant N\},\;\; \Sigma^x_{\ell+1}=\{[\sigma^x_{.k}, Z]: 1\leqslant k\leqslant N,\; Z\in \Sigma^x_{\ell}\}$ and $[F, G]$ is the Lie bracket between the vector fields $F$ and $G$ defined by 
\begin{align*}
 [F,G](t,y):= D_yG(t,y)F(t,y)-D_yF(t,y)G(t,y).
\end{align*}
\end{Assum}
The restricted H\"omander brackets condition above generalises uniform ellipticity condition; this follows by noting that uniform ellipticity condition implies that for every point $(t, y)$ in the neigbhourhood of $(0, y),$ that the linear span of $\sigma^x_{.k}(t,y); 1\leqslant k\leqslant N,$ is the whole of $\R^N.$  

For convenience, we briefly recall derivation on the Wiener space $\Om$, we refer interested reader to (e.g.,~\cite{Hairer10, Hairer_Note16, Hairer11, Malliavin, Millet, Nulart, Watanabe} ) for accessible treatment. Let $\mathcal{H} = L^2([0, T];\R^N)$ and let $E$ be a real separable Hilbert space. An $E$-valued random variable $\Theta: \Om\rightarrow E$ will be called smooth if it admit the following representation
\begin{align*}
\Theta = \sum_{j=1}^M \vartheta_j\left(W_{t_1},\cdots,W_{t_m}\right)v_j,
\end{align*}
where $\vartheta_j \in \mathcal{C}_b^{\infty}(\R^{Nm}), \; t_1,\cdots, t_m\in [0, T]$ and $v_1,\cdots, v_M\in E.$ \\ The derivative $\D$ of a smooth random variable is a random variable taking values on the Hilbert space $\mathcal{H}\otimes E = L^2([0, T]; \R^N\otimes E)$ given by 
\begin{align}\label{Mallia}
\D_t^i\Theta  =\sum_{j=1}^M\sum_{k=1}^m\frac{\partial \vartheta_j}{\partial x^{ik}}\left( W_{t_1},\cdots, W_{t_m}\right)\I_{[0, t_k]}(t)v_j, 
\end{align}
for each $t\in [0, T]$ and $i=1, \cdots, N.$\\  Observe that for any $h\in \mathcal{H},$ the random variable $\sum_{i=1}^N\int_0^T\D_t^i\Theta h^i(t)dt$ can be interpreted as the directional derivative
\begin{align*}
\frac{d}{d\varepsilon}\bigg|_{\varepsilon =0}\Theta\left(\om(\cdot)+\varepsilon\int_{0}^{\cdot}h(s)ds\right),
\end{align*}
where $W_t(\om) = \om(t)$.

Moreover, to establish strong Feller property, the following preparatory lemmas would be useful.  We skip the dependence on $x$ in the proofs to simplify notation.
\begin{lem}\label{Malliavin1est}\rm
Let $b^x, \sigma^x_{.k}\in \mathcal{C}_b^l(\R\times\R^N;\R^N), \; 1\leqslant k\leqslant N,\; l\geqslant 2,$ and $\Phi(t,\om,\tilde{y})$ be stochastic flow generated by the lifted SDE (\ref{Fast2}) with the Jacobian matrix $J_{0,t}(\om, \tilde{y})$ at $\tilde{y}\in \Sc\times\R^N$. Then, for all $t>0$ and $p\geqslant 1,$ the following hold, 
\begin{align}
\begin{cases}
\int_0^t\E\Vert \D_r^i\Phi(u,\om,\tilde{y})\Vert^{2p} \leqslant\frac{ (N+1)^{2p}}{pL_{b\sigma}}\exp(-pL_{b\sigma}r)\left[ \exp(pL_{b\sigma}t)-1\right],\\
\int_0^t \E\Vert \D_r^iJ_{0,u}\Vert^{2p}du\leqslant\frac{\tilde{L}^{2p}_{b\sigma}(N+1)^{2p}}{pL_{b\sigma}}\exp(-pL_{b\sigma}r)\Big[\frac{1}{3}\exp(3pL_{b
 \sigma}t)-\exp(pL_{b\sigma}t)\Big].
\end{cases}
\end{align}
where $L_{b\sigma}, \tilde{L}_{b\sigma}$ are constants depending on $\tilde{b}, \tilde{\sigma}$ and their derivatives up to second order.

Moreover, there exist constants $0<\tilde{C}_j(N,p,L_{b\sigma})<\infty, $ $j=1,2.,$ depending on $N,\; p,$ the coefficients $b,\sigma$ and their derivatives w.r.t $\tilde{y} =(s,y)$ up to second order such that  
 \begin{align}
 \notag &\hspace{-1cm}\E\Vert \D_r^i\hat{C}_u\Vert^{2p}\leqslant \exp(-2pL_{b\sigma}r)\bigg\{C_1(N,p,L_{b\sigma})\left[\frac{1}{12pL_{b\sigma}}\exp(12pL_{b\sigma}u)-\frac{1}{4pL_{b\sigma}}\exp(4pL_{b\sigma}u)\right]^{1/2}\\ &\left[\exp(8pL_{B\sigma}u)-1\right]^{1/2}
 +C_2(N,p, L_{b\sigma})\left[\exp(4pL_{b\sigma}u)-1\right]^{1/2}\left[\exp(4pL_{b\sigma}u)-4pL_{b\sigma}u-1\right]^{1/2}\bigg\},
 \end{align}
 where $\hat{C}_u(\om,\tilde{y}) := \left\langle \D\Phi(u,\om,\tilde{y}), \D\Phi(u,\om,\tilde{y})\right\rangle$ is the Malliavin covariance of the flow map $\Phi(u,\om,\tilde{y})$, also defined for $\xi\in \Sc\times\R^N,$ as follows 
\begin{align*}
\mathcal{A}_u^T\xi &= \tilde{\sigma}^T(\Phi(u,\om,\tilde{y}))(J_{0,u}^{-1})^T\xi \\
\hat{C}_u\xi&= \mathcal{A}_u\mathcal{A}_u^T\xi = \int_0^u J_{0,\ell}^{-1}\tilde{a}(\Phi(\ell,\om,\tilde{y}))(J_{0,\ell}^{-1})^T\xi d\ell.
\end{align*}
\end{lem}
\noindent {\it Proof.}
First, use that $b, \sigma_{.k}\in \mathcal{C}_b^{l}(\R\times\R^N; \R^N), \; l\geqslant 2,$ and Gronwall's inequality, for all $u\in [0, t]$ and $p\geqslant 1,$ to obtain 
\begin{align}\label{p2}
\begin{cases} \E\Vert J_{0,u}^{-1}\Vert^{2p} \leqslant (N+1)^p\exp(pL_{b,\sigma}u),\\  \E\Vert J_{0,u}\Vert^{2p} \leqslant (N+1)^p\exp(pL_{b,\sigma}u)
\end{cases}
\end{align}
 Let $u, r\in [0, t]$, then the Malliavin derivative of $\Phi(t,\om,\tilde{x})$ for $u>r$ is 
 \begin{align*}
&\hspace{-2cm} \D_r^{i}\Phi(u,\om,\tilde{y}) = \int_r^uD\tilde{b}(\Phi(\ell,\om,\tilde{y}))\D^i_r\Phi(\ell,\om,\tilde{y})d\ell \\ &\hspace{1cm} +\sum_{k=1}^m \int_r^uD\tilde{\sigma}_k(\Phi(\ell,\om,\tilde{y}))\D^i_r\Phi(\ell,\om,\tilde{y})dW_{\ell}^k+\tilde{\sigma}_i(\Phi(r,\om,\tilde{y})),
 \end{align*}
 and  for $u\leqslant r$, we have $\D_r^{i}\Phi(u,\om,\tilde{y}) = 0.$ In fact, the Malliavin derivative is simply the composition of the derivative flows $J_{0,u}\circ J_{0,r}^{-1}\tilde{\sigma}_i(\Phi(r,\om,\tilde{y})) = J_{r,u}\tilde{\sigma}_i(\Phi(r,\om,\tilde{y}))$ starting from  $\tilde{\sigma}_i(\Phi(r,\om,\tilde{y}))$ at time $u=r$.
So, by the regularity of the coefficients $\tilde{b}^x, \tilde{\sigma}^x$ and Gronwall's inequality on $J_{r,u}$, one arrive at  
\begin{align*}
\E\Vert \D_r^i\Phi(u,\om,\tilde{y})\Vert^{2p} = \E\Vert J_{r,u}\tilde{\sigma}_i(\Phi(r,\om,\tilde{y}))\Vert^{2p}\leqslant (N+1)^{2p}\exp(pL_{b\sigma}(u-r)).
\end{align*}
 Integrating over $u\in [0,t]$, we have 
 \begin{align*}
 \int_0^t\E\Vert \D_r^i\Phi(u,\om,\tilde{y})\Vert^{2p}du \leqslant\frac{ (N+1)^{2p}}{pL_{b\sigma}}\exp(-pL_{b\sigma}r)\left[ \exp(pL_{b\sigma}t)-1\right].
 \end{align*}
 For the second assertion, we recall from (e.g.,~Proposition 5.2 in \cite{Hairer_Note16b}) that the Malliavin derivative of the Jacobian matrix $J_{0,u}, \; u\in [0, t],$  solves the following variational equation
 \begin{align*}
  \D^i_rJ_{0,u}
&=  \int_0^u J_{\ell,u}D^2\tilde{b}(\Phi(\ell,\om,\tilde{y}))\Big(J_{0,\ell}, \D_r^i\Phi(\ell,\om,\tilde{y})\Big)d\ell\\
&\hspace{1cm} +\int_0^u J_{\ell,u}D^2\tilde{\sigma}(\Phi(\ell,\om,\tilde{y}))\Big(J_{0,\ell}, \D_r^i\Phi(\ell,\om,\tilde{y})\Big)dW_{\ell}.
 \end{align*}
 So that,
 \begin{align*}
 \E\Vert \D_r^iJ_{0,u}\Vert^{2p}&\leqslant 2^{2p-1}\E\int_0^u\Vert J_{\ell,u}D^2\tilde{b}(\Phi(\ell,\om,\tilde{y}))\Big(J_{0,\ell}, \D^i_r\Phi(\ell,\om,\tilde{y})\Big)\Vert^{2p}d\ell \\ &\hspace{1cm}+2^{2p-1}\E\int_0^u\Vert J_{\ell,u}D^2\tilde{\sigma}(\Phi(\ell,\om,\tilde{y}))\Big(J_{0,\ell}, \D^i_r\Phi(\ell,\om,\tilde{y})\Big)\Vert^{2p}d\ell\\
& \leqslant 2^{2p}\tilde{L}^{2p}_{b\sigma}\E\left(\Vert J_{0,u}\Vert^{2p}\int_0^u\Vert \D_r^i\Phi(\ell,\om,\tilde{y})\Vert^{2p} d\ell\right)\\
& \leqslant 2^{2p}\tilde{L}^{2p}_{b\sigma}(N+1)^p\exp(pL_{b\sigma}u)\int_0^u\E\Vert \D_r^i\Phi(\ell,\om,\tilde{y})\Vert^{2p} d\ell\\
&\leqslant \frac{2^{2p-1}\tilde{L}^{2p}_{b\sigma}(N+1)^{3p}}{pL_{b\sigma}}\exp(pL_{b\sigma}(u-r))\big[\exp(2pL_{b\sigma}u)-1\big].
 \end{align*}
 It then follows that 
 \begin{align*}
\int_0^t \E\Vert \D_r^iJ_{0,u}\Vert^{2p}du\leqslant \frac{\tilde{L}^{2p}_{b\sigma}(N+1)^{2p}}{pL_{b\sigma}}\exp(-pL_{b\sigma}r)\Big[\frac{1}{3}\exp(3pL_{b
 \sigma}t)-\exp(pL_{b\sigma}t)\Big].
 \end{align*}
 For the final assertion, by product rule of differentiation, we have 
 \begin{align*}
 \D^i_r\hat{C}_u = \int_0^u\D^i_r[J_{0,\ell}^{-1}]\tilde{a}(\Phi(\ell,\om,\tilde{y}))(J_{0,\ell}^{-1})^Td\ell+ \int_0^tJ_{0,\ell}^{-1}\D_r^i[\tilde{a}(\Phi(\ell,\om,\tilde{y}))](J_{0,\ell}^{-1})^Td\ell \\ + \int_0^tJ_{0,\ell}^{-1}\tilde{a}(\Phi(\ell,\om,\tilde{y}))\D_r^i[(J_{0,\ell}^{-1})^T]d\ell.
 \end{align*}
 So that 
 \begin{align*}
 \E\Vert \D_r^i\hat{C}_u\Vert^{2p} &\leqslant 4^{2p-1}\int_{0}^u\E\Vert J_{0,\ell}^{-1}\D_r^i[J_{0,\ell}]J_{0,\ell}^{-1}\tilde{a}(\Phi(\ell,\om,\tilde{y}))(J_{0,\ell}^{-1})^T\Vert^{2p}d\ell \\ &\hspace{1cm} + 4^{2p-2}\int_0^u\E\Vert J_{0,\ell}^{-1}D\tilde{a}(\Phi(\ell,\om,\tilde{y}))\D_r^i[\Phi(\ell, \om,\tilde{y})] J_{0,\ell}^{-1}\Vert^{2p} d\ell\\
& \leqslant 4^{2p-1}\int_0^u\E\Big[\Vert J_{0,\ell}^{-1}\Vert^{2p}\Vert \D_r^iJ_{0,\ell}\Vert^{2p}\Vert \D_r^i\Phi(\ell,\om,\tilde{y})\Vert^{4p}\Big]d\ell\\
&\hspace{1cm} +4^{2p-2}\int_0^u\E\Big[\Vert J_{0,\ell}^{-1}\Vert^{4p}\Vert D\tilde{a}(\Phi(\ell,\om,\tilde{y}))\Vert^{2p}\Vert \D_r^i\Phi(\ell,\om,\tilde{y})\Vert^{2p}\Big]d\ell\\
 &\leqslant 4^{2p-1}\int_0^u\Big[\E\Vert J_{0,\ell}^{-1}\Vert^{4p}\Vert \D_r^iJ_{0,\ell}\Vert^{4p}\Big]^{1/2}\Big[\E\Vert \D_r^i\Phi(\ell,\om,\tilde{y})\Vert^{4p}\Big]^{1/2}d\ell\\
&\hspace{1cm} + 4^{2p-2}\int_0^u\Big[\E\Vert J_{0,\ell}^{-1}\Vert^{8p}\Vert\D\tilde{a}(\Phi(\ell,\om,\tilde{y}))\Vert^{4p}\Big]^{1/2}\Big[\E\Vert \D_r^i\Phi(\ell,\om,\tilde{x})\Vert^{4p}\Big]^{1/2}d\ell\\
& \leqslant 4^{2p-1}\left(\int_0^u\E\Big[\Vert J_{0,\ell}^{-1}\Vert^{8p}\Vert\D_r^iJ_{0,\ell}\Vert^{4p}\Big]d\ell\right)^{1/2}\left(\int_0^u\E\Vert \D_r^i\Phi(\ell,\om,\tilde{y})\Vert^{4p}d\ell\right)^{1/2}\\
&\hspace{.5cm} +4^{2p-2}\left(\int_0^u\Big[\E\Vert J_{0,\ell}^{-1}\Vert^{8p}\Vert\D\tilde{a}(\Phi(\ell,\om,\tilde{y}))\Vert^{4p}\Big]d\ell\right)^{1/2}\left(\int_0^u\Big[\E\Vert \D_r^i\Phi(\ell,\om,\tilde{y})\Vert^{4p}\Big]\right)^{1/2}d\ell
 \end{align*}
 In a similar fashion as in the first and second assertions, there exist constants $\tilde{C}_j(N,p,L_{b\sigma})<\infty, $ $j=1,2.,$ depending on $N,\; p,$ the coefficients $b,\sigma$ and their derivatives up to second order, such that  
 \begin{align*}
 &\hspace{-1cm}\E\Vert \D_r^i\hat{C}_u\Vert^{2p}\leqslant \exp(-2pL_{b\sigma}r)\bigg\{C_1(N,p,L_{b\sigma})\left[\frac{1}{12pL_{b\sigma}}\exp(12pL_{b\sigma}u)-\frac{1}{4pL_{b\sigma}}\exp(4pL_{b\sigma}u)\right]^{1/2}\\ &\left[\exp(8pL_{B\sigma}u)-1\right]^{1/2}
 +C_2(N,p, L_{b\sigma})\left[\exp(4pL_{b\sigma}u)-1\right]^{1/2}\left[\exp(4pL_{b\sigma}u)-4pL_{b\sigma}u-1\right]^{1/2}\bigg\}
 \end{align*}
 
  \qed
 
 \begin{lem}\label{Malliavin est}\rm
 Let $b^x, \sigma^x_{.k}\in \mathcal{C}_b^{l,\delta}(\R\times\R^N; \R^N)\; 1\leqslant k\leqslant m$ for $l\geqslant 2,$ be $\tau$-periodic in time. If the assumptions \ref{A2.1} and \ref{A2.1bc}~are satisfied. Then, there exists $0<K^1_t, K_t^2<\infty$ such that 
\begin{align}
\begin{cases}\E\left(\int_0^{t}\vert h_\eta(u)\vert^2du\right)\leqslant K_t^1\vert \eta\vert^{2},\\
\E\left(\int_0^{t}\int_0^{t}\text{trace}[\D_uh_\eta(r)\D_rh_\eta(u)]dudr\right)\leqslant K_t^1\vert \eta\vert^{2}.
\end{cases}
\end{align}
for all $\eta\in \Sc\times\R^N,$ where  $h_\eta(u):=\mathcal{A}_u^T\hat{C}_u^{-1}J_{0,u}\eta.$

 \end{lem}
 \noindent {\it Proof.} 
  To derive the first estimate, we notice that  
 \begin{align}\label{p0}
 \vert h_\eta(u)\vert^2\leqslant \Vert \mathcal{A}_u^{T}\Vert^2\Vert \hat{C}^{-1}_u\Vert^2\Vert J_{0,t}\Vert^2\vert \eta\vert^2
 \end{align}
 The restricted H\'omander condition \ref{A2.1bc}, implies there exists $M>0$ such that (e.g.,~\cite{Hairer10, Nulart}),
 \begin{align}\label{p1}
\sup_{u\in [0, t]} \Vert \hat{C}_u^{-1}\Vert^2 <M, 
\end{align}   
and since $\sigma_{.k}\in \mathcal{C}_b^l(\R\times\R^N;\R^N), \; 1\leqslant k\leqslant N, \; l\geqslant 2$, we have 
\begin{align}\label{pq}
\Vert \mathcal{A}^T\Vert^2\leqslant L^2_{\sigma}\Vert J_{0,u}^{-1}\Vert^2.
\end{align}

Finally, substituting (\ref{p1}),  \ref{pq}) and (\ref{p2}) into (\ref{p0}), then for all $t>0$, we have 
\begin{align}\label{p3}
\notag \E\left(\int_0^t\vert h_\eta(u)\vert^2du\right)&\leqslant (N+1)^2M L_{\sigma}^2\frac{1}{2L_{b\sigma}}\left(\exp(2L_{b,\sigma}t)-1\right)\vert \eta\vert^2\\
&:=K_t^1\vert \eta\vert^2, \qquad \eta\in \Sc\times\Rd.
\end{align}
 Now, we turn to the second estimate, we note that 
 \begin{align}\label{DerivativeES}
 \notag \E\left(\int_0^t\int_0^t\text{trace}[\D_uh_\eta(r)\D_rh_\eta(u)]drdu
\right) &= \E\left(\int_0^t\int_0^t\Vert \D_rh_\eta(u)\Vert^2_{\text{HS}}dr du\right)\\
&= \E\left(\sum_{i=1}^m\int_0^t\int_0^t\vert \D^i_rh_\eta(u)|^2dudr\right),
 \end{align}
 and
 \begin{align}\label{Ibound}
\notag \E\vert \D^i_rh_\eta(u)\Vert^2 &= \E\vert \D^i_r[(\mathcal{A}^T\hat{C}^{-1} J_{0,u}\eta]\vert^2 \\ & \notag \leqslant  4\E\vert\D^i_r[\mathcal{A}^T]\hat{C}^{-1}J_{0,u}\eta\vert^2+4\E\vert\mathcal{A}^T\hat{C}^{-1}\D_r^i[J_{0,u}]\eta\vert^2 +4\E\vert \mathcal{A}^T\D^i_r[\hat{C}^{-1}] J_{0,u}\eta\vert^2\\
 &=: 4\E\vert\Ic_1\vert^2+4\E\vert\Ic_2\vert^2+4\E\vert\Ic_3\vert^2.
 \end{align}
 We derive bound on  $\E\vert\Ic_1\vert^2$ as follows
   \begin{align}\label{Ma36}
 \notag  \E[ \vert \Ic_1\vert^2]& \leqslant \E\left[\Vert \D_u^i\mathcal{A}_t\Vert^2\Vert \hat{C}_t^{-1}\Vert^2\Vert J_{0,t}\Vert^2\vert \eta\vert^2\right] \leqslant M\E[\Vert \D_r^i\mathcal{A}_t\Vert^4]^{1/2}\E[\Vert J_{0,t}\Vert^4]^{1/2}\vert \eta\vert^2\\ &\hspace{2cm} \leqslant M(N+1)^2\exp(2L_{b\sigma}t)[\E\Vert \D_r^i\mathcal{A}_t\Vert^4]^{1/2}\vert \eta\vert^2.
   \end{align}
 It remains to find bound for $\E\Vert \D_r^i\mathcal{A}_t\Vert^4.$ To this end, we have,
 \begin{align}\label{Ma37}
 \E\left(\Vert \D_u^i\mathcal{A}_t\Vert^4\right)\leqslant 8L^4_{\sigma}\int_0^t\E\Vert \D^i_rJ_{0,u}^{-1}\Vert^4 du +8L^4_{\sigma}\int_0^t\E\Vert \D_r^i\Phi(u,\om,\tilde{x})J_{0,u}^{-1}\Vert^4du
 \end{align}
  Recall that the Fr\'echet derivative of a square matrix $A$ in the direction of $H$ is given by $D_HA^{-1} = -A^{-1}HA^{-1}$, so that by chain rule, we have 
  \begin{align*}
 \int_0^t\E\Vert \D^i_rJ_{0,u}^{-1}\Vert^4 du &=  \int_0^t\E\Vert J_{0,u}^{-1}\D_r^i[J_{0,u}]J_{0,u}^{-1}\Vert^4du \\[.2cm]
& \leqslant \int_0^t[\E\Vert J_{0,u}^{-1}\Vert^8]^{1/2}[\E\Vert \D_r^iJ_{0,u}\Vert^4]^{1/2}du\\[.2cm]
&\leqslant (N+1)^2\int_0^t\exp(2L_{b\sigma}u)[\E\Vert \D_r^iJ_{0,u}\Vert^4]^{1/2}du\\[.2cm]
&\leqslant (N+1)^2\left(\int_0^t\exp(4L_{b\sigma}u)du\right)^{1/2}\left(\int_0^t\E\Vert \D_r^iJ_{0,u}\Vert^4du\right)^{1/2}\\[.2cm]
&\leqslant \frac{(N+1)^2}{2L_{b\sigma}^{1/2}}\left(\exp(4L_{b\sigma}t)-1\right)^{1/2}\left(\int_0^t\E\Vert \D_r^iJ_{0,u}\Vert^4du\right)^{1/2},
\end{align*}
and, 
\begin{align*}
  \int_0^t\E\Vert \D_r^i\Phi(u,\om,\tilde{y}) J_{0,u}^{-1}\Vert^4du & \leqslant \frac{(N+1)^2}{2L^{1/2}_{b,\sigma}}\Big[\exp(4L_{b\sigma}t) -1\Big]^{1/2}\E\left(\int_{0}^t\Vert \D_r^i\Phi(u,\om,\tilde{y})\Vert^8du\right)^{1/2}.
  \end{align*}
   Now, we apply Lemma \ref{Malliavin1est} to obtain
   \begin{align}\label{Ma38}
  \notag &\hspace{-1cm} \int_0^t\E\Vert \D^i_rJ_{0,u}^{-1}\Vert^4 du =  \int_0^t\E\Vert \D_r^i\Phi(u,\om,\tilde{y}) J_{0,u}^{-1}\Vert^4du\\ & \leqslant \frac{\tilde{L}_{b\sigma}^4(N+1)^6}{2L_{b\sigma}}\exp(-\frac{1}{2}L_{b\sigma}r)\Big[\frac{1}{3}\exp(12L_{b\sigma}t)-\exp(4L_{b\sigma}t)\Big]^{1/2}\Big[\exp(4L_{b\sigma}t)-1\Big]^{1/2}
   \end{align}
   and,
   \begin{align}\label{Ma39}
 \int_0^t\E\Vert \D_r^i\Phi(u,\om,\tilde{y}) J_{0,u}^{-1}\Vert^4du \leqslant \frac{(N+1)^4}{2\sqrt{2}L_{b,\sigma}}\Big[\exp(4L_{b\sigma}t) -1\Big]^{1/2}\exp(-L_{b\sigma}r)\Big[\exp(2L_{b\sigma}t)-1\Big]^{1/2}
   \end{align}
   Substituting (\ref{Ma38}) and (\ref{Ma39}) into (\ref{Ma37}), we have 
   \begin{align}\label{Ma40}
\notag   \E\left(\Vert \D_u^i\mathcal{A}_t\Vert^4\right) &\leqslant 
4\tilde{L}_{\sigma}^4L^3_{b\sigma}(N+1)^6\exp(-\frac{1}{2}L_{b\sigma}r)\Big[\exp(12L_{b\sigma}t)-\exp(4L_{b\sigma}t)\Big]^{1/2}\Big[\exp(4L_{b\sigma}t)-1\Big]^{1/2}\\
&\hspace{.5cm} +\frac{4(N+1)^4}{\sqrt{2}L_{b,\sigma}}\exp(-\frac{1}{2}L_{b\sigma}r)\Big[\exp(4L_{b\sigma}t) -1\Big]^{1/2}\Big[\exp(2L_{b\sigma}t)-1\Big]^{1/2}
   \end{align}
   Put the inequality (\ref{Ma40}) in (\ref{Ma36}), we have 
   \begin{align}\label{Ma41}
\notag  & \E\vert \Ic_1\vert^2\leqslant 2M(N+1)^4 \exp(-\frac{1}{4}L_{b\sigma}r)\bigg((N+1)^2\tilde{L}^4_{b\sigma}L_{b\sigma}\Big[\exp(12L_{b\sigma}t)-\exp(4L_{b\sigma}t)\Big]^{1/2} \\ &\hspace{2cm}+\Big[\exp(2L_{b\sigma}t)-1\Big]^{1/2} \frac{1}{\sqrt{2}L_{b\sigma}}\bigg)^{1/2} \Big[\exp(16L_{b\sigma}t)-\exp(8L_{b\sigma}t)\Big]^{1/4}\vert \eta\vert^2
   \end{align}
   Next, we derive bound on $\E\vert \Ic_2\vert^2$ as follows
   \begin{align}\label{Ma42}
  \notag \E\vert \Ic_2\vert^2 &= \E\vert \mathcal{A}^T\hat{C}^{-1}\D_r^i[J_{0,u}]J_{0,u}\eta\vert^2\\ & \leqslant\E\Big[\Vert \mathcal{A}\Vert^2\Vert \hat{C}^{-1}\Vert^2\Vert \D_r^iJ_{0,u}\Vert^2\Vert J_{0,u}\Vert^2\Big]\vert \eta\vert^2 \\ \notag &\leqslant M^2 L_\sigma^2\E\Big[\Vert J_{0,u}^{-1}\Vert^2\Vert J_{0,u}\Vert^2\Vert D_r^iJ_{0,u}\Vert^2\Big]|\eta|^2\\ \notag &\leqslant M^2L_{\sigma}^2\Big[\E\Vert J_{0,u}^{-1}\Vert ^4\Vert J_{0,u}\Vert^4\Big]^{1/2}\Big[\E\Vert \D_r^iJ_{0,u}\Vert^4\Big]^{1/2}|\eta|^2\notag\\
 & \leqslant 2M^2(N+1)^2\frac{\tilde{L}^2_{b\sigma} (N+1)^3}{L_{b\sigma}^{1/2}}\exp(L_{b\sigma}(2u-r))\Big[\exp(2L_{b\sigma}u)-1\Big]^{1/2}|\eta|^2
   \end{align}
   Finally, for $\E\vert \Ic_3\vert^2$, again, we use that the Fr\'echet derivative of a square matrix $A$ in the direction of $H$ is given by $D_HA^{-1} = -A^{-1}HA^{-1}$, we have 
 \begin{align*}
\D^i_r[\hat{C}^{-1}]= -\hat{C}^{-1}[\D^i_r\hat{C}]\hat{C}^{-1} 
 \end{align*}
 We obtain the following estimate,
 \begin{align*}
 \E\vert \Ic_3\vert^2 &= \E\vert \mathcal{A}\hat{C}^{-1}\D_r^i[\hat{C}_u]\hat{C}^{-1}J_{0,u}\eta\vert^2\\
& \leqslant \Big[\E\Vert \mathcal{A}\Vert^2\Vert \hat{C}^{-1}\Vert^2 \Vert\D_r^i\hat{C}_u\Vert^2\Vert \hat{C}^{-1}\Vert^2\Vert J_{0,u}\Vert^2\vert \eta\vert^2\Big]|\eta|^2\\
& \leqslant \Big[\E\Vert \mathcal{A}_u\Vert^4 \Vert \hat{C}_u^{-1}\Vert^{8}\Vert J_{0,u}\Vert^4\Big]^{1/2}\Big[\E\Vert \D_r^i\hat{C}_u\Vert^4\Big]^{1/2}|\eta|^2\\
& \leqslant M^4\Big[\E\Vert \mathcal{A}_u\Vert^8\Big]^{1/4} \Big[\E\Vert J_{0,u}\Vert^8\Big]^{1/4} \Big[\E\Vert \D_r^i\hat{C}_u\Vert^4\Big]^{1/2}|\eta|^2\\
& \leqslant M^2L_{\sigma}^2 \Big[\E\Vert J_{0,u}^{-1}\Vert^8\Big]^{1/4} \Big[\E\Vert J_{0,u}\Vert^8\Big]^{1/4} \Big[\E\Vert \D_r^i\hat{C}_u\Vert^4\Big]^{1/2}|\eta|^2\\
& \leqslant M^2L_{\sigma}^2(N+1)^4\exp(4L_{b\sigma}u) \Big[\E\Vert \D_r^i\hat{C}_u\Vert^4\Big]^{1/2} |\eta|^2.
 \end{align*}
Then, by Lemma \ref{Malliavin1est}, we have constants $0<\tilde{C}_j(N,p,L_{b\sigma})<\infty, $ $j=1,2.,$ depending on $N,\; p =2,$ the coefficients $b,\sigma$ and their derivatives up to second order such that  
 \begin{align}\label{Ma43}
 \notag &\hspace{-1cm}\E\vert \Ic_3\vert^2\leqslant \exp(2L_{b\sigma}(2u-r))\bigg\{M^2C_1(N,2,L_{b\sigma})\left[\frac{1}{24L_{b\sigma}}\exp(24L_{b\sigma}u)-\frac{1}{8L_{b\sigma}}\exp(8L_{b\sigma}u)\right]^{1/2}\\ &\left[\exp(16L_{B\sigma}u)-1\right]^{1/2}
 +C_2(d,2, L_{b\sigma})\left[\exp(8L_{b\sigma}u)-1\right]^{1/2}\left[\exp(8L_{b\sigma}u)-8L_{b\sigma}u-1\right]^{1/2}\bigg\}^{1/2}|\eta|^2.
 \end{align}
 If we put the estimates (\ref{Ma41}), (\ref{Ma42}) and (\ref{Ma43}) into the inequality (\ref{Ibound}), and then subsitute into the double integral (\ref{DerivativeES}), we obtain the desired bound. \qed
 
 \medskip
 
 With the above preparatory estimates at hand, we derive the strong Feller property of the Markov semigroup $(\tilde{\mathcal{P}}_{t})_{t\geqslant 0};$ we shall skip the dependence on $x$ to simplify notation.
  \begin{prop}\label{Strong Feller}\rm
 Let $b^x, \sigma^x_{.k}\in \mathcal{C}_b^{l,\delta}(\R\times\R^N; \R^N)\; 1\leqslant k\leqslant N,$ for $l\geqslant 2,$ be $\tau$-periodic in time. If the assumptions \ref{A2.1} and \ref{A2.1bc}~are satisfied. Then, there exists $0<K_t<\infty$ such that for $\tilde{x}, \tilde{y}\in \Sc\times\Rd$ with $t\in [0, \tau),$ and $\varphi\in \mathcal{C}_b(\Sc\times\R^N)$, we have
\begin{align}
\vert \tilde{\mathcal{P}}_{t}\varphi(\tilde{y}) -\tilde{\mathcal{P}}_{t}\varphi(\tilde{z})\vert \leqslant {K_t}||\varphi||_{\infty}\vert \tilde{y}-\tilde{z}\vert.
\end{align}
 \end{prop}
 \noindent {\it Proof.}  
 First, recall the Malliavin integration by part formula (e.g.,~\cite{Hairer10, Hairer_Note16, Hairer11, Malliavin, Nulart, Arnaudon}) yield 
\begin{align}\label{Int_by}
D(\tilde{\mathcal{P}}_t\varphi)(\tilde{y})\eta= \E\left(\varphi(\Phi(t,\om,\tilde{y}))\int_0^th_\eta(u)\star d\tilde{W}_u\right), \quad \varphi \in \mathcal{C}_b^1(\Sc\times\R^N),
\end{align}
where $h_\eta(u):=\tilde{\sigma}^TJ_{u,t}^T\hat{C}_{t}^{-1}J_{0,t}\eta,$ for $\eta\in \Sc\times\Rd,$ and $J_{u,t} =  J_{0,t}\circ J_{0, u}^{-1}$. Note that $h_\eta$ need not be adapted to $\sigma\{W_t: t\geqslant 0\}$, however, one can interprete the stochastic integral (\ref{Int_by}) in the sense of Skorokhod (cf.~\cite{Hairer10}).

Next, since $\mathcal{C}_b^1(\Sc\times\R^N)$ is dense in $\mathcal{C}_b(\Sc\times\R^N)$, we have $(\varphi_n)_{n\in \N}\subset \mathcal{C}_b^1(\Sc\times\R^N)$ such that 
\begin{align}\label{Squence}
\begin{cases}\lim_{n\rightarrow\infty}\tilde{\mathcal{P}}_t\varphi_n(\tilde{y}) = \tilde{\mathcal{P}}_t\varphi(\tilde{y}), \\\
\lim_{n\rightarrow\infty}D(\tilde{\mathcal{P}}_t\varphi_n)(\tilde{y})\eta = \E\left(\varphi(\Phi(t,\om,\tilde{y}))\int_0^th_\eta(u)\star d\tilde{W}_u\right),
\end{cases}
\end{align}
for all $\varphi \in \mathcal{C}_b(\Sc\times\R^N),$ convergence being uniform (e.g.,~\cite{Daprato}). On the other hand, the regularity of the coefficients $\tilde{b}, \; \tilde{\sigma}$ and assumption \ref{A2.1bc} ensure the existence of a function $0<\tilde{\rho}_{t}\in \mathcal{C}_b^2(\Sc\times\R^N)\times\mathcal{C}_b^2(\Sc\times\R^N),$ such that $\tilde{P}(0,\tilde{y};t,d\tilde{z}) = \tilde{\rho}_{t}(\tilde{y},\tilde{z})d\tilde{z}$. This implies that 
\begin{align}\label{densitywise}
\notag \lim_{n\rightarrow\infty}D(\PT_t\varphi_n)(\tilde{y}) &= \lim_{n\rightarrow\infty}\int_{\Sc\times\R^N}\varphi_n(\tilde{y})D_{\tilde{y}}\tilde{\rho}_t(\tilde{y},\tilde{z})d\tilde{z} \\ &= \int_{\Sc\times\R^N}\varphi(\tilde{y})D_{\tilde{y}}\tilde{\rho}_t(\tilde{y},\tilde{z})\eta d\tilde{z} = D(\PT_t\varphi)(y)\eta .
\end{align}
The equalities (\ref{Squence}) and (\ref{densitywise}), yield (\ref{Int_by}) for all $\varphi\in \mathcal{C}_b(\Sc\times\R^N).$ By Cauchy Schwartz inequality, we have 
\begin{align}\label{Cauch}
\vert D(\tilde{\mathcal{P}}_t\varphi)(\tilde{x})\eta\vert\leqslant \sqrt{(\tilde{\mathcal{P}}_t\varphi^2)(\tilde{x})}\left(\E\left\vert\int_0^t h_\eta(u) \star dW_u\right\vert^2\right)^{1/2}, \quad \varphi\in \mathcal{C}_b(\Sc\times\R^N)
\end{align}
 
 Then by generalised It\^o isometry (cf.~\cite{Hairer10, Nulart}) and Lemma \ref{Malliavin est}, we have 
 \begin{align*}
 \E\left\vert\int_0^t h_\eta(u) \star dW_u\right\vert^2 &= \E\left(\int_0^t\vert h_\eta(u)\vert^2du\right)+\E\left(\int_0^t\int_0^t\text{trace}[\D_uh_\eta(r)\D_rh_\eta(u)]dudr\right) \\ &\hspace{2cm}\leqslant \tilde{K}_t\vert \eta\vert^2
 \end{align*}
 It follows from  (\ref{Cauch}), that 
 \begin{align*}
 \vert D(\tilde{\mathcal{P}}_t\varphi)(\tilde{y})\eta\vert\leqslant K_t\Vert \varphi\Vert_{\infty}\vert \eta\vert , \qquad \tilde{y},\eta\in \Sc\times\R^N.
 \end{align*}
 Finally, let $\tilde{\xi}_{\ell} = \ell \tilde{y}+ (1-\ell)\tilde{x}, \; \tilde{y}, \tilde{z}\in \Sc\times\R^N, \; \ell\in [0, 1]$ and $\eta = \tilde{y}-\tilde{z}$, by the mean value theorem, we have 
 \begin{align*}
 \vert \tilde{\mathcal{P}}_t\varphi(\tilde{y})- \tilde{\mathcal{P}}_t\varphi(\tilde{z})\vert \leqslant \int_0^1\vert D(\tilde{\mathcal{P}}_t\varphi)(\tilde{\xi}_\ell)\eta\vert d\ell \leqslant K_t\Vert\varphi\Vert_{\infty}\vert \tilde{y}-\tilde{z}\vert  
 \end{align*}
 \qed

 \begin{theorem}\label{Ps_erg}\rm\;
Suppose that the coefficients of the SDE (\ref{Fast}) are $\tau$-periodic in time and satisfy Assumptions \ref{A2.1} and \ref{A2.1bc}. Suppose further that 
\begin{align}\label{Extra}
\sup_{t\in [0, 1]}\E\left[V(\tilde{y} - \Phi(t,\om, \tilde{y}))\right]<\infty,
\end{align}
then, for fix $x\in \Nb,$ the family of periodic measures $(\tilde{\mu}^x_r)_{r\in [0,\tau)}$ induced by the random periodic path $\tilde{S}^x(r,\om)$ is PS-ergodic.
\end{theorem}
\noindent {\it Proof.}  The idea of the proof is to use the construction carried out in Theorem \ref{R_per1} to verify the conditions of Lemma \ref{Expon}. To this end, we skip the dependence on $x$ to simplify notation and divide the proof into three steps.\\
{\bf Step I:} We show that for any $\tilde{Y} ,\tilde{Z}\in L^p(\Om, \F_{-\infty}^r, \p),\;p\geqslant 1$ and for $k\in \Z,$ the following inequality holds
\begin{align}\label{Est_gen}
\notag &\E\vert \Phi(r+k\tau,\om, \tilde{Y}(\om))- \Phi(r+k\tau, \om, \tilde{Z}(\om))\vert^p\\ &\hspace{3cm}\leqslant C\exp\Big(\int_{-r}^{k\tau}\lambda(u+r)du\Big)\E\vert \tilde{Y}(\om)-\tilde{Z}(\om)\vert^p.
\end{align}
To derive this, set $\alpha(t) = \exp\Big(-\int_0^t\lambda(r)ds\Big)$ and $$M(t,\om,\tilde{y})= M(t,\om, (s,y))=\int_0^t\begin{pmatrix}
 0& 0\\ 0 &\sum_{k=1}^N\sigma_k(r+s,x, y)
  \end{pmatrix}\begin{pmatrix} dW_r^0\\ dW_r^k\end{pmatrix},$$ then, by It\^o's formula (\cite{Hkunita, Kunita, Schmal01}), we obtain 
\begin{align*}
&\hspace{-.5cm} d\bigg(\alpha(t)V(\Phi(t,\om,\tilde{Y})-\Phi(t,\om,\tilde{Z}))\bigg)\\[.1cm] &\hspace{.5cm}=-\lambda(t)\alpha(t) V(\Phi(t,\om,\tilde{Y})-\Phi(t,\om,\tilde{Z}))dt +\alpha(t)\mathcal{L}^{(2)}V(\Phi(t,\om,\tilde{Y})-\Phi(t,\om,\tilde{Z}))dt\\[.1cm] &\hspace{1.5cm}
+\alpha(t)V(\Phi(t,\om,\tilde{Y})-\Phi(t,\om,\tilde{Z}))d\Big(M(\Phi(t,\om,\tilde{Y})-M(\Phi(t,\om,\tilde{Z}))\Big)\\[.1cm] & \hspace{2.5cm}
\leqslant   -\lambda(t)\alpha(t) V(\Phi(t,\om,\tilde{Y})-\Phi(t,\om,\tilde{Z}))dt
  +\lambda(t)\alpha(t) V(\Phi(t,\om,\tilde{Y})-\tilde{X}(t,\om,\tilde{Z}))dt\\[.1cm] & \hspace{3.5cm}+\alpha(t)V(\Phi(t,\om,\tilde{Y})-\Phi(t,\om,\tilde{Z}))d\Big(M(\Phi(t,\om,\tilde{Y})-M(\Phi(t,\om,\tilde{Z}))\Big).
\end{align*}
This implies for $k\in \N,$ 
\begin{align*}
\E\bigg(V(\Phi(r+k\tau,\om,\tilde{Y}(\om)) -\Phi(r+k\tau,\om,\tilde{Z}(\om))\bigg)\leqslant \exp\Big(\int_{-r}^{k\tau}\lambda(u+r)du\Big)\E[V(\tilde{Y}(\om)-\tilde{Z}(\om))]
\end{align*}
We use that $\vert y\vert^p \leqslant V(t,y), \; p\geqslant 1,$ to obtain for $k\in\N$,
\begin{align*}
\E\bigg(\vert \Phi(r+k\tau,\om,\tilde{Y}(\om))- \Phi(r+k\tau,\om,\tilde{Z}(\om))\vert^p \bigg)&\leqslant \E\bigg(V(\Phi(r+k\tau,\om,\tilde{Y}(\om)) -\Phi(r+k\tau,\om,\tilde{Y}(\om))\bigg)\\
&\leqslant C\exp\Big(\int_{-r}^{k\tau}\lambda(u+r)du\Big)\E\vert \tilde{Y}(\om)-\tilde{Z}(s,\om)\vert^p.
\end{align*}
{\bf Step II:} In this step, we show that there exists $0<\tilde{C}<\infty$ such that for any $\tilde{Y}\in L_p(\Om, \F_{-\infty}^r, \p),$ the following inequality holds
 \begin{align}\label{Lim}
\E\vert \Phi(r+k\tau, \om, \tilde{Y}(\om))-\tilde{S}(r, \theta_{k\tau}\om)\vert^p\leqslant\tilde{C}H_r(k\tau)
\end{align}
where  $H_r(k\tau) = \eta_r(k\tau)+\exp\Big( \frac{1}{2}\beta k\tau\Big)$ and 
\begin{align*}
\eta_r(k\tau) = \exp\left(\int_{-r}^{k\tau}\lambda(u+r)du\right)\sum_{j=1}^{k-1}\exp\left(\int_{0}^{j\tau}\lambda(u)du\right) , \qquad k\in\N.
\end{align*}
To see this, recall a construction done in Theorem \ref{R_per1}, precisely, equation (\ref{PS_L}) reads 
\begin{align*}
\tilde{S}(r,\om) = \lim_{k\rightarrow\infty}\Phi(r+k\tau,\theta_{-k\tau}\om,\tilde{y}), \quad r\in [0, \tau].
\end{align*}
As $\tilde{S}(r+k\tau,\om) = \tilde{S}(r,\theta_{k\tau}\om),$ then by equation (\ref{Est_gen}) in step I, estimate (\ref{Expfas}) and triagngle inequality, we have 
\begin{align*}
&\E\Big(\vert \Phi(r+k\tau,\om,\tilde{y})-\tilde{S}(r,\theta_{k\tau}\om)\vert^p\Big)\\ & \hspace{2cm}\leqslant \E\Big(\vert\Phi(r+k\tau, \om,\tilde{y})-\Phi(k\tau,\om,\Phi(r+k\tau, \theta_{-k\tau}\om,\tilde{y}))\vert^p\Big)+ \hat{\gamma}(\om)\exp\Big( \frac{1}{2}\beta k\tau\Big)\\
&\hspace{2cm}= \E\left(\vert\Phi(r+k\tau,\om,\tilde{y})-\Phi(r+2k\tau,\theta_{-k\tau}\om,\tilde{y})\vert^p\right)+  \hat{\gamma}(\om)\exp\Big( \frac{1}{2}\beta k\tau\Big) \\
&\hspace{2cm}= \E\left(\vert \Phi(r+k\tau,\om,\tilde{y})-\Phi(r+k\tau,\om,\Phi(k\tau,\theta_{-k\tau}\om,\tilde{y}))\vert^p\right)+\hat{\gamma}(\om)\exp\Big( \frac{1}{2}\beta k\tau\Big)\\
&\hspace{2cm}\leqslant C\exp\left(\int_{-r}^{k\tau}\lambda(u+r)du\right)\E\vert \tilde{y}-\Phi(k\tau,\theta_{-k\tau}\om,\tilde{y})\vert^p +  \hat{\gamma}(\om)\exp\Big( \frac{1}{2}\beta k\tau\Big),
\end{align*}
where $\hat{\gamma}$ and $\beta$ are random variables appearing in (\ref{Expfas}).
Now, set $\beta_r(t) := \exp\left(\int_{-r}^{t}\lambda(u+r)du\right),$ by the cocycle property of $\Phi$ and step I again, we have  
\begin{align*}
 &C\beta_r(k\tau)\E\vert \tilde{y}-\Phi(k\tau, \theta_{-k\tau}\om, \tilde{y})\vert^p\\&\hspace{.5cm}  \leqslant C\beta_r(k\tau)\bigg\{\E\bigg[\vert \tilde{y}-\Phi(1, \theta_{-1}\om, \tilde{y})\vert^p+\vert \Phi(1, \theta_{-1}\om,\tilde{y})-\Phi(2, \theta_{-2}\om, \tilde{y})\vert^p \\ &\hspace{2cm}  +\cdots + \vert \Phi([k\tau]-1, \theta_{1-[k\tau]}\om, \tilde{y})-\Phi([k\tau],\theta_{-[k\tau]}\om, \tilde{y})\vert^p\bigg]\bigg\}\\
 &= C\beta_r(k\tau)\bigg\{\E\bigg[\vert \tilde{y}-\Phi(1, \theta_{-1}\om, \tilde{y})\vert^p+\vert \Phi(1, \theta_{-1}\om,\tilde{y})-\Phi(1, \theta_{-1}\om, \Phi(1, \theta_{-1}\om, \tilde{y}))\vert^p \\ &\hspace{1cm}+\cdots + \vert \Phi([k\tau]-1, \theta_{1-[k\tau]}\om, \tilde{y})-\Phi([k\tau]-1,\theta_{1-[k\tau]}\om, \Phi(1, \theta_{-1}\om,\tilde{y}))\vert^p\bigg]\bigg\}\\
&\hspace{.5cm}\leqslant C\beta_r(k\tau)\bigg\{ \Big[C\beta_0(1)+\cdots +C\beta_0([k\tau]-1)\Big]\E\vert \tilde{y}-\Phi(1, \theta_{-1}\om,\tilde{y})\vert^p\bigg\}\\
&\hspace{.5cm}=C^2\beta_r(k\tau)\sum_{j=1}^{[k\tau]-1}\beta_0(j)\E\vert \tilde{y}-\Phi(1, \theta_{-1}\om,\tilde{y})\vert^p.
\end{align*}
Consequently, we have 
\begin{align*}
\E\Big(\vert \Phi(r+k\tau,\om,\tilde{y})-\tilde{S}(r,\theta_{\tau}\om)\vert^p\Big)\leqslant C^2\eta_r(k\tau)\E\vert \tilde{y}-\Phi(1,\theta_{-1}\om,\tilde{y})\vert^p+ \hat{\gamma}(\om)\exp\Big( \frac{1}{2}\beta k\tau\Big)
\end{align*}
and using that $\vert \tilde{y}\vert^p\leqslant V(\tilde{y})$, we obtain
\begin{align*}
\E\Big(\vert \Phi(k\tau,\om,\tilde{y})-\tilde{S}(s,\theta_{k\tau}\om)\vert^p\Big)\leqslant C^2\eta_r(k\tau)\Big(\E V( \tilde{y}-\Phi(1, \theta_{-1}\om,\tilde{y}))\Big) + \hat{\gamma}(\om)\exp\Big( \frac{1}{2}\beta k\tau\Big).
\end{align*}
 Finally, by the condition (\ref{Extra}) together with the $\p$-preserving property of $\theta$, we arrive at the required estimate (\ref{Lim}).\\ 
 {\bf Step III:}
 In this step, we show that the conditions of Lemma \ref{Expon} are satisfied. For this purpose, we use step II and invariance of $\tilde{\mu}_r$ w.r.t. the discrete Markov semigroup $(\tilde{\mathcal{P}}_{r+k\tau})_{k\in \N}$ and for any $h\in \text{Lip}_b(\Bt_r),\; r\in [0,\tau),$ to obtain
\begin{align*}
\left\vert \tilde{\mathcal{P}}_{r+k\tau}h-\int_{\Bt_r}hd\tilde{\mu}_r\right\vert&\leqslant \int_{\Sc\times\Yt}\vert \E[h(\Phi(r+k\tau,\om,.)- h(\tilde{S}(r+k\tau,\om))\vert d\tilde{\mu}_0\\
&\hspace{.5cm}\leqslant \text{Lip}(h)\int_{\Sc\times\Yt}\E\vert \Phi(r+k\tau, \om, .)-\tilde{S}(r+k\tau,\om)\vert d\tilde{\mu} \\
&\hspace{1cm}\leqslant \hat{C}\text{Lip}(h)H_r(k\tau),
\end{align*}
i.e., there exists $0<\tilde{C}<\infty$ such that 
\begin{align}\label{LipSCH_CQ}
\left\vert \tilde{\mathcal{P}}_{k\tau+r}h(\tilde{x})-\langle h,\tilde{\mu}_r\rangle\right\vert\leqslant \tilde{C} \text{Lip}(h)\bigg[\exp\left(p^{-1}\int_0^{k\tau+r}\lambda(u)du\right) +\exp\Big( \frac{1}{2p}\beta k\tau\Big)\bigg].
\end{align}
Finally, let $f\in \mathcal{C}_b(\Bt_r)$ be given. Setting $h = \PT_{r+\tau}f $ in (\ref{LipSCH_CQ}), we have by Proposition \ref{Strong Feller}
and by the invariance of $\tilde{\mu}_t$ under the Markov semigroup $(\PT_{r+n\tau})_{k\in \N},$ we have 
\begin{align*}
\left\vert \PT_{r+\tau+k\tau}f(\tilde{y})- \langle\PT_{r+\tau}f, \tilde{\mu}_r\rangle\right\vert = \left\vert \int_{\Bt_r}\left(\PT_{r+\tau+ k\tau}f(\tilde{y}) -\PT_{t+k\tau +\tau}f(\tilde{y})\right)\tilde{\mu}_r(d\tilde{y})\right\vert\\
\leqslant \tilde{C}\text{Lip}
(\PT_{r+\tau}f)\bigg[\exp\left(p^{-1}\int_0^{r+k\tau}\lambda(u)du\right)+\exp\Big( \frac{1}{2p}\beta k\tau\Big)\bigg].\\
\leqslant \tilde{K}_{\tau}\Vert f \Vert_{\infty}\bigg[\exp\left(p^{-1}\int_0^{r+k\tau}\lambda(u)du\right)+\exp\Big( \frac{1}{2p}\beta k\tau\Big)\bigg].
\end{align*}
where $\tilde{K}_\tau= K_\tau \tilde{C}$.
It then follows that for any $f\in \mathcal{C}_b(\Bt_r)$ and $k>1,$ 
\begin{align*}
\left\vert \PT_{t+k\tau}f(\tilde{y})- \int_{\Bt_r}f d\tilde{\mu}_r\right\vert\leqslant \tilde{K}_{\tau}\Vert f\Vert_{\infty}\bigg[\exp\left(p^{-1}\int_{0}^{r+k\tau-\tau}\lambda(u)du\right)+\exp\Big( \frac{1}{2p}\beta k\tau\Big)\bigg],
\end{align*}
this implies that there $0<\tilde{K}<\infty$ such that (\ref{C2}) holds.
Then, by Lemma \ref{Expon}, we have the convergence of Krylov--Bogolyubov scheme for periodic measures.

\qed

 \begin{rem}\rm\;
 The approach adopted in this subsection relies on the construct periodic measure from the random periodic path $\tilde{S}(r,\om)$. It is also possible to construct random periodic paths given periodic measures. To see this, we proceed as in \cite{Feng18} and suppose that we are given the periodic measures $(\tilde{\mu}_r)_{r\in\R}\subset \mathcal{P}(\Sc\times\Yt)$ of period $\tau>0,$ for a Markovian RDS $\Phi$ on $\Sc\times\Yt.$  Let the skew product  flow $\chi: \Om\times\Sc\times\Yt \rightarrow\Om\times\Sc\times\Yt$ be defined by 
\begin{align}\label{EnlD}
\chi_t(\om, \tilde{y})&= (\theta_t\om, \Phi(t,\om,\tilde{y})),\quad t\in\R.
\end{align}
Note that the two sided skew product in (\ref{EnlD}) is defined by running two independent Brownian motions $\{W^k_t: t\geqslant 0, \; 1\leqslant k\leqslant N\}$ and $\{W^k_{-t}: t\geqslant 0, \; 1\leqslant k\leqslant N\}$ (e.g.,~\cite{Arnold, Kunita}). Set
\begin{align*}
&\hat{\om}=(\om,\tilde{y}),\quad \Hat{\Om} := \Om\times\Sc\times\Yt, \quad \hat{\F} :=\F\otimes\Bb(\Sc\times\Yt)\\
&\bar{\tilde{\mu}}(\tilde{A}) :=\frac{1}{\tau}\int_0^{\tau}\tilde{\mu}_r(\tilde{A})dr,\quad 
\Psi(t,\hat{\om}, \cdot) := \Phi(t,\om, \cdot), \quad t\in \R.
\end{align*}
It is relatively straightforward to check that $\Psi$ defines an RDS over the ergodic base dynamical system $\bar{\Om}=(\Hat{\Om}, \hat{\F}, \bar{\tilde{\mu}}, (\chi)_{t\in \R}).$
A random periodic path $R:\Rp\times\hat{\Om}\rightarrow\Sc\times\Yt$ for the enlarged RDS $\Psi$ is constructed as follows: for any $\hat{\om}^{*} = (\om^*, \tilde{y}^*(\theta_{-r}\om^*))\in \hat{\Om},\; r\in \R$ 
\begin{align}
R(t,\om^*) := \Phi(t+r,\theta_{-r}\om^*, \tilde{y}^*(\theta_{-r}\om^*)), \quad t\in \Rp.
\end{align}
Finally, the transition probability of $\Psi(t,\hat{\om},\tilde{y})$ is the same as $\tilde{P}(0,\tilde{y};t,.)$ and the law of $R$ is 
\begin{align*}
\hat{\p}\circ R^{-1}(t,.) =\tilde{ \mu}_t, \quad \text{for all $t$}\in\Rp.
\end{align*} 
Note the resemblance of the process $R(t,\om^*)$ and the limit $\tilde{S}(t,\om)$ of the Cauchy sequence in equation (\ref{PS_L}) from Theorem \ref{R_per1}, the major difference is that the construction involve different ergodic base dynamical system. Related idea was previously employed in \cite{Uda16} to construct stable random periodic solutions for a Markovian RDS on a cylinder $\Sc\times\R^N$. 
\end{rem}

\section{Stochastic averaging in the random periodic regime}\label{Aver_pp}
One of the difficulties in averaging for non-autonomous SDEs is to give a suitable notion of ergodic evolution system of probability measures. This issue has been resolved in subsection \S\ref{Ranp_M} for some class of SDEs with time periodic forcing. Having established the ergodicity of periodic measures for a class of systems satisfying the Assumptions \ref{A2.1} and \ref{A2.1bc}, the proof of averaging principle will follow the same line of argument as in (e.g.,~\cite{ Fredlinbook, Hasm66, Kifer09, Loch88}) once we know the manner in which the PS-ergodic measures depend on the frozen slow variables $x\in \Nb\Subset\Rd$. First, we deal with dependency of the PS-ergodic measures $\bar{\tilde{\mu}}^x$ on the parameter $x\in \Nb\Subset \Rd$ in \S\ref{Ran_Lp} and subsequently, identify the averaged equation.  We prove an averaging principle in the random periodic regime in \S\ref{AV_limit} and conclude with a toy example arising from kinetic theories of turbulent flows (e.g.,~\cite{Bouchet15}).
\subsection{The averaged equation }\label{Ran_Lp}
A crucial step in averaging principle is to identify the averaged equation via the vector field of the slow motion and the ergodic measures of the fast subsystem. 
As the averaging framework considered in this paper is fully coupled, some degree of regularity of the invariant measures $x\mapsto\bar{\tilde{\mu}}^x$ will be required. In fact, Lipschitz continuity with respect to the parameter $x\in\Nb\Subset \Rd$ is enough for our purpose. This Lipschitz property together with the regularity of the slow vector field $F$ would ensure the unique local solvability of the averaged equation. We employ similar argument as the one in  \S\ref{Ranp_M} to prove Lipschitz dependence of periodic measures on parameters.

\medskip
We recall the notion of Lipschitz continuity of probability measures in the sense of the narrow topology on $\mathcal{P}(\Sc\times\Yt).$ Consider a real-valued function $f: \Sc\times\Yt\rightarrow\R$ and the space of $1$-Lipschitz functions defined by 
\begin{align*}
\text{Lip}(1,\Sc\times\Yt) =\{f\in \text{Lip}_b(\Sc\times\Yt): \Vert f\Vert_{BL}\leqslant 1\}.
\end{align*}
The space of $1$-Lipschitz functions $\text{Lip}(1,\Sc\times\Yt)$ generate a narrow topology on $\mathcal{P}(\Sc\times\Yt)$ (e.g.,~\cite{Crau2}) with the metric
\begin{align}
\text{d}_{BL}(\mu,\nu) = \sup\bigg\{\int_{\Sc\times\Yt}fd(\mu-\nu): f\in \text{Lip}(1,\Sc\times\Yt)\bigg\}.
\end{align}
\medskip
\begin{definition}[Lipschitz continuity of probability measures]\rm\;
 We say that a probability measure $(x_1,\cdots, x_d)= x\mapsto\tilde{\mu}^{x}\in \mathcal{P}(\Sc\times\Yt) $ is Lipschitz continuous in the sense of narrow topology, if there exists a constant $\tilde{C}>0$ such that 
\begin{align*}
\textrm{d}_{BL}(\tilde{\mu}^{x}, \tilde{\mu}^{z})\leqslant \tilde{C}\vert x-z\vert. 
\end{align*}
\end{definition}

With the above notations in place, we have the following result on the Lipschitz dependence of PS-ergodic periodic measures on parameters. 
\begin{prop}\label{differentiable}\rm\;
Let $\Nb$ be a non-empty relative compact subset of $\Rd$, i.e., $\Nb\Subset\Rd$ and let $\{\tilde{\mu}^{x}_t: t\geqslant 0, x\in \Nb \}$ be a family of periodic measures generated by the SDE (\ref{Fast}) on $\Sc\times\Yt$ with coefficients satifying Assumptions \ref{A2.1} and \ref{A2.1bc}. Then, there exists $\tilde{C}>0$ such that for all $x, z\in\Nb,$ we have 
\begin{align}
\text{d}_{BL}(\bar{\tilde{\mu}}^{x}, \bar{\tilde{\mu}}^{z})\leqslant \tilde{C} \vert x -z\vert.
\end{align}
\end{prop}
\noindent {\it Proof.}
Recall from the definition of periodic measure that for $f\in \mathcal{C}_b(\Sc\times\Yt),$
\begin{align*}
\langle f, \tilde{\mu}_t^{x}\rangle =\int_{\Sc\times\Yt}fd\tilde{\mu}^{x}_t =  \int_{\Sc\times\Yt}\E[f(\tilde{S}^x(t+k\tau, \om))]d\mu_0, \quad k\in \N.
\end{align*}
In particular, for any $f\in \text{Lip}(1,\Sc\times\Yt),$ we have 
\begin{align*}
\langle f, \tilde{\mu}_t^{x}-\tilde{\mu}_t^z\rangle &= \int_{\Sc\times\Yt}\E[f(\tilde{S}^x(t+k\tau, \om))]d\tilde{\mu}_0- \int_{\Sc\times\Yt}  \E[f(\tilde{S}^z(t+k\tau, \om))] d\tilde{\mu}_0\\[.1cm]
& \hspace{1cm}\leqslant  \int_{\Sc\times\Yt}\left\vert \E[f(\tilde{S}^x(t+k\tau, \om))]- \E[f(\Phi^x(t+k\tau, \om , .))]\right\vert d\tilde{\mu}_0\\[.1cm] & \hspace{1.5cm} +  \int_{\Sc\times\Yt}\left\vert \E[ f(\tilde{S}^z(t+k\tau, \om))]- \E[f(\Phi^z(t+k\tau, \om , .))]\right\vert d\tilde{\mu}_0\\[.1cm] & \hspace{1.7cm}+
\int_{\Rd}\Big\vert \E[f(\Phi^x(t+k\tau, \omega,.))] - \E[f(\Phi^z(k\tau, \omega,.)]\Big\vert d\tilde{\mu}_0\\[.1cm]
&\hspace{1cm} \leqslant \int_{\Sc\times\Yt}\E\Big\vert \tilde{S}^x(t+k\tau, \om)- \Phi^x(t+k\tau, \om, .)\Big\vert d\tilde{\mu}_0\\
&\hspace{1.5cm} +  \int_{\Sc\times\Yt}\E\left\vert \tilde{S}^z(t+k\tau, \om)- \Phi^z(k\tau, \om, .)\right\vert d\tilde{\mu}_0\\[.1cm]
&\hspace{1.7cm} + \int_{\Sc\times\Yt}\E\Big\vert \Phi^x(t+k\tau, \om, .)- \Phi^z(t+k\tau, \om, .)\Big\vert d\tilde{\mu}_0.
\end{align*}
Next, by Step I in Theorem \ref{Ps_erg}, i.e., equation (\ref{Est_gen}) and the stability of random periodic solution, we have for any $t\in [0, \tau),$  
\begin{align}\label{ExB}
\langle f, \tilde{\mu}_t^{x}-\tilde{\mu}_t^{z}\rangle \leqslant\limsup_{k\rightarrow\infty} C\exp\left(\frac{1}{2k\tau}\int_{0}^{r+k\tau}\lambda(u)du\right)^{2k\tau}\int_{\Sc\times\Yt}\E\vert (x,\tilde{y})-(z,\tilde{y})\vert \tilde{\mu}_0(d\tilde{y}).
\end{align}
 Integrating both sides of (\ref{ExB}) over $t\in [0, \tau)$ and dividing by the length of the interval, we obtain
\begin{align}\label{ExB2}
 \frac{1}{\tau}\int_0^\tau\langle f, \tilde{\mu}_t^{x}-\tilde{\mu}_t^{z}\rangle dt\leqslant \tilde{C}\vert x-z\vert.
\end{align}
Finally, taking supremum of (\ref{ExB2}) over $f\in \text{Lip}(1,\Bt),$ we arrive at the desired inequality
\begin{align*}
\text{d}_{BL}(\bar{\tilde{\mu}}^{x}, \tilde{\tilde{\mu}}^{z})=\sup\left\{ \frac{1}{\tau}\int_0^\tau\langle f, \tilde{\mu}_t^{x}-\tilde{\mu}_t^{z}\rangle dt:  f\in \text{Lip}(1,\Bt)\right\}\leqslant \tilde{C}\vert x-z\vert,
\end{align*}
\qed 

\medskip
 Now, given the ergodic property of periodic measures $\{\tilde{\mu}^x_r: r\in [0, \tau),\; x\in\Nb\Subset\Rd \}$ and Lipschitz dependence on parameter $x\in \Nb$, we identify the averaged vector field $\bar{F}$ from the slow vector filed $F$. Importantly, we show that $\bar{F}$ is regular enough to ensure local solvability of the corresponding averaged equation. The argument here is fairly straightforward, due to the regularity of the slow vector field and that of the periodic measures with respect to parameters $x\in \Nb.$ Consider a vector field $\bar{F}: \Rd\rightarrow\Rd$ defined by 
\begin{align}
\bar{F}(x) = \frac{1}{\tau}\int_{0}^{\tau}\int_{\Sc\times\Yt}F(x,\pi_2(\tilde{y}))\tilde{\mu}_r^x(d\tilde{y}) dr= \int_{\Sc\times\Yt}F(x,\pi_2(\tilde{y}))\bar{\tilde{\mu}}^x(d\tilde{y})
\end{align}
where $\pi_2:\Sc\times\Yt\rightarrow\Yt$ is a projection map, i.e., $\pi_2(\tilde{y}) = \pi_2(s,y)= y.$ 
\begin{prop}\label{Averaged}\rm\;
Let $F\in \mathcal{C}^{l-1}_b(\Rd\times\Yt;\Rd),\; l\geqslant 2$ and let $\{\tilde{\mu}_r^x: r\in [0,\tau),\; x\in \Nb\Subset\Rd\}$ be a family of periodic measures generated by the fast SDE satisfying Assumptions \ref{A2.1} and \ref{A2.1bc}. Then, the family of vector fields
\begin{align}
\left\{\int_{\Sc\times\Yt}F(x, \pi_2(\tilde{y}))\}\tilde{\mu}_r^x(d\tilde{y}): r\in [0, \tau), \; x\in \Nb\right\} 
\end{align}
is $\tau$-periodic. Moreover, there exists $\tilde{C}_F>0$ such that  for all $x, z\in \Nb$
\begin{align}\label{Local}
\begin{cases}
\vert \bar{F}(x)- \bar{F}(z)\vert \leqslant\tilde{C}_F\vert x-z\vert, \\ \vert F(x)\vert \leqslant\tilde{C}_F.
\end{cases}
\end{align}
\end{prop}
\noindent {\it Proof.} The periodicity part follows from the $\tau$-periodicity of the periodic measure $\tilde{\mu}_r^x$ and the fact that $F$ is time independent. It only remains to prove  the local Lipschitz and boundedness property (\ref{Local}). As  $F\in \mathcal{C}_b^{l-1}(\Rd\times\Yt;\Rd), \; l\geqslant 2$, then by the mean value theorem and Lipschitz property of $x\mapsto \tilde{\mu}_r, \; r\in [0, \tau),$ we have the existence $L_F, \tilde{C}\geqslant 0,$ such that for $x, z\in\Nb\Subset\Rd$
\begin{align*}
\vert \bar{F}(x)-\bar{F}(z)\vert &= \left\vert \frac{1}{\tau}\int_0^\tau \int_{\Bt}F(x,\pi_2(\tilde{y}))\tilde{\mu}^x_r(d\tilde{y})dr -\frac{1}{\tau}\int_0^\tau \int_{\Bt}F(z,\pi_2(\tilde{y}))\tilde{\mu}^z_r(d\tilde{y})dr \right\vert \\
&=\left\vert \frac{1}{\tau}\int_0^\tau\int_{\Bt}\E\left[F(x, \Pi\tilde{S}^x(r,\om))-F(z, \Pi\tilde{S}^z(r,\om))\right]d\tilde{\mu}_0dr\right\vert\\ 
&\hspace{1cm} \leqslant L_F\vert x-z\vert +L_F\frac{1}{\tau}\int_0^\tau \int_{\Bt}\E\vert \Pi\tilde{S}^x(r,\om)-\Pi\tilde{S}^z(r,\om)\vert d\tilde{\mu}_0dr\\
&\hspace{1.5cm}\leqslant L_F\vert x-z\vert +L_F\tilde{C}\vert x-z\vert \leqslant \tilde{C}_F\vert x-z\vert.
\end{align*}
 For the second assertion, since $F\in \mathcal{C}_b^{l-1}(\Rd\times\Yt;\Rd)$,  there exists $L_F>0$ such that 
 \begin{align}\label{Linear_grow}
 \vert F(x,y)\vert \leqslant L_F(1+\vert x\vert +\vert y\vert).
 \end{align}
 Now, consider the vector fields $F_n(x,y), \; n\in \N,$ defined by 
 \begin{align*}
 F_n(x,\pi_2(\tilde{y})) = \begin{cases} F(x,\pi_2(\tilde{y})), \quad  &\text{if}\;\; \vert \tilde{y}\vert \leqslant 2^n\\
 F\left(x, \frac{2^n\pi_2(\tilde{y})}{\vert\pi_2(\tilde{y})\vert}\right ),\quad & \text{if}\;\; \vert \tilde{y}\vert >2^n
 \end{cases}.
\end{align*}  
It follows that $F_n(x,.): \Yt\rightarrow\Yt$ is bounded Lipschitz continuous and for any $R>0$, we have
\begin{align*}
\sup_{\vert x\vert \leqslant R}\vert F_n(x,.)\vert \leqslant \tilde{C}_F.
\end{align*}
We use the inequality (\ref{Linear_grow}) and the moment bound on the periodic measures or random periodic solutions, i.e., equation (\ref{Est_gen}), to obtain 
\begin{align*}
 \left\vert \int_{\{\tilde{y}: \vert\tilde{y}\vert>2^n \}}\left[F(x,\pi_2(\tilde{y}))- F_n(x,\pi_2(\tilde{y}))\right]\bar{\tilde{\mu}}^x(d\tilde{y})\right\vert
& \leqslant L_F\int_{\{\tilde{y}: \vert \tilde{y}\vert>2^n\}}(1+\vert x\vert + \vert \pi_2(y)\vert)\bar{\tilde{\mu}}^x(d\tilde{y})\\ &\hspace{1cm} \leqslant 2^{-n}K_F(1+ \vert x\vert^2).
\end{align*}
This implies that, for any $\delta>0$ and $R>0,$ we can choose $\tilde{n} =n(\delta.R)\in \N\setminus\{0\}$ large enough such that 
\begin{align*}
\sup_{x\in B_R}\left\vert \int_{\{\tilde{y}: \vert\tilde{y}\vert>2^{\tilde{n}} \}}\left[F(x,\pi_2(\tilde{y}))-F_{\tilde{n}}(x,\pi_2(\tilde{y}))\right]\bar{\tilde{\mu}}^x(d\tilde{y})\right\vert<\delta.
\end{align*}
So, we obtain 
\begin{align*}
\vert \bar{F}(x)\vert &\leqslant \left\vert\int_{\Bt}F_n(x,\pi_2(\tilde{y}))\bar{\tilde{\mu}}^x(d\tilde{y})\right\vert + \left\vert \int_{\{\tilde{y}: \vert\tilde{y}\vert>2^n \}}\left[F(x,\pi_2(\tilde{y}))- F_n(x,\pi_2(\tilde{y}))\right]\bar{\tilde{\mu}}^x(d\tilde{y})\right\vert\\
& \hspace{1cm}\leqslant \tilde{C}_F+\delta.
\end{align*}
Finally, since $\delta>0$ is arbitrary, we may take limit as $\delta \rightarrow 0$ to arrive at the required bound. 

\qed 
\subsection{The Averaging limit}\label{AV_limit}
In this subsection, we present an averaging result in the random periodic regime; we follow Hasminskii's time discretisation scheme (e.g., \cite{Fredlinbook, Hasm66, Kifer09}). The idea of the scheme is to decompose the time interval $[0, t]$ into subintervals of length $\frac{1}{n(\varepsilon)} = \Delta(\varepsilon)$, such that in each subinterval the slow variable $X^{\varepsilon}$ is almost everywhere constant and the fast variable $\tilde{Y}_t^{\varepsilon}$ after appropriate rescaling is well captured by the long time behaviour of the one-point motion $\tilde{Y}^{\varepsilon, x,\tilde{y}}_t$ for fixed $x\in \Nb\Subset\Rd$.
First, we write the slow-fast SDE in the fast time scale $t\mapsto t/\varepsilon$  as follows,
\begin{align*}
\begin{cases} dX_t^\varepsilon = F(X_t^\varepsilon, \Pi\tilde{Y}_t^\varepsilon)dt\\
d\tilde{Y}_t^{\varepsilon} = \frac{1}{\varepsilon}\tilde{b}(X^\varepsilon_t,\tilde{Y}_t^{\varepsilon})dt+\frac{1}{\sqrt{\varepsilon}}\tilde{\sigma}(X_t^\varepsilon, \tilde{Y}_t^{\varepsilon})d\tilde{W}_t,
\end{cases}
\end{align*}
where $\Pi: \Sc\times\Yt\rightarrow\Yt$  is a projection map.
Next, we write $[0, t]$ as a union of subintervals  $[\frac{jt}{n}, \frac{(j+1)t}{n}]: j=0, 1,\cdots n-1$ and set $\hat{\tilde{Y}}_u^{\varepsilon}$ as an auxillary process defined on each subinterval as follows (c.f.~\cite{Fredlinbook})
\begin{align*}
\begin{cases}
\hat{\tilde{Y}}^{\varepsilon}_{jt/n} = \tilde{Y}^\varepsilon_{jt/n},\\
d\hat{\tilde{Y}}_{u}^{\varepsilon,x,\tilde{y}} = \frac{1}{\varepsilon}\tilde{b}(X^{\varepsilon}_{jt/n},\hat{\tilde{Y}}_u^{\varepsilon})du +\frac{1}{\sqrt{\varepsilon}}\tilde{\sigma}(X^{\varepsilon}_{jt/n},\hat{\tilde{Y}}_u^{\varepsilon})d\tilde{W}_u, \quad jt/n<u\leqslant (j+1)t/n,
\end{cases}
\end{align*}
where $\tilde{W}_t$ is the same Brownian path used in the definition of $\tilde{Y}^{\varepsilon}$ as a stochastic flow on $\Sc\times\Yt.$

\medskip
Next, we give some preparatory lemmas on how to choose the integer $n(\varepsilon)$ such that the interval spacing is not too small as PS-ergodicity need some time to occur.
\begin{lem}[c.f. \cite{Wain13}]\label{NLema}\rm\;
Let $F\in \mathcal{C}^{l-1}_b(\Rd\times\Yt;\Rd),\; b, \sigma_k\in \mathcal{C}_b^{l}(\R\times\Rd\times\Yt;\Yt),\; l\geqslant 2, \; 1\leqslant k\leqslant N$. Then, there exists $K>0$ such that 
\begin{align}
\sup_{0\leqslant u\leqslant t}\E\vert\tilde{Y}_{u}^{\varepsilon,x,\tilde{y}}-\hat{\tilde{Y}}^{\varepsilon,x,\tilde{y}}_{u}\vert^2  \leqslant K\left(\frac{1}{\varepsilon^2n^3}+\frac{1}{\varepsilon n^2}\right)\exp\left(\frac{K}{\varepsilon^2n^2}+\frac{K}{\varepsilon n}\right).
\end{align}
\end{lem}
\noindent {\it Proof.}
Let $u\in [0, t],$ then there exists $j = j(u)$ such that for any $jt/n\leqslant u\leqslant (j+1)t/n$ we have 
\begin{align*}
\tilde{Y}_{u}^{\varepsilon,x,\tilde{y}}-\hat{\tilde{Y}}_{u}^{\varepsilon,x,\tilde{y}} &= \frac{1}{\varepsilon}\int_{jt/n}^u\left[\tilde{b}(X_\xi^{\varepsilon,x,\pi_2(\tilde{y})},\tilde{Y}_\xi^{\varepsilon,x,\tilde{y}})-\tilde{b}(X_{jt/n}^{\varepsilon,x,\pi_2(\tilde{y})},\hat{\tilde{Y}}_\xi^{\varepsilon,x,\tilde{y}})\right]d\xi\\
&\hspace{.5cm}+ \frac{1}{\sqrt{\varepsilon}}\int_{jt/n}^u\left[\tilde{\sigma}(X_\xi^{\varepsilon,x,\pi_2(\tilde{y})},\tilde{Y}_\xi^{\varepsilon,x,\tilde{y}})-\tilde{\sigma}(X_{jt/n}^{\varepsilon,x,\pi_2(\tilde{y})},\hat{\tilde{Y}}_\xi^{\varepsilon,x,\tilde{y}})\right]d\tilde{W}_\xi.
\end{align*}
Next, we use Cauchy Schwartz inequality on the drift part and It\^o isometry on the diffusion part to obtain
\begin{align*}
\E\left\vert \tilde{Y}_{u}^{\varepsilon,x,\tilde{y}}-\hat{\tilde{Y}}_{u}^{\varepsilon,x,\tilde{y}}\right\vert^2 &\leqslant \frac{1}{\varepsilon^2n}\int_{jt/n}^u\E\left\vert \tilde{b}(X_\xi^{\varepsilon,x,\pi_2(\tilde{y})},\tilde{Y}_\xi^{\varepsilon,x,\tilde{y}})-\tilde{b}(X_{jt/n}^{\varepsilon,x,\pi_2(\tilde{y})},\hat{\tilde{Y}}_\xi^{\varepsilon,x,\tilde{y}})\right\vert^2d\xi\\
&\hspace{.5cm} + \frac{1}{\varepsilon}\int_{jt/n}^u\E\left\vert \tilde{\sigma}(X_\xi^{\varepsilon,x,\pi_2(\tilde{y})},\tilde{Y}_\xi^{\varepsilon,x,\tilde{y}})-\tilde{\sigma}(X_{jt/n}^{\varepsilon,x,\pi_2(\tilde{y})},\hat{\tilde{Y}}_\xi^{\varepsilon,x,\tilde{y}})\right\vert ^2d\xi.
\end{align*}
By the regularity of $b, \sigma,$ then by the mean value theorem, there exist $L_b, L_\sigma>0$ such that 
\begin{align*}
\E\left\vert \tilde{Y}_{u}^{\varepsilon,x,\tilde{y}}-\hat{\tilde{Y}}_{u}^{\varepsilon,x,\tilde{y}}\right\vert^2 &\leqslant\frac{L_b}{\varepsilon^2n}\int_{jt/n}^u\left\{ \E\vert X_{\xi}^{\varepsilon,x,\pi_2(\tilde{y})}- X_{jt/n}^{\varepsilon,x,\pi_2(\tilde{y})}\vert^{2} + \E\vert \tilde{Y}_\xi^{\varepsilon,x,\tilde{y}}- \hat{\tilde{Y}}_\xi^{\varepsilon,x,\tilde{y}}\vert^{2}\right\}d\xi\\
&\hspace{.5cm} +\frac{L_\sigma}{\varepsilon}\int_{jt/n}^u\left\{ \E\vert X_{\xi}^{\varepsilon,x,\pi_2(\tilde{y})}- X_{jt/n}^{\varepsilon,x,\pi_2(\tilde{y})}\vert^{2} + \E\vert \tilde{Y}_\xi^{\varepsilon,x,\tilde{y}}- \hat{\tilde{Y}}_\xi^{\varepsilon,x,\tilde{y}}\vert^{2}\right\}d\xi.
\end{align*}
Next, by the regularity of $F$, there exist $L_F>0$ such that 
\begin{align*}
\E\vert X_{\xi}^{\varepsilon,x,\pi_2(\tilde{y})}- X_{jt/n}^{\varepsilon,x,\pi_2(\tilde{y})}\vert^{2}\leqslant L_F \vert \xi-jt/n\vert,
\end{align*}
so that 
\begin{align*}
\E\left\vert \tilde{Y}_{u}^{\varepsilon,x,\tilde{y}}-\hat{\tilde{Y}}_{u}^{\varepsilon,x,\tilde{y}}\right\vert^2 &\leqslant L_F\left(\frac{L_b}{\varepsilon^2n}+ \frac{L_\sigma}{\varepsilon}\right)\int_{jt/n}^u\vert \xi-jt/n\vert d\xi \\ &\hspace{1cm}+ \left(\frac{L_b}{\varepsilon^2 n}+\frac{L_\sigma}{\varepsilon}\right)\int_{jt/n}^u\E\vert \tilde{Y}_\xi^{\varepsilon,x,\tilde{y}}-\hat{\tilde{Y}}_{jt/n}^{\varepsilon,x,\tilde{y}}\vert^2d\xi\\
&\leqslant K\left( \frac{1}{\varepsilon^2 n^3}+\frac{1}{\varepsilon n^2}\right)+ K\left(\frac{1}{\varepsilon^2 n}+ \frac{1}{\varepsilon}\right)\int_{jt/n}^u\E\vert \tilde{Y}_\xi^{\varepsilon,x,\tilde{y}}-\hat{\tilde{Y}} _\xi^{\varepsilon,x,\tilde{y}}\vert^2d\xi,
\end{align*}
 Finally, by Gronwall lemma, we have the desired bound
\begin{align*}
\sup_{0\leqslant u\leqslant t}\E\vert\tilde{Y}_{u}^{\varepsilon,x,\tilde{y}}-\hat{\tilde{Y}}^{\varepsilon,x,\tilde{y}}_{u}\vert^2  \leqslant K\left(\frac{1}{\varepsilon^2n^3}+\frac{1}{\varepsilon n^2}\right)\exp\left(\frac{K}{\varepsilon^2n^2}+\frac{K}{\varepsilon n}\right)
\end{align*}
\qed

We choose $n(\varepsilon)= \frac{1}{\varepsilon\sqrt[4]{\ln(\varepsilon^{-1})})},$ thanks to the above Lemma \ref{NLema} (see, also~\cite{Fredlinbook, Kifer09}).  The next lemma deals with the interval spacing. 
We take $t_0>0$ and consider the solution $\tilde{Y}_u^{\varepsilon,x,\tilde{y}}$ of the equation 
\begin{align*}
\begin{cases}
d\tilde{Y}_u^{\varepsilon,x,\tilde{y}} = \frac{1}{\varepsilon}\tilde{b}(x,\tilde{Y}_u^{\varepsilon,x,\tilde{y}})du+ \frac{1}{\sqrt{\varepsilon}}\tilde{\sigma}(x,\tilde{Y}_u^{\varepsilon,x,\tilde{y}})d\tilde{W}_u, \quad u>t_0\\
\tilde{Y}_{t_0}^{\varepsilon,x,\tilde{y}} = \tilde{y}\in \Sc\times\Yt.
\end{cases}
\end{align*}

\begin{lem}\label{Discret}\rm\;
Let $F\in \mathcal{C}^{l-1}_b(\Rd\times\Yt;\Rd),\; b, \sigma_k\in \mathcal{C}_b^{l}(\R\times\Rd\times\Yt;\Yt),\; l\geqslant 2, \; 1\leqslant k\leqslant N$ and Assumptions \ref{A2.1} and \ref{A2.1bc} hold.
Then, there exist constants $C>0$ such for any $T>0,\; x\in\Rd, \; \tilde{y}\in \Sc\times\Yt$, we have   
\begin{align*}
\E\left\vert \frac{1}{T}\int_{t_0}^{t_0+T}F(x, \Pi\tilde{Y}_u^{\varepsilon,x,\tilde{y}})du - \bar{F}(x)\right\vert^2\leqslant \frac{C}{T}\left(1+\vert x\vert^{2}+\vert\tilde{y}\vert^{2}\right)^2+\alpha(T,x)
\end{align*}
for some mapping $\alpha: \Rp\times\Rd\rightarrow\Rp$ with 
\begin{align*}
\sup_{T\in \Rp}\alpha(T,x) \leqslant C(1+ \vert x\vert^2)
\end{align*}
and, for any $\Nb\Subset\Rd$
\begin{align*}
\lim_{T\rightarrow\infty}\sup_{x\in \Nb}\alpha(T,x) =0.
\end{align*}
\end{lem}
\noindent {\it Proof.} We proceed as in the proof of Lemma 7.2 in \cite{Cera17} or as in the proof of Lemma 3.2 in \cite{Wain13}, we have 
\begin{align*}
&\E\left(\frac{1}{T}\int_{t_0}^{t_0+T}F(x,\Pi\tilde{Y}_u^{\varepsilon,x,\tilde{y}})du-\bar{F}(x)\right)^2\\
&\hspace{1cm}=\frac{1}{T^2}\int_{t_0}^{t_0+T}\int_{t_0}^{t_0+T}\E\left[ \left\langle F(x, \Pi\tilde{Y}_u^{\varepsilon,x,\tilde{y}})-\bar{F}(x),  F(x, \Pi\tilde{Y}_r^{\varepsilon,x,\tilde{y}})-\bar{F}(x)\right\rangle\right]dudr\\
&\hspace{1cm} = \frac{2}{T^2}\int_{t_0}^{t_0+T}\int_{r}^{t_0+T}\E\left[ \left\langle F(x, \Pi\tilde{Y}_u^{\varepsilon,x,\tilde{y}})-\bar{F}(x),  F(x, \Pi\tilde{Y}_r^{\varepsilon,x,\tilde{y}})-\bar{F}(x)\right\rangle\right]dudr.
\end{align*}
Next, by Markov property of $\Pi\tilde{Y}_t^{\varepsilon,x,\tilde{y}}$, we have for $r\leqslant u,$
\begin{align*}
& \E\left[ \left\langle F(x, \Pi\tilde{Y}_u^{\varepsilon,x,\tilde{y}})-\bar{F}(x),  F(x, \Pi\tilde{Y}_r^{\varepsilon,x,\tilde{y}})-\bar{F}(x)\right\rangle\right] \\ &\hspace{1cm}= \E\left[ \left\langle F(x, \Pi\tilde{Y}_r^{\varepsilon,x,\tilde{y}})-\bar{F}(x),  \tilde{P}_{u-r}(F(x,\Pi\tilde{Y}_r^{\varepsilon,x,\tilde{y}})-\bar{F}(x))\right\rangle\right].
\end{align*}
By Cauchy-Schwartz inequality, we have 
\begin{align}\label{DW1}
\notag &\E\left(\frac{1}{T}\int_{t_0}^{t_0+T}F(x,\Pi\tilde{Y}_u^{\varepsilon,x,\tilde{y}})du-\bar{F}(x)\right)^2\\
\notag &\hspace{2cm} \leqslant \frac{2}{T^2}\int_{t_0}^{t_0+T}\int_{r}^{t_0+T}\bigg\{\left(\E\vert F(x,\Pi\tilde{Y}_r^{\varepsilon,x,\tilde{y}})-\bar{F}(x)\vert^2\right)^{1/2}\\ &\hspace{3cm} \left(\E\vert \tilde{P}_{u-r}F(x,\Pi\tilde{Y}_r^{\varepsilon,x,\tilde{y}})-\bar{F}(x)\vert^2\right)^{1/2}\bigg\}dudr.
\end{align}
Next, we estimate the integrand using the regularity of vector fields $F,\; \bar{F}$ and PS-ergodicity of the fast motion (see \S\ref{Ranp_M}), these yield the existence of $C>0$ such that 
\begin{align}\label{GQ}
\notag \E\vert F(x,\Pi\tilde{Y}_r^{\varepsilon,x,\tilde{y}})-\bar{F}(x)\vert^2 & \leqslant C\left( 1+\vert x\vert^2+ \E\vert \tilde{Y}_r^{\varepsilon,x,\tilde{y}}\vert^2\right)\\
&\hspace{.5cm} \leqslant C\left( 1+\vert x\vert^2+\exp\left(\int_{t_0}^{r}\lambda(\xi)d\xi\right)\vert \tilde{y}\vert^2\right).
\end{align}
Also, as $F\in \mathcal{C}_b^{l-1}(\Rd\times\Yt;\Rd),\; l\geqslant 2,$ there exists $L_F>0$ such that 
\begin{align}\label{GQ1}
\vert F(x,y)-F(x,z)\vert \leqslant L_F\vert y-z\vert.
\end{align}
The inequality (\ref{GQ1}) together with PS-ergodicity of the fast variables yield 
\begin{align}\label{GQ2}
\E\vert \tilde{P}_{u-r}F(x,\Pi\tilde{Y}_u^{\varepsilon,x,\tilde{y}})-\bar{F}(x)\vert^2 \leqslant C\left( 1+\vert x\vert^2+\vert\tilde{y}\vert^2\right)\exp\left(\int_{t_0}^{u-r}\lambda(\xi)d\xi\right).
\end{align}
We substitute the above estimates (\ref{GQ}) and (\ref{GQ2}) into the integral (\ref{DW1}) to obtain
\begin{align*}
&\E\left(\frac{1}{T}\int_{t_0}^{t_0+T}F(x,\Pi\tilde{Y}_u^{\varepsilon,x,\tilde{y}})du-\bar{F}(x)\right)^2\\
&\hspace{1cm} \leqslant \frac{2C}{T^2}\int_{t_0}^{t_0+T}\int_{r}^{t_0+T}(1+\vert x\vert^2+\vert \tilde{y}\vert^2)^2\exp\left(\int_{t_0}^{r}\lambda(\xi)d\xi\right) dudr\\
&\hspace{1.5cm}\leqslant 2C\left( 1+\vert x\vert^2+\vert \tilde{y}\vert^2\right)^2\frac{1}{T^2}\int_{t_0}^{t_0+T}\int_{r}^{t_0+T}\exp\left(\int_{t_0}^{u-r}\lambda(\xi)d\xi\right)dudr\\
&\hspace{2cm}\leqslant \frac{\tilde{C}}{T} \left( 1+\vert x\vert^2+\vert \tilde{y}\vert^2\right)^2.
\end{align*}
Finally, from Proposition \ref{Averaged}, we see that the family of vector fields
\begin{align*}
\left\{\int_{\Sc\times\Yt}F(x, \pi_2(\tilde{y}))\}\tilde{\mu}_r^x(d\tilde{y}): r\in [0, \tau), \; x\in \Nb\right\} 
\end{align*}
is $\tau$-periodic, then by PS-ergodicity of the fast motion,  we know that 
\begin{align*}
\lim_{k\rightarrow\infty}\frac{1}{k\tau}\int_{t_0}^{t_0+k\tau}\int_{\Bt}F(x,\pi_2(\tilde{y}))\tilde{\mu}_t^x(d\tilde{y})dt = \bar{F}(x).
\end{align*}
In view of this, we may define the mapping $\alpha: \Rp\times\Rd\rightarrow\Rp$ by 
\begin{align*}
\alpha(T,x) = 2\left\vert \frac{1}{T}\int_{t_0}^{t_0+T}\int_{\Bt}F(x,\pi_2(\tilde{y}))\tilde{\mu}^x_tdt - \bar{F}(x)\right\vert^2 
\end{align*}
\qed

With the above preparatory lemmas and notations in place, we are ready to state and prove a stochastic averaging principle in the random periodic regime.
\begin{theorem}\rm\;
Assume that $F\in \mathcal{C}_b^{l-1}(\Rd\times\Yt; \Rd)$ and $ b,\sigma_k\in \mathcal{C}_b^{l,\delta}(\R\times\Rd\times\Yt; \Yt), \; 1\leqslant k\leqslant N, \; l\geqslant 2$. If for fixed $x\in\Nb\Subset \Rd,$ the coefficients $b, \sigma_k,\; 1\leqslant k\leqslant N$ are  time periodic and satisfy Assumptions \ref{A2.1} and \ref{A2.1bc}. Moreover, if $\{\Bt_r: r\in [0, \tau)\}=:\Bt\subset \Sc\times\Yt$ is such that each $\Bt_r$ is a $k_r\tau$-irreducible Poincar\'e section, then
\begin{equation}\label{F2Q11}
\lim_{\varepsilon\rightarrow 0}\int_{\Nb}\int_{\Bt}\E\left(\sup_{0\leqslant t\leqslant T/\varepsilon}\Big\vert X_t^{\varepsilon,x,\pi_2(\tilde{y})}-\bar{X}_t^{\varepsilon,x}\Big\vert\right)\bar{\tilde{\mu}}^x(d\tilde{y})\nu(dx) =0,
\end{equation}
where $\nu\in \mathcal{P}(\Nb)$ and $\bar{\tilde{\mu}}^x\in \mathcal{P}(\Bt).$
\end{theorem}
\noindent {\it Proof.} 

We start by using the Lipschitz property of the averaged vector field $\bar{F}$ to obtain 
\begin{align*}
\bigg\vert X_t^{\varepsilon,x,\pi_2(\tilde{y})}-\bar{X}_t^{\varepsilon,x}\bigg\vert &= \varepsilon\bigg\vert \int_0^t\Big[ F\left(X_u^{\varepsilon,x,\pi_2(\tilde{y})},\Pi\tilde{Y}_u^{\varepsilon,x,\tilde{y}} \right)- \bar{F}(\bar{X}_u^{\varepsilon,x})\Big]du\bigg\vert\\
&\leqslant \varepsilon \bigg\vert \int_0^t\Big[ F\left(X_u^{\varepsilon,x,\pi_2({y})},\Pi\tilde{Y}_u^{\varepsilon,x,\tilde{y}}\right)- \bar{F}(X_u^{\varepsilon,x, \tilde{y}})\Big]du\bigg\vert+\varepsilon\tilde{C}_F\int_0^t\Big\vert X_u^{\varepsilon,x,\tilde{y}}-\bar{X}_u^{\varepsilon,x}\Big\vert du
\end{align*}
and by Gronwall's lemma, we have 
\begin{align}\label{F3Q1}
\sup_{0\leqslant t\leqslant T/\varepsilon}\Big\vert X_t^{\varepsilon,x,\pi_2(\tilde{y})}-\bar{X}_t^{\varepsilon,x}\Big\vert \leqslant \varepsilon e^{T\tilde{C}_F} \sup_{0\leqslant t\leqslant T/\varepsilon}\bigg\vert \int_0^t\Big[ F\left(X_u^{\varepsilon,x,\pi_2(\tilde{y})},\Pi\tilde{Y}_u^{\varepsilon, x,\tilde{y}}\right)-\bar{F}(X_u^{\varepsilon,x,\pi_2(\tilde{y})})\Big]du\bigg\vert.
\end{align}
Next, from Lemmas \ref{NLema} and \ref{Discret}, we can choose $n(\varepsilon) = \frac{1}{\varepsilon\sqrt[4]{\ln(\varepsilon^{-1})})}$ and set 
$t(\varepsilon) = \frac{T}{\varepsilon n(\varepsilon)}$ (see, also~\cite{Fredlinbook, Kifer09}) and write the integral in RHS of the inequality (\ref{F3Q1}) in the form 
\begin{align*}
&\varepsilon \bigg\vert \int_{0}^{t(\varepsilon)} \left[F\Big(X_{jt(\varepsilon)}^{\varepsilon,x,\pi_2(\tilde{y})}, \Pi\tilde{Y}_{jt(\varepsilon)+u}^{\varepsilon,x,\tilde{y}}\Big)du-\bar{F}\big(X_{jt(\varepsilon)}^{\varepsilon,x,\pi_2(\tilde{y})}\big)\right]\bigg\vert\\[.2cm]
&\hspace{.3cm}=\varepsilon \bigg\vert \int_0^{t(\varepsilon)}\bigg\{\Big[ F\left(X_{jt(\varepsilon)}^{\varepsilon,x,\pi_2({y})}, \Pi\tilde{Y}_{jt(\varepsilon)+u}^{\varepsilon, x,\tilde{y}}\right)-F\left(X_{jt(\varepsilon)+u}^{\varepsilon,x,\pi_2(\tilde{y})}, \Pi\tilde{Y}_{jt(\varepsilon)+u}^{\varepsilon, x,\tilde{y}}\right)\Big]\\[.2cm] & \hspace{.7cm}+\Big[F\left(X_{jt(\varepsilon)+u}^{\varepsilon,x,\pi_2(\tilde{y})}, \Pi\tilde{Y}_{jt(\varepsilon)+u}^{\varepsilon, x,\tilde{y}}\right)-\bar{F}\left(X_{jt(\varepsilon)+u}^{\varepsilon,x, \pi_2(\tilde{y})}\right)\Big]
+\Big[\bar{F}\left(X_{jt(\varepsilon)+u}^{\varepsilon,x, \pi_2(\tilde{y})}\right)
-\bar{F}\left(X_{jt(\varepsilon)}^{\varepsilon,x, \pi_2(\tilde{y})}\right)\Big]\bigg\}du\bigg\vert \\[.2cm]
&\hspace{1.3cm} =: \mathcal{I}_1+\mathcal{I}_2+\mathcal{I}_3.
\end{align*}
Now, by the regularity of the vector field $F,$ we have the existence of $L_F>0,$ such that 
\begin{align}\label{F2Q6}
\notag \mathcal{I}_1= \varepsilon \bigg\vert\int_{0}^{t(\varepsilon)} \Big[ F\left(X_{jt(\varepsilon)+u}^{\varepsilon,x,\pi_2(\tilde{y})},\tilde{Y}_{jt(\varepsilon)+u}^{\varepsilon,x,\tilde{y}}\right)-F\left(X_{jt(\varepsilon)}^{\varepsilon,x,\pi_2(\tilde{y})}, \tilde{Y}_{jt(\varepsilon)+u}^{\varepsilon,x,\tilde{y}}\right)\Big]du\bigg\vert \\ \leqslant \varepsilon L_{F}\int_0^{t(\varepsilon)}\left\vert X_{jt(\varepsilon)+u}^{\varepsilon,x,\pi_2(\tilde{y})}-X_{jt(\varepsilon)}^{\varepsilon,x,\pi_2(\tilde{y})}\right\vert\leqslant \tilde{L}_{F}(\varepsilon t(\varepsilon))^2.
\end{align}
In a similar fashion, by the regularity of $F$ and that of averaged vector field $\bar{F}$, we obtain
\begin{align}\label{F2Q7}
\mathcal{I}_3= \varepsilon \bigg\vert \int_{0}^{t(\varepsilon)}\Big[ \bar{F}\left(X_{jt(\varepsilon)+u}^{\varepsilon,x, \pi_2(\tilde{y})}\right)-\bar{F}\left(X_{jt(\varepsilon)}^{\varepsilon,x, \pi_2(\tilde{y})}\right)\Big]du\bigg\vert \leqslant [\tilde{L}_{F}+ \tilde{C}_{F}L_F](\varepsilon t(\varepsilon))^2.
\end{align}
Given $n(\varepsilon)$ with $t(\varepsilon) = \frac{T}{\varepsilon n(\varepsilon)},$ we use (\ref{F2Q6}) and (\ref{F2Q7}) to obtain the following inequality (see also \cite{Kifer09}),
\begin{align}\label{F3Q2}
\notag &\varepsilon\E\left( \sup_{0\leqslant t\leqslant T/\varepsilon}\bigg\vert \int_0^t\Big[F\left(X_u^{\varepsilon,x,\pi_2(\tilde{y})}, \Pi\tilde{Y}_u^{\varepsilon,x,\tilde{y}}\right)-\bar{F}(X_u^{\varepsilon,x,\pi_2(\tilde{y})})\Big]du\bigg\vert\right)\\[.1cm]
\notag &\hspace{.5cm}\leqslant 2L_F\varepsilon t(\varepsilon)+\varepsilon\E\left(\sum_{j=0}^{n(\varepsilon)-1}\bigg\vert \int_{jt(\varepsilon)}^{(j+1)t(\varepsilon)}\Big[ F\left(X_u^{\varepsilon,x,\pi_2(\tilde{y})},\Pi\tilde{Y}_u^{\varepsilon,x,\tilde{y}}\right)-\bar{F}(X_u^{\varepsilon,x,\pi_2(\tilde{y})})\Big]du\bigg\vert\right)\\[.1cm]
\notag &\hspace{.5cm} = 2L_F\varepsilon t(\varepsilon)+\varepsilon\E\left(\sum_{j=0}^{n(\varepsilon)-1}\bigg\vert \int_0^{t(\varepsilon)}\Big[ F\left(X_{jt(\varepsilon)+ u}^{\varepsilon,x,\pi_2(\tilde{y})},\Pi\tilde{Y}_{jt(\varepsilon)+u}^{\varepsilon,x,\tilde{y}}\right)-\bar{F}(X_{jt(\varepsilon)+u}^{\varepsilon,x,\pi_2(\tilde{y})})\Big]du\bigg\vert\right)\\[.1cm]
\notag &\hspace{.7cm}\leqslant 2L_F\varepsilon t(\varepsilon)+ \varepsilon^2n(\varepsilon)t(\varepsilon)[\tilde{L}_F+\tilde{C}_FL_F]\\[.1cm] \notag &\hspace{1.5cm}+\varepsilon t(\varepsilon)\left(\sum_{j=0}^{n(\varepsilon)-1}\E\bigg\vert \frac{1}{t(\varepsilon)}\int_0^{t(\varepsilon)}F\left(X_{jt(\varepsilon)}^{\varepsilon,x,\pi_2(\tilde{y})}, \Pi\tilde{Y}_{jt(\varepsilon)+u}^{\varepsilon,x,\tilde{y}}\right)du-\bar{F}(X_{jt(\varepsilon)}^{\varepsilon,x,\pi_2(\tilde{y})})\bigg\vert\right)\\[.1cm]
\notag &\hspace{.7cm}\leqslant 2L_F\varepsilon t(\varepsilon)+ \varepsilon^2n(\varepsilon)t(\varepsilon)[\tilde{L}_F+\tilde{C}_FL_F]\\[.1cm]  &\hspace{1.5cm}+\varepsilon t(\varepsilon)\sum_{j=0}^{n(\varepsilon)-1}
\left(\E\left\vert \frac{1}{t(\varepsilon)}\int_{0}^{t(\varepsilon)}F\left( X_{jt(\varepsilon)}^{\varepsilon,x,\pi_2(\tilde{y})},\Pi\tilde{Y}_{jt(\varepsilon)+u}^{\varepsilon,x,\tilde{y}}\right)du-\bar{F}(X_{jt(\varepsilon)}^{\varepsilon,x,\pi_2(\tilde{y})})\right\vert^2\right)^{1/2}.
\end{align}
Recalling from the discretisation scheme that the slow motion $X_{jt(\varepsilon)}^{\varepsilon,x,\tilde{y}}, j=0,\cdots n(\varepsilon)-1,$ is $\p(d\om)\bar{\tilde{\mu}}^x(d\tilde{y})\nu(dx)$ almost everywhere constant on each sub-interval of length $t(\varepsilon),$ then by Lemma \ref{Discret}, we have 
\begin{align}\label{D34e}
\lim_{\varepsilon\downarrow 0}\max_{0\leqslant j\leqslant n(\varepsilon)-1}\int_{\Nb}\int_{\Bt}\E\vert \Lambda_\varepsilon(j, x,\tilde{y})\vert^2\bar{\tilde{\mu}}^x(d\tilde{y})\nu(dx)=0
\end{align}
where 
\begin{align}\label{D34f}
\notag 2\E\vert \Lambda_\varepsilon(j,x,\tilde{y})\vert^2&= 2\E\left\vert \frac{1}{t(\varepsilon)}\int_0^{t(\varepsilon)} F\left( X_{jt(\varepsilon)}^{\varepsilon,x,\pi_2({\tilde{y}})},\Pi\tilde{Y}_{jt(t\varepsilon)+u}^{\varepsilon,x,\tilde{y}}\right)du-\bar{F}(X_{jt(\varepsilon)}^{\varepsilon,x,\pi_2(\tilde{y})})\right\vert^2\\
&\hspace{2cm}=\E\left[ \alpha\left(t(\varepsilon),X_{jt(\varepsilon)}^{\varepsilon,x,\pi(\tilde{y})}\right)\right].
\end{align} 
Integrating both sides of (\ref{F3Q1}) with respect to $\p(d\om)\bar{\tilde{\mu}}^x(d\tilde{y})\nu(dx)$ and taking note of (\ref{F3Q2}) and  (\ref{D34f}), we have 
\begin{align*}
\notag &\int_{\Nb}\int_{\Bt}\E\left(\sup_{0\leqslant t\leqslant T/\varepsilon}\Big\vert X_t^{\varepsilon,x,\pi_2(\tilde{y})}-\bar{X}_t^{\varepsilon,x}\Big\vert\right)\bar{\tilde{\mu}}^x(d\tilde{y})\nu(dx)\leqslant e^{T\tilde{C}_F}\bigg( \frac{2L_FT}{n(\varepsilon)}+\varepsilon T[\tilde{L}_F +\tilde{C}_FL_F]\\[.1cm] &\hspace{.5cm}
+T\max_{0\leqslant j\leqslant n(\varepsilon)-1}\int_{\Nb}\int_{\Bt}\left(\E\vert \Lambda_\varepsilon(j, x,\tilde{y})\vert^2\right)^{1/2}\bar{\mu}^x(d\tilde{y})\nu(dx)\bar{\tilde{\mu}}^x(d\tilde{y})\nu(dx)\bigg).
\end{align*}
Taking limit as $\varepsilon\rightarrow 0,$ and recalling that $n(\varepsilon)\rightarrow\infty$ as $\varepsilon\downarrow 0$ and in view of (\ref{D34e}), the proof is complete.
\qed

\begin{exa}\rm\;
We consider the following time periodic forcing SDE in two dimensions with $d=1$ and $N=1$.
\begin{equation}
\begin{cases}
dX= \varepsilon\left[ Y^2+\alpha X+\vartheta X^3\right]dt\\
dY = \left[-\gamma(X)Y+\beta\cos(2\pi t)\right]dt+ \sigma dW_t
\end{cases}
\end{equation}
where $0<\varepsilon\ll 1,\; \alpha, \vartheta\in \R, \; \beta\neq 0,\;  \sigma\neq 0$ and $\gamma:\R\rightarrow \Rp\setminus\{0\}$ is a polynomial. 
\end{exa}
\noindent {\it Proof.}
Here, the solution of fast subsystem $Y_t^{x,y}$ can be written explicitly as   
\begin{align*}
Y^{x,y}_t = e^{-\gamma(x)t}y+\beta\int_{0}^t e^{-\gamma(x)(t-s)}\cos(2\pi s)ds+ \sigma\int_{0}^t e^{-\gamma(x)(t-s)}dW_s,
\end{align*}
and the random periodic solution $S^x(t,\om)$ of period $\tau =1$ is given by 
\begin{align*}
S^x(t,\om) = \beta \int_{-\infty}^t e^{-\gamma(x)(t-s)}\cos(2\pi s)ds+  \sigma\int_{-\infty}^t e^{-\gamma(x)(t-s)}dW_s.
\end{align*}
The lifted random periodic solution $\tilde{S}^x(t,\om)$ on the cylinder $[0, 1)\times\R$ is 
\begin{align*}
\tilde{S}^x(t,\om) = (t\mod 1, S^x(t,\om)).
\end{align*}
The corresponding periodic measure $\tilde{\mu}^x_t$ is given by 
\begin{align*}
\tilde{\mu}^x_t(dyds) = \frac{1}{\sqrt{ 2\pi C(x)}}\exp\left(-\frac{(y-\beta v(t,x))^2}{2C(x)}\right)\I_{\{ t\mod 1\}}(s)dyds
\end{align*}
where 
\begin{align*}
v(t,x) = \int_{-\infty}^te^{-\gamma(x)(t-s)}\cos(2\pi s)ds
\end{align*}
and $C(x)$ solves the Lyapunov equation $2\gamma(x)C(x)=\sigma^2$.

Given the periodic measure $\tilde{\mu}^x_{t}$, we calculate he averaged vector field $\bar{F}(x)$  as follows 
\begin{align*}
\bar{F}(x) &= \int_{0}^1\int_{[0,1)\times\R}y^2\bar{\mu}^x_t(dyds)dt+\alpha x+\vartheta x^3\\
&= C(x)+\int_0^1\beta^2v^2(t,x)dt+\alpha x+\vartheta x^3.
\end{align*}
We compute the process $v(t,x)$ as 
\begin{align*}
v(t,x) = \frac{\gamma(x)\cos 2\pi t+ 2\pi\sin 2\pi t}{4\pi^2+\gamma^2(x)}
\end{align*}
and the integral of $v^2(t,x)$ over $[0, 1)$ is 
\begin{align*}
\int_0^1v^2(t,x)dt = \int_0^1\left(  \frac{\gamma(x)\cos 2\pi t+ 2\pi\sin 2\pi t}{\gamma^2(x)+4\pi^2}\right)^2 dt = \frac{2}{\gamma^2(x)+4\pi^2}.
\end{align*}
Finally, the rescaled slow motion $X_{t/\varepsilon}^{x,y}$ as $\varepsilon\downarrow 0$, can be approximated by the ODE 
\begin{align}
\frac{d\bar{X}_t}{dt} = \frac{\sigma^2}{2\gamma(\bar{X}_t)}+\frac{2\beta^2}{\gamma^2(\bar{X}_t)+4\pi^2}+\alpha\bar{X}_t+\vartheta \bar{X}^3_t.
\end{align}
\qed
\section{Conclusions and future work}
We proved a stochastic averaging principle for fully coupled SDEs with time periodic forcing. The ergodicity of evolution system of periodic probability measures on an irreducible Poincar\'e section is a key step to the identification and validity of the averaging limit. This key step is natural in applications like climate change studies and neural networks, for example, only the statistics and not a particular path of the weather (fast variables)  are important in climate evolution (e.g., \cite{Arnold01}). 

Averaging is one of the important steps for climatological applications, but it does not give the whole applicable recipes required in climate change studies. Some moderate and large deviation results for the slow motion from averaged motion are required for sensitive and rare events analysis. 
It is possible to use our framework to develop large deviation results and fluctuation dissipation theorems (FDTs) (e.g., \cite{Bran12, Andy10, Andy16}),  particularly important  in climate and neuronal models governed by multi-scale SDEs with time periodic forcing, where usual detailed balance condition (symmetry of stationary measures) fails. We shall investigate this aspect in the future publications and hope to carry out Hasselmann's program (e.g.,~\cite{Arnold01, Kifer04b, Kifer04}), i.e., approximation of slow subsystem by a nonlinear SDE where random fluctuations are taken into account, providing description for hopping of the slow variables (climate) between local basins of attractors/metastable states (extreme climate events).


\begin{thebibliography}{plain}
\setlength{\itemsep}{-1mm}
\bibitem{Arnold} L. Arnold, \emph{Random Dynamical Systems,} Springer, 1998.
\bibitem{Arnoldp} L. Arnold and M. Scheutzow, \emph{Perfect cocycles through stochastic differential equations,} Probability Theory and Related Fields Vol. 101 (1995), 65-88.
\bibitem{Arnold01} L. Arnold, \emph{Hasselmann's program revisted: the analysis of stochasticity in deterministic climate models}, Stochastic Climate Models: Progress in Probability, Vol. 49 (2001) , 141 -- 157, Birkh\"auser, Basel.
\bibitem{Arnaudon} M. Arnaudon and A. Thalmaier, \emph{The differentiation of hypoelliptic diffusion semigroups}, Illinois Journal of Mathematics, 54 (2010), 1285 -- 1311.
\bibitem{Kifer04b} V. Bakhtin and Y. Kifer, \emph{Diffusion approximation for slow motion in fully coupled averaging,} Probability Theory and Related Fields, Vol. 129 (2004), 157 -- 181.
\bibitem{Kry-Bog} N. N. Bogolyubov, Y. A. Mitropol'skii, \emph{Asymptotic methods in theory of nonlinear oscillations}, Gordon and Breach, Delhi, 1961. 
\bibitem{Bran12} M. Branicki and A. J. Majda, \emph{Quantifying uncertainty for predictions with model error in non-Gaussian systems with intermittency}, Nonlinearity, Vol. 25 (2012) 2543 -- 2578.
\bibitem{Bouchet15} F. Bouchet, T. Grafke, T. Tangarife and E. Vanden-Eijnden, \emph{Large deviations in fast-slow systems}, Journal of Statistical Physics, Vol. 162 (2016), 793 -- 812.
\bibitem{Cera17} S. Cerra and A. Lunardi, \emph{Averaging principle for nonautonomous slow-fast systems of stochastic reaction-diffusion equations: The almost periodic case}, SIAM Journal of Mathematical Analysis, Vol. 49 (2017), 2843 -- 2884.
\bibitem{Schu11} X. Chen, J. Duan and M. Scheutzow, \emph{Evolution systems of measures for stochastic flows}, Dynamical Systems, Vol. 26 (2011), 323 -- 334.
\bibitem{Checkron11} M. D. Checkroun, E. Simonnet and M. Ghil, \emph{Stochastic climate dynamics: Random attractors and time dependent invariant measures}, Physica D, Vol. 240 (2011), 1685 -- 1700.
\bibitem{Crau2} H. Crauel, \emph{Random Probability Measures on Polish spaces}, Taylor and Francis,  2002.
\bibitem{Daprato} G. Da Prato, \emph{Introduction to Stochastic Analysis and Malliavin Calculus}, Edizioni della Normale, 2008.
\bibitem{Elw78} K. D. Elworthy, \emph{Stochastic dynamical systems and their flows.} In: Stochastic Analysis, ed. A. Friedman, M. Pinsky, London-New York press, (1978) 79 - 95.
\bibitem{Feng11} C. Feng, H. Zhao and B. Zhou, \emph{Pathwise random periodic solutions of stochastic differential equations,} Journal of Differential Equations, 251(2011) 119 - 149.
\bibitem{Feng12} C. Feng and H. Zhao, \emph{Random Periodic Solutions of SPDEs via Integral Equations and Wiener-Sobolev Compact Embedding,} Journal of Functional Analysis, Vol. 262 (2012)  4377 - 4422.
\bibitem{Feng18} C. Feng and H. Zhao, \emph{Random periodic processes, periodic measures and ergodicity,} available at arXiv:1408.1897v4 [math.PR] 
\bibitem{Fredlinbook} M. I. Freidlin and A. D. Wentzell, \emph{Random perturbations of Dynamical Systems}, Springer, 1998.
\bibitem{Wain12} M. Galtier and G. Wainrib, \emph{Multiscale analysis of slow-fast neuronal learning models with noise}, The Journal of Mathematical Neuroscience, Vol. 2 (2012), 1 -- 64.
\bibitem{Hairer10} M. Hairer and A. J. Majda, \emph{A simple framework to justify linear response theory,} Nonlinearity 23 (2010), 909 - 922.
\bibitem{HairerD} M. Hairer, J. C. Mattingly and M. Scheutzow, \emph{Asymptotic coupling and a general form of Harris' theorem with applications to stochastic delay equation}, Probab. Theory Relat. Fields, 149 (2011), 223 - 259.
\bibitem{Hairer11a} M. Hairer and J. C. Mattingly, \emph{ A theory of hypoellipticity and unique ergodicity for semilinear stochastic PDEs}, Electronic Journal of Probability, Vol. 16 (2011), 658 - 738. 
\bibitem{Hairer_Note16b} M. Hairer, \emph{Advanced stochastic analysis}, Lecture Note Univ. Warwick. Available: http://www.hairer.org/Course.pdf.
\bibitem{Hairer_Note16} M. Hairer, \emph{Convergence of Markov processes}, Lecture Note Univ. Warwick. Available: http://www.hairer.org/notes/Convergence.pdf.
\bibitem{Hairer11} M. Hairer, \emph{On Malliavin's proof of H\"ormander's theorem}, Bull. Sci. Math. 135 (2011), 650 - 666.
\bibitem{Hairer11c} M. Hairer and N. S. Pillai, \emph{Ergodicity of hypoelliptic SDEs driven by fractional Brownian motion,} Annales de I'Institut Henri Poincar\'e - Prob. Stat. Vol. 47 (2011), 601 - 625.
\bibitem{Hasm80} R. Z. Has\'minskii, \emph{Stochastic stability of differential equations,} Sijthoff and Noordhoff, 1980 (Trnaslated from Russia).
\bibitem{Hasm66} R. Z. Has\'minskii, \emph{On the averaging principle for stochastic differential It\^o equations, } Kibernetika 4, (1968), 260--279.
\bibitem{Ikeda81} N. Ikeda and S. Watanabe, \emph{Stochastic Differential Equations and Diffusion Processes,} North Holland - Kodansha, Tokyo, 1981. 
\bibitem{Kifer01} Y. Kifer, \emph{Averaging and climate models}, Stochastic Climate Models: Progress in Probability, Vol. 49 (2001), 171 -- 188, Birkh\"auser, Basel.
\bibitem{Kifer09} Y. Kifer, \emph{Large deviations and adiabatic transitions for dynamical systems and Markov processes in fully coupled averaging}, American Mathematical Society, 2009.
\bibitem{Kifer04} Y. Kifer, \emph{Averaging principle for fully coupled dynamical systems and large deviations}, Ergodic Theory and Dynamical Systems, Vol. 24 (2004), 847 -- 871.
\bibitem{Kifer03} Y. Kifer, \emph{$L^2$ diffusion approximation for slow motion in averaging}, Stochastics and Dynamics, Vol. 3 (2003), 213 -- 246.
\bibitem{Kifer95} Y. Kifer, \emph{Limit theorems in averaging for dynamical systems}, Ergodic Theory and Dynamical Systems, Vol. 15 (1995), 1143 -- 1172.
\bibitem{Hkunita} H. Kunita, \emph{Stochastic differential equations and stochastic flow of diffeomorphisms,} \'Ecole d'\'et\'e de Probabilit\'es de Saint-Flour 12, 1982 Lecture Notes in Mathematics, Vol. 1097 (1984), 143 - 303. 
\bibitem{Kunita} H. Kunita, \emph{Stochastic flows and stochastic differential equations,} Cambridge University press, Cambridge, 1990.
\bibitem{Andy10}  A. J. Majda and X. Wang, \emph{Linear response theory for statistical ensembles in complex systems with time-periodic forcing}, Communication of Mathematical Sciences, Vol. 8 (2010), 145 -- 172.
\bibitem{Andy16} A. J. Majda, \emph{Introduction to Turbulent Dynamical Systems in Complex Systems}, Springer,  2016.
\bibitem{Malliavin} P. Malliavin, \emph{Stochastic analysis,} Grundlehren der Mathematischen Wissenschaften 313, Springer, 1997.
\bibitem{Mao94} X. Mao, \emph{Exponential stability of stochastic differential equations}, Marcel Dekker, New York, 1994.


\bibitem{Millet} A. Millet, D. Nualart and M. Sanz, \emph{Integration by parts and time reversal for diffusion processes,} The Annals of probability, Vol. 17 (1989), 208 - 238.
\bibitem{Loch88} P. Lochak and C. Meunier, \emph{Multiphase Averaging for Classical Systems with Application to Adiabatic Theorems}, Springer, 1988.  
\bibitem{Nulart} D. Nulart, \emph{The Malliavin calculus and related topics}, Springer, 2006.
\bibitem{Stuart08} G. A. Pavliotis and A. M. Stuart, \emph{Averaging and Homogenization}, Springer, 2008.
\bibitem{Schmal01} B. Schmalfuss, \emph{Lyapunov functions and non-trivial stationary solutions of stochastic differential equations,} Dynamical Systems, Vol. 16 N\underbar{o}: 8 (2001) 303 - 317.
\bibitem{Sturm00} R. Sturman and J. Stark, \emph{Semi-uniform ergodic theorems and applications to forced systems}, Nonlinearity, Vol. 13 (2000) 113 -- 143. 
\bibitem{Uda14} K. Uda, \emph{A qualitative approach to the existence of random periodic solutions}, Loughborough University Doctoral Thesis, 2015.
\bibitem{Uda16} K. Uda, \emph{Existence of random invariant periodic curves via random semiuniform ergodic theorem,} Stochastics and Dynamics, Vol. 17 (2017), 1750007.
\bibitem{Uda18} K. Uda and H. Zhao, \emph{Random periodic solutions and ergodicity for stochastic differential equations}. Submitted, 2018. Available at arXiv:1811.05503 [math.PR].
\bibitem{YU99}  A. Y. Veretennikov, \emph{On large deviations in the averaging principle for SDEs with full dependence}, The Annals of Probability, Vol. 27 (1999), 284 -- 296.
\bibitem{YU92} A. Y. Veretennikov, \emph{On large deviations in the avergaing principle for stochastic differential equations with periodic coefficients II}, Math. USSR Izvestiya, Vol. 39 (1992), 677.
\bibitem{Wain13} G. Wainrib, \emph{Double averaging for periodically forced slow-fast systems}, Electronic Communications in Probability, Vol. 18 (2013), 1 --12.
\bibitem{Wan151} B. Wang, \emph{Periodic and almost periodic random inertial manifolds for non-autonomous stochastic equations,}  page 189 - 208, Continuous and Distributed Systems II: Theory and Applications, Edited by V. A. Sadovnichiy and M. Z. Zgurovsky, Springer, 2015.
 \bibitem{Wan152} B. Wang, \emph{Stochastic bifurcation of pathwise almost random periodic and almost automorphic solutions for random dynamical systems,} Discrete and Continuous Dynamical Systems, Vol. 36, issue 1, (2015), 3745 - 3769.
 \bibitem{Wan14} B. Wang, \emph{ Existence, stability and bifurcation of random complete and periodic solutions of stochastic parabolic equations,} Nonlinear Analysis: Theory, Methods and Applications, Vol. 103 (2014), 9 - 25.
 \bibitem{Watanabe} S. Watanabe, \emph{Lectures on stochastic differential equations and Malliavin calculus}, Springer, 1984.
\bibitem{Zhao09} H. Zhao and Z. Zheng, \emph{Random Periodic Solutions of Random Dynamical systems,} Journal of differential equations, Vol. 246, N\underbar{o}: 5 (2009), 2020 - 2038. 
\end{thebibliography}
 \end{document}